\documentclass[11pt,final]{amsart}
\usepackage[utf8]{inputenc}
\usepackage[T1]{fontenc}
\usepackage[english]{babel}
\usepackage[letterpaper,centering]{geometry}
\usepackage{lmodern}
\usepackage{graphics}
\usepackage{mathrsfs}
\usepackage{amsmath,amssymb,amsfonts}
\usepackage{comment}

\usepackage{lineno}
\usepackage[textwidth=20mm, textsize=tiny]{todonotes}
\usepackage[unicode,pdfborder={0 0 0},final]{hyperref}
\usepackage[all]{xy}
{\setbox0\hbox{$ $}}\fontdimen16\textfont2=\fontdimen17\textfont2
\entrymodifiers={+!!<0pt,\the\fontdimen22\textfont2>}
\SelectTips{cm}{11}
\usepackage{quiver-5}
\mathcode`A="7041 \mathcode`B="7042 \mathcode`C="7043 \mathcode`D="7044
\mathcode`E="7045 \mathcode`F="7046 \mathcode`G="7047 \mathcode`H="7048
\mathcode`I="7049 \mathcode`J="704A \mathcode`K="704B \mathcode`L="704C
\mathcode`M="704D \mathcode`N="704E \mathcode`O="704F \mathcode`P="7050
\mathcode`Q="7051 \mathcode`R="7052 \mathcode`S="7053 \mathcode`T="7054
\mathcode`U="7055 \mathcode`V="7056 \mathcode`W="7057 \mathcode`X="7058
\mathcode`Y="7059 \mathcode`Z="705A

\theoremstyle{plain}
\newtheorem{thm}{Theorem}[section]

\newtheorem{lem}[thm]{Lemma}

\newtheorem{cor}[thm]{Corollary}
\newtheorem*{cor*}{Corollary}
\newtheorem{prop}[thm]{Proposition}

\theoremstyle{definition}

\newtheorem{defn}[thm]{Definition}
\newtheorem{rmk}[thm]{Remark}

\numberwithin{equation}{section}

\input cyracc.def
\DeclareFontFamily{U}{russian}{}
\DeclareFontShape{U}{russian}{m}{n}
        { <5><6> wncyr5
        <7><8><9> wncyr7
        <10><10.95><12><14.4><17.28><20.74><24.88> wncyr10 }{}
\DeclareSymbolFont{Russian}{U}{russian}{m}{n}
\DeclareSymbolFontAlphabet{\mathcyr}{Russian}
\makeatletter
\let\@math@cyr\mathcyr
\renewcommand{\mathcyr}[1]{\@math@cyr{\cyracc #1}}
\makeatother
\newcommand{\Sha}{{\mathcyr{Sh}}}
\newcommand{\ab}{{\mathrm{ab}}}

\newcommand{\A}{{\mathbf A}}
\newcommand{\C}{{\mathbf C}}
\newcommand{\F}{{\mathbf F}}

\renewcommand{\Im}{{\mathrm{Im}}}
\renewcommand{\P}{{\mathbf P}}
\newcommand{\Q}{{\mathbf Q}}
\newcommand{\R}{{\mathbf R}}
\newcommand{\Z}{{\mathbf Z}}
\newcommand{\CH}{\mathrm{CH}}
\newcommand{\nr}{\mathrm{nr}}
\newcommand{\Gm}{\mathbf{G}_\mathrm{m}}
\newcommand{\et}{{\text{ét}}}
\newcommand{\mot}{{\text{mot}}}

\newcommand{\SL}{\mathrm{SL}}
\newcommand{\Gal}{\mathrm{Gal}}
\newcommand{\Pic}{\mathrm{Pic}}
\newcommand{\Br}{\mathrm{Br}}
\newcommand{\Brnr}{\Br_\nr}
\newcommand{\Spec}{\mathrm{Spec}}

\newcommand{\rmH}{\mathrm{H}}
\renewcommand{\phi}{\varphi}
\renewcommand{\epsilon}{\varepsilon}
\renewcommand{\emptyset}{\varnothing}
\newcommand{\colim}{\rm colim}
\newcommand{\Ker}{{\mathrm{Ker}}}
\newcommand{\Coker}{{\mathrm{Coker}}}
\newcommand{\Hom}{{\mathrm{Hom}}}
\newcommand{\GL}{\mathrm{GL}}

\newcommand{\mmu}{\boldsymbol{\mu}}
\newcommand{\norm}{{\mathrm{nor}}}

\date{\today}
\title[Degree three unramified cohomology of homogeneous spaces]{Degree three unramified cohomology of homogeneous spaces with finite stabilisers}

\author{Lucas Lagarde}
\address{Institut Galil\'ee, Universit\'e Sorbonne Paris Nord, 99~avenue Jean-Baptiste Cl\'ement, 93430 Villetaneuse, France}
\curraddr{Institut für Algebraische Geometrie, Leibniz Universität Hannover, Welfengarten 1, 30167 Hannover, Germany}
\email{lagarde@math.uni-hannover.de}
\begin{document}
\begin{abstract} We provide an explicit formula for the third unramified cohomology group of homogeneous spaces of $\SL_n$ with finite stabilisers and a rational point, valid over an arbitrary field of characteristic prime to the order of the stabilisers. As a concrete application, we compute this group explicitly for every generalised quaternion group over an arbitrary field of characteristic 0, and settle Noether’s problem over the rational numbers in the negative for those of order at least 16, thereby generalising a classical result of Serre. Over an algebraically closed field, our formula recovers and sharpens a theorem of Peyre; we further show that universal triviality of unramified cohomology in degree at most three is equivalent to its vanishing in the same degrees. Over the function field of a $p$-adic curve, we extend a result of Harari–Scheiderer–Szamuely. Over a finite field, the formula is effective and yields an algorithm for computing the group explicitly.
\end{abstract}

\maketitle

\tableofcontents
\section{Introduction}
\label{sec:intro}
Let $G$ be a finite algebraic group over a field $k$ and let $G\to \GL(V)$ be a faithful $k$\nobreakdash-representation. A natural question is whether the field of invariant rational functions $k(V)^G$ is a purely transcendental extension over $k$. This problem was first studied by Emmy Noether, in the case where $k=\C$ is the field of complex numbers. It turns out in general that the field of invariant functions $k(V)^G$ coincides with the function field of a unirational variety defined over $k$. An obstruction for $k(V)^G$ to be purely transcendental lies in the so-called \emph{unramified cohomology groups} of this field relative to $k$: for purely transcendental extensions of $k$, these groups are reduced to Galois cohomology groups of the field $k$. Moreover, these groups are invariant under stable equivalence of field extensions of $k$. In degree two, unramified cohomology coincides with the unramified Brauer group. Artin and Mumford \cite{artinmumford} first used an incarnation of this group in order to disprove the rationality of certain unirational complex varieties. The unramified Brauer group was then used by Saltman to give the first example of a group $G$ for which the field $\C(V)^G$ is not stably rational \cite{saltmannoether}. Later, Bogomolov provided a general formula for this group \cite[Theorem 3.1]{bogomolovbrauer} over the complex numbers, which only depends on the cohomology of the group $G$. This formula allowed him to find other examples of groups $G$ for which the unramified Brauer group of $\C(V)^G$ is non-trivial. 

Unramified cohomology groups of arbitrary degree have then been introduced by Colliot-Thélène and Ojanguren \cite{CTO89} in order to give further examples of fields of unirational varieties over $\C$ that are not stably purely transcendental. However, these fields are not fields of invariant functions. It was only very much later that Peyre managed to find an explicit formula \cite[Theorem 3.1]{peyreh3noether} for the third unramified cohomology group of $\C(V)^G$ in terms of the cohomology of the group $G$. This enabled him to construct the first examples of groups $G$ such that $\C(V)^G$ is not stably rational but for which the unramified Brauer group is trivial. Kahn and Nguyen \cite[Théorème~5.10]{kahnnguyen} later provided a formula similar to that of Peyre over every field $k$ that contains all the roots of unity of order dividing the exponent of the group $G$.

Over fields of arithmetic interest, unramified cohomology is known in certain cases to provide obstructions to the existence of rational points or to their approximation properties. For instance, if $k$ is a global field and $X$ is a homogeneous space of a semi-simple simply connected group over $k$ (e.g.\ $\SL_{n,k}$) with finite geometric stabiliser, and if $k(X)$ denotes its function field, then the unramified Brauer group $\Brnr(k(X))$ can be used to define the Brauer-Manin obstruction to the Hasse principle on $X$, see \cite{wslc} for a complete exposition. Over the function field $k$ of a $p$-adic curve, Colliot-Thélène pointed out that one can make use of the third unramified cohomology group $\rmH^3_\nr(k(X),\Q/\Z(2))$ in order to define a reciprocity obstruction to the existence of rational points as well as their approximation properties, see \cite[§4]{hararischeidererszamuely}. (See also \cite{ctps} and \cite{wittenbergICM} for elaborations on this point.)

Going back to an arbitrary field $k$, in the specific case where $X$ already admits a $k$-point, the unramified cohomology groups of $X$ coincide with the unramified cohomology groups of $k(V)^G$ for any faithful representation $V$ of the stabiliser $G$ of this $k$-point. However, supplying practical formulas for the unramified cohomology groups of $k(V)^G$ is quite out of reach in full generality. Already for $\Brnr(k(V)^G)$, providing a concrete description of this group is notoriously much harder than over the complex numbers. Partial formulas for the unramified Brauer group of $k(V)^G$ (i.e.\ under specific conditions on the group $G$ or on the field $k$) have been extensively investigated in the modern literature, see e.g.\ \cite[Theorem~5.1]{borovoidemarcheharari}, \cite[Théorème 2.6]{ctnonramifie} and \cite[Theorem 5.1]{lucchiniunramifiedbrauer}. Explicit formulas for the entire unramified Brauer group have been provided very recently in \cite[Lemma~6.10]{loughransantens} and in \cite[Theorem 4.4]{lagardeunramifiedbrauer}. 

As for unramified cohomology in degree three, the situation is much more mysterious. For a torus $T$ over an arbitrary field $k$, or equivalently a homogeneous space with finite abelian geometric stabiliser (see the proof of \cite[Theorem 3.15]{LinhArithmetics}), Blinstein and Merkurjev \cite[Theorem B]{blinsteinmerkurjev} provided an exact sequence describing the group $\rmH^3_\nr(k(T),(\Q/\Z)'(2))$. However, one does not seem to be able to practically use this description to precisely exhibit elements in the unramified cohomology group. For a torus $T$ over a $p$-adic function field $k$, Harari, Scheiderer and Szamuely managed to relate the group $\rmH^3_\nr(k(T),(\Q/\Z)'(2))$ to a Tate-Shafarevich group $\Sha^2_\omega(k,T')$ which is more tractable \cite[Theorem 4.2]{hararischeidererszamuely}. More precisely, they showed that the image of this Tate-Shafarevich group in $\rmH^3_\nr(k(T),(\Q/\Z)'(2))$ controls the reciprocity obstruction to weak approximation on $T$. More recently, Linh showed \cite[Theorem B]{LinhArithmetics} that this image also determines the reciprocity obstruction to the existence of rational points on torsors of $T$. 

\subsection{Main results} A natural question posed in \cite{kahnnguyen} is whether one can find a general formula, for an arbitrary field $k$ of exponential characteristic $p>0$ and an arbitrary $k$\nobreakdash-group $G$ of order coprime to $p$, for the unramified cohomology group $\rmH^3_\nr(k(V)^G,(\Q/\Z)'(2))$ associated to any faithful $k$-representation $V$ of $G$. More generally, if $X$ is a homogeneous space of a semi-simple simply connected $k$-group with finite geometric stabiliser, one would like to find a general formula for the group $\rmH^3_\nr(k(X),(\Q/\Z)'(2))$. In this paper, we provide such a formula under the assumption that $X$ admits a $k$-point, see §\ref{introtheoremeprincipal} for more detail.

\subsubsection{Applications of the formula} Let us first mention significant consequences of this formula when specialised to various fields of arithmetico-geometric interest. 

If $G$ is a finite $k$-group, we denote by $\CH^2(B_kG)$ the second \emph{$G$-equivariant Chow group} of $k$ defined by Totaro (see §\ref{recollectionequivariant} for a recollection). This group is generated by Chern classes $c_2(\rho)$ of $k$-representations $\rho:G\to \GL(V)$, see~§\ref{rappelsKthéorie}. 
In \cite{garibaldimerkurjevserre}, Serre proved that the group $Q_{16}$ admits a non-trivial unramified and normalised cohomological invariant in degree~$3$ defined over $\Q$ (see §\ref{paragrapheexempleserre} for details), thereby providing a counter-example to Noether's problem over the field of rational numbers. Using the same cohomological invariant, Hoshi, Kang and Yamasaki \cite{hoshikangyamasaki} extended Serre’s result to real quadratic fields $\Q(\sqrt{d})$, where $d$ is a square-free integer of the form $1+8m$ with $m\ge2$. We apply our formula to explicitly determine the corresponding unramified cohomology class; specifically, we show that this class arises from a $\Gamma_\Q$-invariant element of $\CH^2(B_{\overline{\Q}}Q_{16})$. Furthermore, we establish that this element generates the entire unramified cohomology group modulo the constant classes. 

If we instead consider generalised quaternion groups $Q_{2^n}$ with $n\ge 5$, then Noether's problem over $\Q$ remained open, see e.g.\ \cite{hukangpgroups}. A result of Demarche \cite[Corollaire 3, Remarque 7]{demarchebrnr} asserts that for all such groups, the unramified Brauer group $\Brnr(\Q(V)^{Q_{2^n}})$ is reduced to the constants and hence cannot obstruct stable rationality. Moreover, Serre’s argument, which relies on an embedding $Q_{16}\hookrightarrow\mathrm{Spin}_{7,\Q}$ and on the theory of trace forms in dimension at most $7$, does not extend to $Q_{2^n}$ for $n\ge5$; see §\ref{paragrapheexempleserre2}. Our formula bypasses this issue: the generator is realised by pulling back the Rost invariant of a Weil restriction of $\SL_1(A)$ for a suitable quaternion algebra $A$, cf.\ §\ref{invariantdeRost}. It turns out that our computation for $Q_{16}$ can be conducted for $Q_{2^n}$ for any $n\ge 4$. Consequently, we settle Noether's problem over $\Q$ in the negative for all such groups:

\begin{thm}[Theorem \ref{Q_16}] Let $n\ge 4$ and let $Q_{2^n}\to \GL(V)$ be a faithful $\Q$\nobreakdash-representation. There is a split short exact sequence
    \[0\longrightarrow \rmH^3(\Q,\Q/\Z(2))\longrightarrow \rmH^3_\nr(\Q(V)^{Q_{2^n}},\Q/\Z(2)) \longrightarrow \Z/2\longrightarrow 0.\]

Notably, Noether's problem for the field of invariants $\Q(V)^{Q_{2^n}}$ admits a negative answer. More precisely, the quotient of $\rmH^3_\nr(\Q(V)^{Q_{2^n}},\Q/\Z(2))$ by the subgroup of constant classes identifies with the cokernel of the base change map $\CH^2(B_{\Q}Q_{2^n})\to \CH^2(B_{\overline{\Q}}Q_{2^n})^{\Gamma_\Q}$. Moreover, it is generated by the image of $2^{n-3}c_2(\rho)$ where $\rho$ denotes any irreducible faithful $\overline{\Q}$-representation $Q_{2^n}\to \SL_2(\overline{\Q})$. 
\end{thm}

More generally, over any field $K$ of characteristic $0$, if we let $E:= K(\zeta_{2^{n-1}}+\zeta_{2^{n-1}}^{-1})$, then the quotient of $ \rmH^3_\nr(K(V)^{Q_{2^n}},\Q/\Z(2))$ by the subgroup of constant classes is isomorphic to $\Z/2$ if $K\neq E$ and the class of the quaternion algebra $(-1,-1)_E$ is non-zero in $\Br(E)$. The aforementioned quotient vanishes otherwise; see Theorem \ref{Q_16} for a more precise statement. To our knowledge, Theorem \ref{Q_16} is the first complete and explicit computation of the third unramified cohomology group of linear quotients for an infinite family of classical non\nobreakdash-abelian finite groups, over an arbitrary field of characteristic $0$, in particular over a number field. Theorem \ref{Q_16} also recovers \cite[Theorem 1.4]{hoshikangyamasaki} for $Q_{16}$, since their assumptions imply that $\sqrt2\notin K$ and $(-1,-1)_{K(\sqrt{2})}\neq0$, and extends it: for instance, for every field $K$ of characteristic $0$ which admits an ordering and does not contain $\sqrt{2}$, the field $K(V)^{Q_{16}}$ is not retract rational over $K$ for any faithful $K$-representation $V$, without any congruence condition.

For smooth and projective varieties $X$ over an algebraically closed field $k$, a question of fundamental interest is whether the unramified cohomology group $\rmH^3_\nr(k(X),(\Q/\Z)'(2))$ is \emph{universally trivial} as soon as it is trivial. That is, for any field $L$ that contains $k$, we ask whether the natural map $\rmH^3(L,(\Q/\Z)'(2)) \to \rmH^3_\nr(L(X),(\Q/\Z)'(2))$ is an isomorphism. This property provides a finer obstruction to rationality than mere triviality over the base field; indeed, Voisin \cite{voisindegen, voisinuniversel} demonstrated the existence of complex varieties whose unramified cohomology in degree three is trivial without being universally trivial. 

On the other hand, it is known \cite{auelcolliotparimala} that cubic fourfolds containing a plane have universally trivial unramified cohomology in degrees at most $3$. It is therefore natural to ask whether the same statement holds for homogeneous spaces with finite stabilisers, that is, whether unramified cohomology is universally trivial as soon as it is trivial. We give a complete answer in degrees at most $3$.

If $G$ is an abstract finite group, we denote by $B_0(G)$ its \emph{Bogomolov multiplier}, that is, the subgroup of $\rmH^2(G,\Q/\Z)$ consisting of those classes which vanish along the restriction map $\rmH^2(G,\Q/\Z)\to \rmH^2(A,\Q/\Z)$ for any abelian subgroup $A\subset G$. 
\begin{thm}[Theorem \ref{trivialitépartiealgébrique}, Corollary \ref{trivialitéuniverselle}] Let $k$ be an algebraically closed field of exponential characteristic $p>0$, let $G$ be a finite group of order coprime to $p$ and let $G\to \GL(V)$ be a faithful $k$-representation. If the group $B_0(G)$ is trivial, then for any field $L$ containing $k$, there exists a canonical exact sequence:
\[0\longrightarrow\rmH^3(L,(\Q/\Z)'(2))\longrightarrow\rmH^3_\nr(L(V)^G,(\Q/\Z)'(2))\longrightarrow \rmH^3_\nr(k(V)^G,(\Q/\Z)'(2)).\]

If we further assume that the third unramified cohomology group $\rmH^3_\nr(k(V)^G,(\Q/\Z)'(2))$ is trivial, then it is universally trivial, i.e., for any field $L$ containing $k$, the natural map 
\[\rmH^3(L,(\Q/\Z)'(2))\longrightarrow \rmH^3_\nr(L(V)^G,(\Q/\Z)'(2))\] 
is an isomorphism.
\end{thm}
To our knowledge, such a result is entirely new if we exclude the specific case where $G$ is abelian, see \cite[Corollary 5.8]{blinsteinmerkurjev} and \cite[Remarque 5.4]{ctnonramifie}. Notably, the formula provided in \cite{kahnnguyen} does not seem to allow one to recover our result. We would like to draw the attention of the reader to the fact that the assumption on the vanishing of $B_0(G)$ is necessary for the conclusion to hold. Over the field $L=k((t))$, any non-trivial class in $B_0(G)$ produces a non-zero element of $\rmH^3_\nr(L(V)^G,(\Q/\Z)'(2))$; see Remark~\ref{remarquebogomolov} for more detail.

Over a finite field $k$, the computation of the third unramified cohomology group for various classes of smooth and projective $k$-varieties is a delicate and open problem of fundamental importance. Its significance stems again from its applications to stable rationality, but also its connections to questions regarding algebraic cycles in codimension~$2$; these topics are extensively discussed in the literature, for instance in \cite{colliotkahn}, \cite{colliotscavia}, \cite{colliotvoisin} and \cite{pirutka}. 

However, the existence of an algorithm capable of effectively computing unramified cohomology in degree $3$ for broad classes of varieties remains, in general, out of reach. We provide such an algorithm in the case of homogenous spaces with finite stabilisers over a finite field:

\begin{thm}[Theorem \ref{algorithmecorpsfinis}] There exists an algorithm that, given the datum of a finite field $k$ of characteristic $p$ and a finite $k$-group $G$ of odd order coprime to $p$, computes the group $\rmH^3_\nr(k(V)^G,(\Q/\Z)'(2))$ for any faithful $k$-representation $G\to \GL(V)$ (this does not depend on the prescribed representation).
\end{thm}

In general, the kernel of the natural restriction $\rmH^3_\et(X,(\Q/\Z)'(2))\to \rmH^3_\et(\overline{X},(\Q/\Z)'(2))$ may contribute non-trivially to the unramified cohomology group $\rmH^3_\nr(k(X),(\Q/\Z)'(2))$, e.g.\ it does as soon as the restriction map $\CH^2(X)\to \CH^2(\overline{X})$ is injective; see Theorem~\ref{theoremeprincipal}. Over the function field of a curve defined over a $p$-adic field, such as $\Q_p(t)$, there can be uncountable many non-trivial classes in the unramified cohomology group that arise from the aforementioned kernel and are non-trivial modulo the constants, already for abelian stabilisers:
\begin{thm}[Theorem \ref{contrexempleZ/8}]
    Let $k=\Q_2(t)$ and $G=\Z/8$, and fix an embedding of $k$-groups $G\hookrightarrow \SL_{n,k}$ for some $n\ge 1$. If we denote by $X$ the $k$-variety $\SL_{n,k}/G$, then in the quotient of $\rmH^3_\nr(k(X),\Q/\Z(2))$ by the subgroup $\rmH^3(k,\Q/\Z(2))$ of constant classes, there exist uncountably many non-zero elements that lie in the image of the natural map 
    \[\Ker\Big(\rmH^3_\et(X,\Q/\Z(2))\to \rmH^3_\et(\overline{X},\Q/\Z(2))\Big)\longrightarrow\frac{\rmH^3(k(X),\Q/\Z(2))}{\rmH^3(k,\Q/\Z(2))}.\]
    Moreover, each one of these classes provides a reciprocity obstruction to weak approximation on the $k$-variety $X$.
\end{thm}

The subgroup of $\rmH^3_\et(X,\Q/\Z(2))$ containing the classes from the above result coincides with the Tate\nobreakdash-Shafarevich group $\Sha^2_\omega(k,\mmu_8^{\otimes2})$ defined in \cite{hararischeidererszamuely} and further studied in \cite{harariszamuely}; see Corollary \ref{comparaisonsha2}. The above theorem notably refines and improves Corollary 4.2 from \emph{op.\ cit.}\ 

In \S\ref{paragraphehomogène}, for an arbitrary field $k$ and any finite $k$-group $G$, we define a "Tate-Shafarevich" group $\Sha^2_\nr(k,\widehat{G}^\ab(1))$ which provides a natural source of classes in the unramified cohomology group $\rmH^3_\nr(k(V)^G,(\Q/\Z)'(2))$ for any faithful $k$-representation $G\to \GL(V)$. This construction generalises previous results in several directions. First, when $k$ is the function field of a $p$-adic curve and $G$ is abelian, this group coincides with the aforementioned Tate-Shafarevich group introduced in \cite{hararischeidererszamuely} (see Theorem~\ref{comparaisonsha2}). Furthermore, we prove that this group vanishes whenever $G$ is constant and $k$ contains all the roots of unity dividing the exponent of $G$; this generalizes \cite[Lemma 3.2]{hararischeidererszamuely}, whose proof relies on Lichtenbaum duality for $p$-adic curves. Finally, these vanishing properties demonstrate that this specific subgroup cannot be described using the results of \cite{kahnnguyen} as their main theorem assumes that $G$ is constant and the base field contains all the roots of unity dividing the exponent of $G$.

\subsubsection{A general formula}\label{introtheoremeprincipal} All of the aforementioned results crucially rely on the general formula that we obtain for the unramified cohomology group. Let us succinctly describe it. For the sake of exposition, we assume here that $k$ has characteristic $0$. Let $G$ be a finite $k$-group, let $G\to \GL(V)$ be a faithful $k$-representation and denote by $G(\overline{k})\rtimes\Gamma_k$ the semi-direct product induced by the prescribed action of $\Gamma_k$ on $G(\overline{k})$. In §\ref{paragraphesuitefondamentale}, we show that there is a natural inflation map 
\begin{equation}\label{inflation}\rmH^3(G(\overline{k})\rtimes\Gamma_k,\Q/\Z(2))\to \rmH^3(k(V)^G,\Q/\Z(2)),\end{equation}
whose image contains $\rmH^3_\nr(k(V)^G,\Q/\Z(2))$. In §\ref{definitionglobalementnr}, we provide a construction of residue maps $\partial_{\mu,D}:\rmH^3(G(\overline{k})\rtimes\Gamma_k,\Q/\Z(2))\to \rmH^2(D,\Q/\Z(1))$ for adequate morphisms $\mu:\widehat{\Z}(1)\to G(\overline{k})$ and subgroups $D\subset G(\overline{k})\rtimes\Gamma_k$ satisfying a semi-commutation relation. In the case where the field $k$ contains all the roots of unity of order dividing the exponent of $G$, these residues agree with those defined by Kahn and Nguyen in \cite{kahnnguyen}, as well as those defined by Peyre in \cite{peyreh3noether}, see Remark \ref{remarquerésidus}. We then denote by $\rmH^3_\mathrm{gnr}(G(\overline{k})\rtimes\Gamma_k,\Q/\Z(2))$ the subgroup of \emph{globally unramified classes} of $\rmH^3(G(\overline{k})\rtimes\Gamma_k,\Q/\Z(2)))$ consisting of those classes which die along all of these residue maps. 

In §\ref{paragraphepermutationnégligeable}, we define a group $\rmH^3_\mathrm{per}(G(\overline{k})\rtimes\Gamma_k,\Q/\Z(2))$ which is contained in the kernel of the inflation map~\eqref{inflation}. The definition of this group expands the initial definition, due to Saltman, of \emph{permutation negligible classes} in group cohomology; see \cite{saltmanpermutation} for reference. The following result is a direct generalisation of a theorem of Peyre \cite[Theorem 1]{peyreh3noether}, who treated the case where $k$ is an algebraically closed field of characteristic $0$.

\begin{thm}[Theorem \ref{theoremeprincipal}] The inflation map \eqref{inflation} induces a short exact sequence 
\[0\longrightarrow \CH^2(B_kG)\longrightarrow \rmH^3_\mathrm{gnr}(G(\overline{k})\rtimes\Gamma_k,(\Q/\Z)'(2))\longrightarrow \rmH^3_\nr(k(V)^G,(\Q/\Z)'(2))\longrightarrow 0,\]
If we further assume that $G$ is constant, then we have a surjection:
\[\frac{\rmH^3_\mathrm{gnr}(G\times\Gamma_k,(\Q/\Z)'(2))}{\rmH^3_\mathrm{per}(G\times \Gamma_k,(\Q/\Z)'(2))}\twoheadrightarrow \rmH^3_\nr(k(V)^G,(\Q/\Z)'(2)),\]
the kernel of which is killed by a power of $2$. In particular if $G$ has odd order, then the above map is an isomorphism. 
\end{thm}

\subsection{Notation and conventions} All throughout this text, if $k$ is a field and $\overline{k}$ is a fixed algebraic closure, we denote by $\Gamma_k$ its absolute Galois group $\Gal(\overline{k}/k)$. If $X$ is a $k$-scheme and $n\ge 1$ is an integer that is invertible on $X$, we denote by $\mmu_n$ the étale sheaf of $n$\textsuperscript{th} roots of unity on $X$. For $j\ge 1$, we let $\mmu_n^{\otimes j}$ denote $\mmu_n^{\otimes (j-1)}\otimes \mmu_n$. We also let  $\mmu_n^{\otimes 0}:=\Z/n$ and, for $j<0$, $\mmu_n^{\otimes j}:=\Hom(\mmu_n^{\otimes (-j)},\Z/n)$. We write $\Q/\Z(j):={\colim}_{n\ge 1}\hspace{0.05cm} \mmu_n^{\otimes j}$, resp.\ $\widehat{\Z}(j):=\lim_{n\ge 1}\hspace{0.05cm} \mmu_n^{\otimes j}$. Similarly, if $p>0$ denotes the exponential characteristic of $k$, we write $(\Q/\Z)'(j):={\colim}_{\gcd(n,p)=1}\hspace{0.05cm} \mmu_n^{\otimes j}$, resp.\ $\widehat{\Z}'(j):=\lim_{\gcd(n,p)=1}\hspace{0.05cm} \mmu_n^{\otimes j}$.

Given an (abstract) finite group $G$, we denote by $|G|$ its order and $\mathrm{exp}(G)$ its exponent and, if $G$ is abelian and $n\ge 1$ is an integer (resp.\ $p$ is a prime number), we denote by $G[n]$ (resp.\ $G\{p\}$) its $n$-torsion (resp.\ $p$-primary torsion) subgroup. We denote by $G\{p'\}$ the subgroup of torsion coprime to $p$. If $G$ is an arbitrary profinite group and $p$ is a prime number, we denote by $G^{(p')}$ its pro-$p'$ completion. If $H\subset G$ is a subgroup, we denote by $N_G(H)$ (resp.\ $C_G(H)$) its normaliser (resp.\ centraliser) in $G$. 

If $G$ is a finite $k$-group, we will make the following abuse of notation: we will always identify the finite algebraic group with the abstract group of its $\overline{k}$-points $G(\overline{k})$ seen as a $\Gamma_k$-group, i.e.\ as an abstract group equipped with a continuous action of $\Gamma_k$. A finite $k$-group $G$ is said to be constant if $\Gamma_k$ acts trivially on $G(\overline{k})$. For each finite $k$-group $G$ of order coprime to the exponential characteristic of $k$, we denote by $\widehat{G}$ its Cartier dual, i.e.\ the group $\Hom_{\overline{k}}(G(\overline{k}),\overline{k}^*)\simeq \Hom_{\overline{k}}(G(\overline{k}),(\Q/\Z)'(1))$ (morphisms taken in the category of $\overline{k}$-groups) equipped with the continuous action $(\sigma\cdot f)(g):=\sigma\cdot f(\sigma^{-1}\cdot g)$. When $k$ is algebraically closed, this group is (non-canonically) isomorphic to the Pontrjagin dual $\Hom(G,\Q/\Z)$. We let $\Hom_k(G(\overline{k}),(\Q/\Z)'(1))$ be the subgroup of $\Gamma_k$-invariants within $\Hom_{\overline{k}}(G(\overline{k}),(\Q/\Z)'(1))$.

A \emph{variety} over a field $k$ is a separated scheme of finite type over $k$. If $X$ is a variety, we denote by $\overline{X}$ the base change $X\otimes_k \overline{k}$ to the algebraic closure. When $X$ is integral, we denote its function field by $k(X)$. If $\mathscr{F}$ is an étale sheaf on $X$, we denote by $\rmH^*_\et(X,\mathscr{F})$ the corresponding étale cohomology groups. If $x$ is a fixed $k$-point of $X$, then unless explicitly mentioned, we denote by $\rmH^{*,\norm}_\et(X,\mmu_n^{\otimes i})$ the kernel of the pullback map $x^*:\rmH^*_\et(X,\mmu_n^{\otimes i})\to \rmH^*(k,\mmu_n^{\otimes i})$. If $\ell$ is a prime number distinct from the characteristic of $k$, we let $\rmH^*_\et(X,\Z_\ell(i))$ be the limit, over $n\ge 1$, of the groups $\rmH^*_\et(X,\mmu_{\ell^n}^{\otimes i})$. We denote Steenrod operations in étale cohomology in the sense of \cite{voevodskyreducedpowers} as $\mathrm{Sq}^*(-)$.

If $\Lambda$ is a commutative ring and $i,j\in\Z$, we denote by $\rmH^i_\mot(X,\Lambda(j))$ the $\mathrm{cdh}$-motivic cohomology group $\rmH^i(Y_\mathrm{cdh},\Lambda(j))$, where $\Lambda(j):=C_*\Lambda_{\mathrm{tr}}(\Gm^{\wedge j})[-j]$ is the bounded above complex of $\Lambda$- ules defined in \cite{suslinvoevodsky}. When $X$ is smooth, we may replace the cdh topology by the Nisnevich topology. If $\Gamma$ is a profinite group and $M$ is a continuous $\Gamma$-module, we denote by $\rmH^*(\Gamma,M)$ the associated cohomology groups. If $R$ is a commutative ring, we denote by $K_*^M(R)$ its Milnor $K$-theory groups. 

\subsection{Acknowledgements} This work was written as part of my PhD thesis, done under the supervision of Olivier Wittenberg. I would like to express my sincerest gratitude to him for his support and guidance, as well as for his many suggestions and improvements on this text. I am grateful to Jean-Louis Colliot-Thélène, Pierre Guillot, David Harari, Anne Quéguiner-Mathieu, Federico Scavia and Tamás Szamuely for interesting discussions. I would also like to thank Charles De Clercq for raising the question that led to the example discussed in Theorem~\ref{contrexempleZ/8}. 

This work was carried out at Université Sorbonne Paris Nord and completed at Leibniz Universität Hannover.

\section{Preliminaries}

\subsection{Galois covers of classifying schemes}\label{galoiscoversclassifyingschemes}\label{revêtementsétales}
Let $k$ be a field of exponential characteristic $p>0$ and let $G$ be a finite $k$-group of order coprime to $p$. Suppose that $G$ acts linearly and faithfully on $\A^m_k$, let $c\ge 1$ be an integer and let $U\subset \A^m_k$ be an open subscheme on which $G$ acts freely. Let $Z\subset \A^m_k$ be the closed complement of $U$, and assume that $\mathrm{codim}_{\A^m_k}(Z)\ge c+1$. The geometric quotient $X:=U/G$ exists as a smooth and irreducible $k$-variety, and the morphism $U\to X$ is a $G$-torsor. Such a pair $(U\subset \A^m_k,c)$ always exists \cite[Remark~2.7]{totaro}, and we will refer to it as a \emph{classifying pair}. Fix a $k$-point $u$ of $U$; this yields a $k$-point $x$ of $X$, which we view as a geometric point $\overline{x}\mapsto x\in X$. By Zariski-Nagata purity \cite[Exposé~X, Théorème~3.1]{SGA1}, we have that $\pi_1^\et(\overline{U},\overline{u})=\pi_1^\et(\A^m_{\overline{k}},\overline{u})$; moreover, the group $\pi_1^\et(\A^m_{\overline{k}},\overline{u})^{(p')}$ is trivial by virtue of e.g.\ \cite[Exposé X, Théorème 3.10]{SGA2} and \cite[Exposé XIII, Corollaire 2.12]{SGA1}. 

Let $Y\to X$ be an étale cover and $\overline{y}$ a given geometric point. Define $K\subset \overline{k}$ as the subfield fixed by the image of $\pi_1^\et(Y,\overline{y})\to \Gamma_k$. It follows that $Y$ naturally carries the structure of a $K$-scheme; we call $K$ its \emph{field of definition}.
\begin{defn}
A \emph{good} cover $\pi:Y\to X$ is an étale cover such that, if $K$ denotes the field of definition of $Y$, then the Galois closure over $\overline{X}$ of any irreducible component of $Y\otimes_K \overline{k}$ is of degree coprime to $p$ over $\overline{X}$.
\end{defn}
If $k$ notably has characteristic $0$, then any étale cover of $X$ is good. Let $\pi:Y\to X$ be an irreducible good cover and $K$ the field of definition of $Y$. Then $Y\otimes_K \overline{k}\simeq \overline{U}/H$ for some abstract subgroup $H\subset G(\overline{k})$, so we may consider a geometric point $\overline{y}\mapsto y\in Y$ induced by a point of $\overline{U}$. We hence have  $\pi_1^\et(Y\otimes_K \overline{k},\overline{y})^{(p')}\simeq H$ as abstract groups. By virtue of the fundamental exact sequence \cite[Exposé V]{SGA1}
\begin{equation}\label{suitegroupesfondamentaux}1\longrightarrow \pi_1^\et(Y\otimes_K \overline{k},\overline{y})\longrightarrow \pi_1^\text{ét}(Y,\overline{y})\longrightarrow\Gamma_K\longrightarrow 1,\end{equation}
there is an induced outer action of $\Gamma_K$ on $\pi_1^\et(Y\otimes_K \overline{k},\overline{y})$ and in particular on $H$. The isomorphism class of $H$ seen as a group equipped with an outer action of $\Gamma_K$ does not depend on the choice of $\overline{y}$. We define the group $\pi_1^\et(Y,\overline{y})^{[p']}$ as the quotient of $\pi_1^\et(Y,\overline{y})$ by the kernel of $\pi_1^\et(Y\otimes_K \overline{k},\overline{y})\to \pi_1^\et(Y\otimes_K \overline{k},\overline{y})^{(p')}$. 

If $y$ is a $K$-point, then it induces a continuous section of the projection $\pi_1^\text{ét}(Y,\overline{y})\to \Gamma_K$, and hence also a splitting of the short exact sequence 
\[1\longrightarrow H\longrightarrow \pi_1^\text{ét}(Y,\overline{y})^{[p']}\longrightarrow \Gamma_K\longrightarrow 1.\] 
Thus if $Y(K)\neq \emptyset$ then the extension \eqref{suitegroupesfondamentaux} splits and we may write $\pi_1^\text{ét}(Y,\overline{y})\simeq \pi_1^\et(Y\otimes_K \overline{k},\overline{y})\rtimes\Gamma_K$ (resp.\ $\pi_1^\text{ét}(Y,\overline{y})^{[p']}\simeq H\rtimes\Gamma_K$), where the isomorphism depends on the choice of a $K$-point of $Y$.


\subsection{Equivariant Chow groups and étale cohomology}\label{recollectionequivariant} Let $k$ be a field of exponential characteristic $p> 0$ and let $G$ be a finite $k$-group of order coprime to $p$. Let $(U\subset \A^m_k,c)$ be a classifying pair. Following Totaro \cite[Remark~2.7]{totaro}, we define the \emph{equivariant Chow groups} $\CH^r(B_kG):=\CH^r(U/G)$ for $r\le c$ (these may be interpreted as the Chow groups of the classifying stack $B_kG$ of $G$ over $k$, see \emph{loc.\ cit.}, §2.2); this definition does not depend on the choice of $U\subset \A^m_k$ \cite[Theorem~2.5]{totaro}. 

Similarly, we might define the equivariant cohomology of $B_kG$ (or more generally twisted forms thereof), with torsion coefficients, as the étale cohomology of $U/G$ (or of its good covers); this boils down to profinite group cohomology, as the following lemma shows:


\begin{lem}\label{cohomologieequivariante}
    Let $k$ be a field of exponential characteristic $p>0$, let $G$ be a finite $k$-group of order coprime to $p$, let $c\ge 0$ and fix a classifying pair $(U\subset \A^m_k,c)$. Let $n\ge 1$ be an integer which is coprime to $p$. If we denote $X:=U/G,$ then for any irreducible good cover $\pi:Y\to X$, for any $j\in \Z$ and for any $r\le 2c-2$, we have a natural isomorphism:
    \[\rmH^r_\et(Y,\mmu_n^{\otimes j})\simeq\rmH^r(\pi_1^\et(Y,\overline{y})^{[p']},\mmu_n^{\otimes j}),\]
    where $\pi_1^\et(Y,\overline{y})^{[p']}$ acts on $\mmu_n^{\otimes j}$ via its image in $\Gamma_k$.
\end{lem}
\begin{proof} Let $K$ be the field of definition of $Y$, fix a choice of an isomorphism $Y\otimes_K \overline{k}\simeq \overline{U}/H$ and consider the induced Galois cover $\pi':\overline{U}\to Y$. Since $\pi_1^\et(\overline{U},\overline{u})^{(p')}=0$ and $H$ has order coprime to $p$, it follows that the automorphism group of $\pi'$ coincides with $\pi_1^\et(Y,\overline{y})^{[p']}$. The Hochschild-Serre spectral sequence associated to $\pi'$ reads:
\[E_2^{p,q}=\rmH^p(\pi_1^\et(Y,\overline{y})^{[p']},\rmH^q_\et(\overline{U},\mmu_n^{\otimes j}))\Rightarrow \rmH^{p+q}_\et(Y,\mmu_n^{\otimes j}).\]
It therefore suffices to show that $\rmH^q_\et(\overline{U},\mmu_n^{\otimes j})$ is zero for $1\le q \le 2c-2$. The long exact localisation sequence in étale cohomology for the pair $(\overline{U}\subset \A^m_{\overline{k}})$ yields exact portions:
\[\rmH^q_{\et,\A^m_{\overline{k}}\setminus \overline{U}}(\A^m_{\overline{k}},\mmu_n^{\otimes j})\longrightarrow \rmH^q_\et(\A^m_{\overline{k}},\mmu_n^{\otimes j})\longrightarrow \rmH^q_\et(\overline{U},\mmu_n^{\otimes j}) \longrightarrow \rmH^{q+1}_{\et,\A^m_{\overline{k}}\setminus \overline{U}}(\A^m_{\overline{k}},\mmu_n^{\otimes j}).\]
As $n$ is coprime to $p$, we have by semi-purity \cite[§8]{fujiwara} that $\rmH^{r}_{\et,\A^m_{\overline{k}}\setminus \overline{U}}(\A^m_{\overline{k}},\mmu_n^{\otimes j})=0$ for any $r\le 2c-1$. In particular, we get that $\rmH^r_\et(\overline{U},\mmu_n^{\otimes j})\simeq\rmH^r_\et(\A^m_{\overline{k}},\mmu_n^{\otimes j})$, and the latter is isomorphic to $\rmH^r(\overline{k},\mmu_n^{\otimes j})$ thanks to the homotopy invariance of étale cohomology \cite[Chapter VI, Corollary 4.20]{milneetale}. It follows that for $1\le q \le 2c-2$, we have $\rmH^q_\et(\overline{U},\mmu_n^{\otimes j})=0$, and hence $E_2^{p,q}=0$ for $1\le q\le 2c-2$; we thus derive a map $E_2^{r,0}\to \rmH^r_\et(Y,\mmu_n^{\otimes j})$ which is an isomorphism for $1\le r\le 2c-2$. 

Finally, since $\overline{U}$ is an open subset of $\A^m_{\overline{k}}$, we have the equality $\rmH^0(\overline{U},\mmu_n^{\otimes j})=\mmu_n^{\otimes j}$, so that we may identify $E_2^{r,0}=\rmH^r(\pi_1^\et(Y,\overline{y})^{[p']},\mmu_n^{\otimes j})$ where $\pi_1^\et(Y,\overline{y})^{[p']}$ acts on $\mmu_n^{\otimes j}$ through the cyclotomic action of its image in $\Gamma_k$. This concludes the proof of the lemma.
\end{proof}

\subsection{Bloch-Ogus theory} Let $V$ be a smooth and connected variety over $k$ and let $n\ge 1$ be an integer that is invertible in $k$. Let $\xi:V_\et\to V_\text{Zar}$ be the change-of-sites map. For $j\in \Z$ and $q\ge 0$, consider the Zariski sheaf $\mathscr{H}^q_\et(\mmu_n^{\otimes j}):=R^q\xi_*\mmu_n^{\otimes j}$ on $V$. This coincides with the sheafification of the presheaf that sends an open subset $U\subset V$ to $\rmH^q_\et(U,\mmu_n^{\otimes j})$.

Let $k\subset R\subset k(V)$ be a discrete valuation ring whose fraction field is $k(V)$ and denote by $\kappa$ its residue field. One can define a residue map:
\[\partial_R:\rmH^q(k(V),\mmu_n^{\otimes j})\to \rmH^{q-1}(\kappa,\mmu_n^{\otimes (j-1)}),\]
see \cite[§1]{CTO89}.
We define the \emph{unramified cohomology} group $\rmH^q_\mathrm{nr}(k(V),\mmu_n^{\otimes j})$ as the subgroup of $\rmH^q(k(V),\mmu_n^{\otimes j})$ consisting of those classes $\alpha$ such that $\partial_R(\alpha)=0$ for any $R$ as above. If $V$ is proper over $k$, then by \cite[Theorem 6.1]{blochogus} the group $\rmH^0(V,\mathscr{H}^q_\et(\mmu_n^{\otimes j}))$ coincides with $\rmH^q_\mathrm{nr}(k(V),\mmu_n^{\otimes j})$. Furthermore, if we let $V$ and $W$ be two smooth and proper varieties over $k$ which are stably birationnally equivalent, then $\rmH^q_\nr(V,\mmu_n^{\otimes j})\simeq \rmH^q_\nr(W,\mmu_n^{\otimes j})$ \cite[Proposition~1.2]{CTO89}. For a general variety $V$ in characteristic $0$, we may thus let $V^c$ be any smooth compactification and define the unramified cohomology group $\rmH^q_\mathrm{nr}(k(V),\mmu_n^{\otimes j})$ as $\rmH^0(V^c,\mathscr{H}^q_\et(\mmu_n^{\otimes j}))$. 

If $x$ is a $k$-point of $V$, we let $\rmH^{q,\norm}_\nr(k(V),\mmu_n^{\otimes j})$ be the kernel of the evaluation map $x^*:\rmH^0(V^c,\mathscr{H}^q_\et(\mmu_n^{\otimes j}))\to \rmH^q(k,\mmu_n^{\otimes j})$. Over an arbitrary field $k$, we define the groups $\rmH^q_\nr(k(V),(\Q/\Z)'(i))$ in the obvious way.

\subsection{Unramified cohomology and evaluation at points} The central result of this paper is motivated by the following general theorem which is an analogue, in degree $3$, of a criterion proved by Wittenberg for a Brauer class to be unramified; see \cite[Lemma~10.5.12]{CTS21} for details. (Notably, the proof of the following result is directly inspired from the proof of \emph{loc.\ cit.})

\begin{thm}\label{WittenbergH3nr}
Let $k$ be a field of exponential characteristic $p>0$, let $V$ be a smooth variety over a field $k$ and let $\alpha\in \rmH^3_\et(V,(\Q/\Z)'(2))$ be a class whose order is coprime to~$p$. Then the image of $\alpha$ in $\rmH^3(k(V),(\Q/\Z)'(2))$ is unramified if and only if for any perfect field $L$ of cohomological dimension at most~$2$ containing $k$, the evaluation map
\[\mathrm{ev}_\alpha:V(L((t)))\longrightarrow\rmH^3(L((t)),(\Q/\Z)'(2))\]
is the zero map.
\end{thm}

\begin{rmk} This criterion is not restrictive in the sense that it allows one to recover all of the unramified classes in $\rmH^3(k(V),(\Q/\Z)'(2))$. Indeed, a class $\alpha\in \rmH^3(k(V),(\Q/\Z)'(2))$ comes from a class $\alpha_U\in \rmH^3_\et(U,(\Q/\Z)'(2))$ for some open subset $U\subset V$. Hence, up to replacing $V$ by an open subset we may assume that $\alpha$ comes from an element of $\rmH^3_\et(V,(\Q/\Z)'(2))$.
\end{rmk}

\begin{proof}[Proof of Theorem \ref{WittenbergH3nr}] Let $\alpha\in \rmH^3_\et(V,(\Q/\Z)'(2))$ be as above and suppose that its image in the group  $\rmH^3(k(V),(\Q/\Z)'(2))$ is unramified. Fix $L$ as in the statement of the theorem, and let $P:\Spec(L((t)))\to V$ be any morphism. Let $x\in V$ be the image of $P$. Since $V$ is smooth, there exists an open subset $U\subset V$ containing $x$ and an étale morphism $f:U\to \mathbf{A}^n_k$ for some $n\ge 0$. Let $Q=f(P)\in \mathbf{A}^n_k(L((t)))$. There exists $Q'\in \mathbf{A}^n(L((t)))$ arbitrarily close to $Q$ for the topology of $k((t))$ and whose $n$ coordinates are algebraically independent over $L$, so that $Q'$ cannot lie on any proper subvariety of $\mathbf{A}^n_L$. Because $f$ is étale and $Q'$ can be chosen as close to $Q$ as we want, there exists $P'\in U(L((t)))$ arbitrarily close to $P$ such that $f(P')=Q'$. Taking $P'$ close enough to $P$ such that $(P')^*(\alpha)=P^*(\alpha)$, the morphism $P':\Spec (L((t)))\to X$ is dominant and thus induces an embedding $k(V)\subset L((t))$. (See also \cite[Lemma~10.5.11]{CTS21} for a detailed proof of this statement.) 

If we suppose that $k(V)\subset L[[t]]$, then $(P')^*(\alpha)\in \rmH^3(L[[t]],(\Q/\Z)'(2))$. As $L[[t]]$ is henselian and $L$ has cohomological dimension at most $2$, we get that $\rmH^3(L[[t]],(\Q/\Z)'(2))\simeq \rmH^3(L,(\Q/\Z)'(2))=0$ (where the first isomorphism comes from proper base change, cf.\ \cite[Theorem II.3.10]{milneetale}). If $k(V)$ does not inject into $L[[t]]$, then $k(V)\cap L[[t]]$ is a discrete valuation ring containing $k$ and whose fraction field is $k(V)$. Since $\alpha$ is unramified, then by definition it comes from $\rmH^3(k(V)\cap L[[t]],(\Q/\Z)'(2))$, which shows that $(P')^*(\alpha)\in \rmH^3(L[[t]],(\Q/\Z)'(2))=0$ as well. Therefore, we have $P^*(\alpha)=0\in \rmH^3(L((t)),(\Q/\Z)'(2))$.

Conversely, suppose that $P^*(\alpha)=0$ for any $P$ and any $L$ as in the theorem. We aim to show that $\alpha$ is unramified. Since $\alpha$ has order coprime to $p$, then we may assume $k$ to be perfect up to extending scalars to a perfect closure. Let us fix a discrete valuation ring $A$ containing $k$ and with fraction field $k(V)$, and denote by $\kappa$ its residue field. There exists a field $L$ of cohomological dimension at most $2$ containing $\kappa$ in which the latter is algebraically closed and such that the natural map $\Br(\kappa)\to \Br(L)$ is injective, see the main theorem of \cite{ducrosdimensioncohomologique}. Since $k$ is perfect, a theorem of Cohen \cite{cohenstructure} shows that there exists an isomorphism of $k$-algebras between $\kappa[[t]]$ and the completion of $A$. Let us fix such an isomorphism. Taking fraction fields, this provides a $k$-linear embedding of $k(V)$ into $\kappa((t))$, and hence into $L((t))$. This provides an $L((t))$-point $P$ of $V$. 

Since $L$ has cohomological dimension at most $2$, the residue map attached to $A$ induces an isomorphism $\rmH^3(L((t)),(\Q/\Z)'(2))\simeq\rmH^2(L,(\Q/\Z)'(1))$. As $P^*(\alpha)=0$ by assumption, we see that the residue of $\alpha$ with respect to $A$ therefore belongs to the kernel of the restriction map $\rmH^2(\kappa,(\Q/\Z)'(1))\to \rmH^2(L,(\Q/\Z)'(1))$. But the latter coincides with the restriction $\Br(\kappa)\{p'\}\to\Br(L)\{p'\}$, which is injective, hence the desired result. 
\end{proof}

\section{Main theorem}
\subsection{Globally unramified classes}\label{paragraphegéometriquementnonramifié} Let $k$ be a field of exponential characteristic $p> 0$ and let $G$ be a finite $k$-group of order coprime to $p$. Let $(U\subset \A^m_k,c)$ be a classifying pair and let $X:=U/G$ be the geometric quotient. In this section, we define residues on $X$ purely in terms of the structure of the $k$-group $G$ and the arithmetic of the field $k$. 

Let $x\in X(k)$, let $\overline{x}\mapsto x$ be a geometric point and let $\xi:\pi_1^\et(X,\overline{x})\to \Gamma_k$ denote the natural projection. Recall that we may view $\pi_1^\et(X,\overline{x})$ as a semi-direct product $\pi_1^\et(\overline{X},\overline{x})\rtimes\Gamma_k$ as $x$ is a rational point, and that the maximal pro-$p'$ quotient of $\pi_1^\et(\overline{X},\overline{x})$ is isomorphic to $G(\overline{k})$. This moreover exhibits the semi-direct product $G(\overline{k})\rtimes\Gamma_k$ as a quotient of $\pi_1^\et(X,\overline{x})$, denoted $\pi_1^\et(X,\overline{x})^{[p']}$ for simplicity; by abuse of notation, we still write $\xi:\pi_1^\et(X,\overline{x})^{[p']}\to \Gamma_k$ for the natural projection. Let $\chi:\Gamma_k\to \mathrm{Aut}(\widehat{\Z}(1))=\widehat{\Z}^*$ denote the cyclotomic character of $k$. We first introduce the following notation: 
\begin{defn}\label{définitionrésidus} For $\mathrm{u}\in\{\mathrm{any,met}\}$ and $\mathrm{v}\in\{\mathrm{any},\mathrm{{\mathrm{pro-}p'}}, 0\}$, we denote by ${_\mathrm{u}}\mathscr{D}_\mathrm{v}$ the set of pairs $(\mu,D)$ where $\mu:\hat{\Z}'(1)\to G(\overline{k})$ is a morphism (not
assumed to be $\Gamma_k$-equivariant) and $D\subset G(\overline{k})\rtimes \Gamma_k$ is a closed subgroup such that:
\begin{enumerate}
    \item[•]  if $\mathrm{u}=\mathrm{cent}$, then $D\cap G(\overline{k})=C_{G(\overline{k})}(\Im(\mu))$; if $\mathrm{u}=\mathrm{any}$, then we do not impose any condition on $D\cap G(\overline{k})$;
    \item[•] if $\mathrm{v}=0$, then $\xi(D)$ is trivial; if $\mathrm{v}=\mathrm{{\mathrm{pro-}p'}}$, then $\xi(D)$ is a pro-$\ell$-group for some $\ell\neq p$; if $\mathrm{v}=\mathrm{any}$, then we do not impose any condition on $\xi(D)$;
    \item[•] for any $(\gamma,\sigma)\in D$ and for any $\tau\in \mathrm{Im}(\mu)$, we have the equality $(\gamma,\sigma)(\tau,1)(\gamma,\sigma)^{-1}=(\tau^{\chi(\sigma)},1)\in G(\overline{k})\rtimes \Gamma_k$. 
\end{enumerate}
\end{defn}
Let $(\mu,D)\in{_\mathrm{u}}\mathscr{D}_\mathrm{v}$ with $\mathrm{u},\mathrm{v}$ as above. We define a morphism of semi-direct products $\varphi:\widehat{\Z}'(1) \rtimes D\to \langle\mathrm{Im}(\mu), D\rangle$ by letting $\varphi|_{\hat{\Z}'(1)}=\mu$ and $\varphi_{|_{D}}=\mathrm{Id}_D$. Here, the group $D$ acts on $\widehat{\Z}'(1)$ through the cyclotomic character evaluated on $\xi(D)\subset \Gamma_k$. This construction yields a pullback map:
\[\varphi^*:\rmH^3(\pi_1^\et(X,\overline{x})^{[p']},(\Q/\Z)'(2))\longrightarrow\rmH^3(\widehat{\Z}'(1)\rtimes D,(\Q/\Z)'(2)),\]
where $\widehat{\Z}'(1)\rtimes D$ acts on $(\Q/\Z)'(2)$ via $\xi(D)\subset \Gamma_k$. We may now invoke the Hochschild-Serre spectral sequence $E_2^{p,q}=\rmH^p(D,\rmH^q(\widehat{\Z}'(1),(\Q/\Z)'(2)))\Rightarrow\rmH^{p+q}(\widehat{\Z}'(1)\rtimes D,(\Q/\Z)'(2))$, together with the fact that $\widehat{\Z}'(1)$ has cohomological dimension $1$ in order to obtain a residue map:
\[\partial:\rmH^3(\widehat{\Z}'(1)\rtimes D,(\Q/\Z)'(2))\to \rmH^2(D,\rmH^1(\widehat{\Z}'(1),(\Q/\Z)'(2))).\]
As $\rmH^1(\widehat{\Z}'(1),(\Q/\Z)'(2))$ canonically identifies with $\Hom_{\mathrm{cont}}(\widehat{\Z}'(1),(\Q/\Z)'(2))\simeq(\Q/\Z)'(1)$, we thus derive a residue map: 
\[\partial_{\mu,D}: \rmH^3(\pi_1^\et(X,\overline{x})^{[p']},(\Q/\Z)'(2))\longrightarrow \rmH^2(D,(\Q/\Z)'(1)).\]

\begin{defn}\label{definitionglobalementnr}  For $\mathrm{u}\in\{\mathrm{any,met}\}$ and $\mathrm{v}\in\{\mathrm{any},\mathrm{{\mathrm{pro-}p'}}, 0\}$, we define the following subgroups of $\rmH^3(\pi_1^\et(X,\overline{x})^{[p']},(\Q/\Z)'(2))$:
\[H^3_\mathrm{u,v}(\pi_1^\et(X,\overline{x})^{[p']},(\Q/\Z)'(2)):=\bigcap_{(\mu,D)\in {_\mathrm{u}}\mathscr{D}_\mathrm{v}}\Ker\big(\partial_{\mu,D}\big).\]
The subgroup of \emph{globally unramified classes} in $\rmH^3(\pi_1^\et(X,\overline{x})^{[p']},(\Q/\Z)'(2))$ is given by:
\[\rmH^3_\mathrm{gnr}(\pi_1^\et(X,\overline{x})^{[p']},(\Q/\Z)'(2)):=\rmH^3_\mathrm{any,{\mathrm{pro-}p'}}(\pi_1^\et(X,\overline{x})^{[p']},(\Q/\Z)'(2)).\]
\end{defn}

\begin{rmk}\label{remarqueraffinementsrésidus} Let us mention the following two refinements of Definition \ref{definitionglobalementnr}.

(i) If we assume that $k$ is algebraically closed, then we evidently have an inclusion ${_\mathrm{any}\mathscr{D}_\mathrm{{\mathrm{pro-}p'}}}\subset {_\mathrm{any}\mathscr{D}_\mathrm{0}}$.

(ii) There is an equality $\rmH^3_\mathrm{cent,0}(\pi_1^\et(X,\overline{x})^{[p']},(\Q/\Z)'(2))=\rmH^3_\mathrm{any,0}(\pi_1^\et(X,\overline{x})^{[p']},(\Q/\Z)'(2))$. Indeed, for a given pair $(\mu,D)\in {_\mathrm{any}\mathscr{D}_\mathrm{0}}$, the group $D=D\cap G(\overline{k})$ acts trivially on $(\Q/\Z)'(1)$, so that we have an isomorphism $\rmH^2(D,(\Q/\Z)'(1))\simeq \rmH^2(D,(\Q/\Z)')$ where the latter denotes cohomology of finite groups. Let $\tau$ denote a generator of $\Im(\mu)$. As $D$ commutes with $\tau$, we have $D\subset C_{G(\overline{k})}(\tau)$. Conversely, we evidently have $(\mu, C_{G(\overline{k})}(\tau))\in{_\mathrm{any}\mathscr{D}_\mathrm{0}}$. Since the residues are compatible with restrictions, i.e.\ $\partial_{\mu,D'}=\mathrm{Res}^{D}_{D'}\circ\partial_{\mu,D}$ for any pairs $(\mu,D'), (\mu,D)$ with $D' \subset D$, a class $\alpha$ lies in $\rmH^3_\mathrm{any,0}(\pi_1^\et(X,\overline{x})^{[p']},(\Q/\Z)'(2))$ if and only if
$\partial_{\mu,C_{G(\overline{k})}(\tau)}(\alpha)=0 \in \mathrm{H}^2(C_{G(\overline{k})}(\tau),(\Q/\Z)')$ for every such $\mu$, that is, if and only if $\alpha\in \rmH^3_\mathrm{cent,0}(\pi_1^\et(X,\overline{x})^{[p']},(\Q/\Z)'(2))$. Moreover, since residues are compatible with conjugation, it suffices to let $\tau$ run over a set of representatives of the conjugacy classes of $G(\overline{k})$ and $\mu$ over the corresponding morphisms. Hence, deciding whether $\alpha\in\rmH^3_\mathrm{any,0}(\pi_1^\et(X,\overline{x})^{[p']},(\Q/\Z)'(2))$ amounts to checking one condition for each element of $G(\overline{k})$, up to conjugation.
\end{rmk}

\begin{rmk}\label{remarquerésidus} (i) Let $(\mu,D)\in {_\mathrm{any}\mathscr{D}_\mathrm{{\mathrm{pro-}p'}}}$ be a pair such that $\Im(\mu)$ is the trivial subgroup of $G(\overline{k})$. Then the pullback map $\rmH^3(\pi_1^\et(X,\overline{x})^{[p']},(\Q/\Z)'(2))\to \rmH^3(\widehat{\Z}(1)\rtimes D,(\Q/\Z)'(2))$ factors through the inclusion $\rmH^3(D,(\Q/\Z)'(2))\to \rmH^3(\widehat{\Z}(1)\rtimes D,(\Q/\Z)'(2))$ induced by the projection $\widehat{\Z}(1)\rtimes D\to D$. In particular, the associated residue $\partial_{\mu,D}$ is identically zero.

(ii) When $\xi(D)=\{1\}\subset\Gamma_k$ in Definition \ref{définitionrésidus}, for instance if $k=\overline{k}$, then one can check that the construction of the residue map $\partial_{\mu,D}$ refines the group theoretic residue defined by Peyre \cite[Definition 2.7]{peyreh3noether} and Kahn--Nguyen \cite[Définition 1.9]{kahnnguyen}. Indeed, in this case the group $D$ acts trivially on $\widehat{\Z}'(1)$. Denote 
\[N\rmH^3(\mathrm{Im}(\mu)\times D,(\Q/\Z)'):=\Ker\Big(\rmH^3(\mathrm{Im}(\mu)\times D,(\Q/\Z)')\to \rmH^3(\mathrm{Im}(\mu),(\Q/\Z)')\Big)\] and 
\[N\rmH^3(\widehat{\Z}'(1)\times D,(\Q/\Z)'):=\Ker\Big(\rmH^3(\widehat{\Z}'(1)\times D,(\Q/\Z)')\to \rmH^3(\widehat{\Z}'(1),(\Q/\Z)')\Big).\]
The Hochschild-Serre spectral sequences for $\widehat{\Z}'(1)\times D$ and $\Im(\mu)\times D$ yield together a commutative diagram:
\[{\begin{tikzcd}
	{\rmH^3(\Im(\mu)\times D,(\Q/\Z)')} & {N\rmH^3(\mathrm{Im}(\mu)\times D,(\Q/\Z)')} & {\rmH^2(D,\rmH^1(\mathrm{Im}(\mu),(\Q/\Z)'))} \\
	{\rmH^3(\widehat{\Z}'(1)\times D,(\Q/\Z)')} & {N\rmH^3(\widehat{\Z}'(1)\times D,(\Q/\Z)')} & {\rmH^2(D,(\Q/\Z)'),}
	\arrow[from=1-1, to=1-2]
	\arrow[from=1-1, to=2-1]
	\arrow[from=1-2, to=1-3]
	\arrow[from=1-2, to=2-2]
	\arrow[from=1-3, to=2-3]
	\arrow[from=2-1, to=2-2]
	\arrow[from=2-2, to=2-3]
\end{tikzcd}}\]
where the rightmost vertical map is induced by the injection $\rmH^1(\mathrm{Im}(\mu),(\Q/\Z)')\simeq \Hom(\mathrm{Im}(\mu),(\Q/\Z)')\hookrightarrow (\Q/\Z)'$. The bottom composite map is precisely the residue $\partial_{\mu,D}$, while the composite map 
\[\rmH^3(\mathrm{Im}(\mu)\times D,(\Q/\Z)')\to \rmH^2(D,\rmH^1(\mathrm{Im}(\mu),(\Q/\Z)'))\to \rmH^2(D,(\Q/\Z)')\]
agrees with the aforementioned residues.
\end{rmk}

\subsection{Permutation negligible classes}\label{paragraphepermutationnégligeable} In this section, we define a subgroup of distinguished classes in $\rmH^3(\pi_1^\et(X,\overline{x})^{[p']},(\Q/\Z)'(2))$ which arise by pushing forward intersection products of line bundles along some finite good covers of $X$. 

\subsubsection{Modified cup-products} Let $\pi:Y\to X$ be a finite good cover. In the derived category of bounded complexes of abelian étale sheaves on $Y$, let us consider the distinguished triangle $\Z\to \Q\to \Q/\Z\to \Z[1].$ Applying the functor $(-)\otimes^{\mathbf{L}}\Q/\Z$ to the map $\Q/\Z\to \Z[1]$, we obtain a map:
    \[\Q/\Z\otimes^\mathbf{L}\Q/\Z\longrightarrow \Q/\Z[1],\]
    and similarly with Tate twists:
    \[\Q/\Z(i)\otimes^\mathbf{L}\Q/\Z(j)\longrightarrow \Q/\Z(i+j)[1]\]
    for $i,j\in \Z$. Taking cohomology, this yields cup-products:
    \[\rmH^p_\et(Y,\Q/\Z(i))\otimes_\Z \rmH^q_\et(Y,\Q/\Z(j))\longrightarrow \rmH^{p+q+1}_\et(Y,\Q/\Z(i+j)).\]
The same observation can be extended to motivic cohomology by remarking that the groups $\rmH^*_\mot(Y,\Q/\Z(j))$ can be computed as the hypercohomology of the cdh-motivic complexes $\Z(j)\otimes^{\mathbf{L}} \Q/\Z$. Combining this fact with the existence of products $\Z(i)\otimes^{\mathbf{L}} \Z(j)\to \Z(i+j)$ in the derived category of complexes of cdh-sheaves on $Y$ \cite[Proposition 2.8]{suslinvoevodsky}, we obtain cup-products:
 \[\rmH^p_\mot(Y,\Q/\Z(i))\otimes_\Z \rmH^q_\mot(Y,\Q/\Z(j))\longrightarrow \rmH^{p+q+1}_\mot(Y,\Q/\Z(i+j)),\]
which are compatible with the étale realisation maps.

\subsubsection{Permutation classes}\label{definitionnegligible} Let now $K$ be a finite separable field extension of $k$ contained in $\overline{k}$, and let $H\subset G\otimes_k K$ be a $K$-subgroup. The $\Gamma_K$-action on $H(\overline{k})$ and the $\Gamma_k$-action on $G(\overline{k})$ together yield an injective morphism of semi-direct products $H(\overline{k})\rtimes\Gamma_K\to G(\overline{k})\rtimes\Gamma_k$, which defines a finite étale morphism $\pi:Y\to X$ where $Y:=U_K/H$. 

We may view $u\in U(k)$ as a $K$-point of $U$ and consider its image in $Y$. If we denote this point by $y\in Y(K)$, then it provides a retraction of the edge map $\rmH^1(K,(\Q/\Z)'(1))\simeq \rmH^0(Y,\Gm)\otimes_\Z(\Q/\Z)'\to \rmH^1_\et(Y,(\Q/\Z)'(1))$ which appears in the short exact sequence 
\[0\longrightarrow \rmH^0(Y,\Gm)\otimes_\Z(\Q/\Z)'\longrightarrow \rmH^1_\et(Y,(\Q/\Z)'(1))\longrightarrow \Pic(Y)\{p'\}\longrightarrow 0.\]
induced by the Kummer sequence on $Y_\et$. Since $\Pic(U_K)\simeq \Pic(\A^m_K)$ and the latter group is trivial, we see by restriction-corestriction for the map $U_K\to Y$ that $\Pic(Y)$ is torsion, killed by the order of $H$. We hence obtain an inclusion $\Pic(Y)\hookrightarrow \rmH^1_\et(Y,(\Q/\Z)'(1))$ and \emph{a fortiori} a product map:
    \[\Pic(Y)^{\otimes 2}\longrightarrow \rmH^3_\et(Y,(\Q/\Z)'(2))\simeq \rmH^3(\pi_1^\et(Y,\overline{y})^{[p']},(\Q/\Z)'(2)),\]
    where the latter isomorphism comes from Lemma \ref{cohomologieequivariante}. Subsequently, the pushforward map $\pi_*:\rmH^3_\et(Y,(\Q/\Z)'(2))\to \rmH^3_\et(X,(\Q/\Z)'(2))$ coincides with the corestriction associated to the inclusion of the open subgroup $\pi_1^\et(Y,\overline{y})^{[p']}\subset \pi_1^\et(X,\overline{x})^{[p']}$.
    \begin{defn} 
    
We define the group $\rmH^3_\mathrm{per}(\pi_1^\et(X,\overline{x})^{[p']},(\Q/\Z)'(2))$ of \emph{permutation negligible classes} within $\rmH^3(\pi_1^\et(X,\overline{x})^{[p']},(\Q/\Z)'(2))$ as the subgroup generated by the corestrictions
\[\pi_*\left(\mathrm{Im}\big(\Pic(Y)^{\otimes 2}\overset{\smile}{\longrightarrow} \rmH^3(\pi_1^\et(Y,\overline{y})^{[p']},(\Q/\Z)'(2))\big)\right),\] where $\pi:Y\to X$ ranges over the good covers described above.
    \end{defn}
\subsection{Statement of the main theorem} We now conclude this section by stating the main theorem of this paper, which we will prove in §\ref{preuvetheoremeprincipal}. Let us first recall some classical facts about vector bundles on classifying schemes. Throughout the following section, we let $G$ be a finite $k$-group of order coprime to $p$ and $c\ge 3$, and fix a classifying pair $(U\subset \A^m_k,c)$. 
\subsubsection{Reminders on $K$-theory and representations}\label{rappelsKthéorie} Consider the monoidal category $\mathbf{Rep}_k^G$ of linear $k$-representations of $G$, that is, the category whose objects are morphisms of algebraic $k$\nobreakdash-groups $G\to \GL(V)$ with $V$ a $k$-vector space of finite dimension, and whose morphisms are morphisms of $G$-representations defined over $k$. Since $k$ has exponential characteristic $p>0$ and $G$ is finite of order coprime to $p$, this category is semi-simple \cite[Theorem 22.42]{milnealgebraicgroups}. Denote by $\mathrm{Rep}_k(G)$ the Grothendieck group of this category \cite[§14.1]{serrereplin}, that is, the Grothendieck group of virtual $k$-representations of $G$.


Denote by $X:=U/G$ the geometric quotient and let $K_0(X)$ be the Grothendieck group of the exact category of vector bundles on $X$. Since the latter is a smooth variety, we let $F^\bullet K_0(X)$ be the filtration by codimension of the support, with graded quotients $\mathrm{gr}^rK_0(X):=F^rK_0(X)/F^{r+1}K_0(X)$ for $r\ge 0$; for details, see e.g.\ \cite[§7.5]{quillen} and \cite[Exposé II, §2.2]{SGA6}. According to \cite[Example 15.1.5]{fulton}, for each $r\ge 0$ we have a surjective $K$-theoretic cycle class map $\phi^r:\CH^r(X)\to \mathrm{gr}^rK_0(X)$. As explained in \cite[Example~15.3.6]{fulton}, Jouanolou's Riemann-Roch theorem without denominators shows that $\phi^0$, $\phi^1$ and $\phi^2$ are isomorphisms, with inverses given by the rank $\mathrm{rk}$, the determinant $\mathrm{det}$ and the opposite $-c_2$ of the second Chern class. 

Given the assumptions on $X$, one can define a natural morphism $\mathrm{Rep}_k(G)\to K_0(X)$ as follows: if $G\to \GL(V)$ is a $k$-representation, consider the $k$-variety $V\times_k U$ endowed with the linear action $g\cdot(v,u):=(g\cdot v,u\cdot g^{-1})$. The geometric quotient $(V\times_k U)/G$ is a vector bundle over $X$ whose class in $K_0(X)$ only depends on that of the prescribed representation in $\mathrm{Rep}_k(G)$. In the proof of \cite[Theorem~2.25]{totaro}, Totaro defines a filtration $\mathrm{F}^\bullet \mathrm{Rep}_k(G)$ induced by the filtration by codimension of the support on $K_0(X)$; moreover, he shows that this filtration yields well-defined maps $\mathrm{gr}^i \mathrm{Rep}_k(G)\to \mathrm{gr}^i K_0(X)$ for $i\le c$, which are isomorphisms. In particular, any element of $\CH^2(X)$ is the second Chern class of some virtual $k$-representation of $G$.

\begin{thm}\label{theoremeprincipal} Let $k$ be a field of exponential characteristic $p>0$, let $G$ be a finite $k$-group of order coprime to $p$, let $c\ge 3$ and fix a classifying pair $(U\subset \A^m_k,c)$. If we denote $X:=U/G,$ then the natural map $\rmH^3_\et(X,(\Q/\Z)'(2))\to \rmH^3(k(X),(\Q/\Z)'(2))$ induces a short exact sequence 
\[0\longrightarrow \CH^2(X)\longrightarrow \rmH^3_\mathrm{gnr}(\pi_1^\et(X,\overline{x})^{[p']},(\Q/\Z)'(2))\longrightarrow \rmH^3_\nr(k(X),(\Q/\Z)'(2))\longrightarrow 0,\]
where the group $\CH^2(X)$ is generated by Chern classes of $k$-representations of $G$. If we further assume that $G$ is constant, then we have a surjection:
\[\frac{\rmH^3_\mathrm{gnr}(\pi_1^\et(X,\overline{x})^{[p']},(\Q/\Z)'(2))}{\rmH^3_\mathrm{per}(\pi_1^\et(X,\overline{x})^{[p']},(\Q/\Z)'(2))}\twoheadrightarrow \rmH^3_\nr(k(X),(\Q/\Z)'(2)),\]
the kernel of which is killed by a power of $2$. In particular if $G$ has odd order, then the above map is an isomorphism.
\end{thm}

\begin{rmk} If we assume that $k$ is algebraically closed, then Remark \ref{remarqueraffinementsrésidus}, (i) and (ii) imply that we have an equality $\rmH^3_\mathrm{gnr}(\pi_1^\et(X,\overline{x}),(\Q/\Z)'(2))=\rmH^3_{\mathrm{cent},0}(\pi_1^\et(X,\overline{x}),(\Q/\Z)'(2))$. By virtue of Remark \ref{remarquerésidus}, we hence see that for an algebraically closed field $k$, Theorem \ref{theoremeprincipal} is a refinement of \cite[Theorem 1]{peyreh3noether}.\end{rmk}

\section{Proof of the main theorem}\label{preuvetheoremeprincipal}

Our goal in this section is to prove Theorem \ref{theoremeprincipal}.
\subsection{A fundamental exact sequence}\label{paragraphesuitefondamentale} The following proposition is an adapted version, over an arbitrary field, of an argument initially due to Colliot-Thélène \cite[Proposition~1]{peyreh3noether}. We provide a full proof for the sake of completeness, as it does not appear in the literature (at least to our knowledge):
\begin{prop}\label{suiteexactefondamentale}
Let $k$ be a field of exponential characteristic $p>0$, let $G$ be a finite $k$-group of order coprime to $p$, let $c\ge 3$ and fix a classifying pair $(U\subset\A^m_k,c)$. If we denote $X:=U/G,$ then there is a short exact sequence of abelian groups:
\begin{equation}\label{suitefondamentale}
0\longrightarrow\CH^2(X)\longrightarrow \rmH^3(\pi_1^\et(X,\overline{x})^{[p']},(\Q/\Z)'(2))\longrightarrow\rmH^0(X,\mathscr{H}^3_\et((\Q/\Z)'(2)))\longrightarrow 0,
\end{equation}
where the first map factors the $\ell$-adic cycle class map $\CH^2(X)\{\ell\}\to \rmH^4(\pi_1^\et(X,\overline{x})^{[p']},\Z_\ell(2))$ for every prime number $\ell\neq p$.
\end{prop}

The map $\CH^2(X)\to \rmH^3(\pi_1^\et(X,\overline{x}),(\Q/\Z)'(2))$ appearing in \eqref{suitefondamentale} is an incarnation of Bloch's Abel-Jacobi map $\CH^2(X)_\mathrm{tors}\to \rmH^3_\et(X,(\Q/\Z)'(2))$ on torsion cycles, defined over an arbitrary field. Throughout the rest of this text, we will denote it by $\mathrm{cl}_{2,\mathrm{AJ}}$. 

We recall that a variety $X$ over a field $k$ is said to be \emph{cellular} if there exists a finite filtration $\emptyset= Z_n\subset Z_{n-1}\subset \hdots\subset Z_0=X$ such that for each $0\le i\le n-1$, $Z_{i+1}$ is a proper closed subset of $Z_{i}$ and $Z_{i}\setminus Z_{i+1}$ is isomorphic to an affine space over $k$. For $j\ge 0$, let $\mathscr{K}_j$ denote the sheafification, on $X_\mathrm{Zar}$, of the presheaf $V\mapsto K_j(\Gamma(V,\mathscr{O}_X))$ (where $K_j$ denotes Quillen's $K$-theory). In order to prove the above assertion, we crucially make use of the following result (see also the proof of \cite[Proposition 3.7]{arasoninvariant}):
\newcommand{\macitation}{\cite[Corollary 2.2.2]{peyremotiviccomplexes}}
\begin{lem}[\macitation]\label{lemmecelullaire}
Let $X$ be a smooth cellular variety over a field $k$ together with a closed subvariety $Z$ of codimension $c\ge 2$. Then for any $i\le c-2$, there is an isomorphism
\[\rmH^i_\mathrm{Zar}(X\setminus Z,\mathscr{K}_j)\overset{\sim}{\longrightarrow} \CH^i(X)\otimes_\Z K_{j-i}(k)\]
\end{lem}
\begin{proof}[Proof of Proposition \ref{suiteexactefondamentale}] By Bloch-Ogus theory \cite[Corollary 6.2]{blochogus}, the Leray spectral sequence $E_2^{p,q}=\rmH^p_\mathrm{Zar}(X,\mathscr{H}^q_\et((\Q/\Z)'(2)))\Rightarrow\rmH^{p+q}_\et(X,(\Q/\Z)'(2))$ for the change-of-sites map $X_\et\to X_\mathrm{Zar}$ gives rise to an exact sequence of lower terms:
\begin{equation}\begin{tikzcd}
  0 \rar & \rmH^1_\mathrm{Zar}(X,\mathscr{H}^2_\et((\Q/\Z)'(2))) \rar
             \ar[draw=none]{d}[name=X, anchor=center]{}
    & \rmH^3_\et(X,(\Q/\Z)'(2)) \ar[rounded corners,
            to path={ -- ([xshift=2ex]\tikztostart.east)
                      |- (X.center) \tikztonodes
                      -| ([xshift=-2ex]\tikztotarget.west)
                      -- (\tikztotarget)}]{dll}[at end]{} \\      
  \rmH^0(X,\mathscr{H}^3_\et((\Q/\Z)'(2))) \rar & \CH^2(X)\otimes_\Z(\Q/\Z)'. & 
\end{tikzcd}
\end{equation}
On the other hand, a standard argument due to Bloch \cite[§3.2]{colliotthélèneCIME} yields a short exact sequence:
\[0\longrightarrow \rmH^1_\mathrm{Zar}(X,\mathscr{K}_2)\otimes_\Z (\Q/\Z)'\longrightarrow\rmH^1_\mathrm{Zar}(X,\mathscr{H}^2_\et((\Q/\Z)'(2)))\longrightarrow \CH^2(X)\{p'\}\longrightarrow 0.\]
Let $Z$ be the closed complement of $U\subset \A^m_k$. The localisation sequence for Chow groups \cite[Proposition 1.8]{fulton} then reads 
\[\CH_{m-2}(Z)\longrightarrow \CH_{m-2}(\A^m_k)\longrightarrow\CH_{m-2}(U)\longrightarrow 0,\]
and the group $\CH_{m-2}(\A^m_k)$ is trivial; see \emph{loc.\ cit.}, Chapter I, §1.9. Thanks to a restriction-corestriction argument for the finite morphism $U\to X$, we see that $\CH^2(X)$ is killed by the order of $G$. By Lemma \ref{lemmecelullaire}, we have $\rmH^1_\mathrm{Zar}(U,\mathscr{K}_2)\simeq \Pic(\A^m_k)\otimes_\Z k^*$, and the latter is trivial. Hence, we infer by restriction-corestriction that $\rmH^1_\mathrm{Zar}(X,\mathscr{K}_2)$ is also killed by $|G|$. We hence obtain an isomorphism  $\rmH^1(X,\mathscr{H}_\et^2((\Q/\Z)'(2)))\simeq \CH^2(X)$, and we have a short exact sequence:
\[0\longrightarrow \CH^2(X)\longrightarrow\rmH^3_\et(X,(\Q/\Z)'(2))\longrightarrow \rmH^0(X,\mathscr{H}^3_\et((\Q/\Z)'(2)))\longrightarrow 0.\]
Given the assumptions on $X$, it follows from Lemma \ref{cohomologieequivariante} that there is an isomorphism $\rmH^3_\mathrm{ét}(X,(\Q/\Z)'(2))\simeq \rmH^3(\pi_1^\et(X,\overline{x})^{[p']},(\Q/\Z)'(2))$. 

We know from \cite[Proposition~1]{ctss} that the map $\CH^2(X)[\ell^n]\to \rmH^3(\pi_1^\et(X,\overline{x})^{[p']},\mmu_{\ell^m}^{\otimes 2})$ factors the étale $\ell$-adic cycle class map $\CH^2(X)[\ell^n]\to \rmH^4(\pi_1^\et(X,\overline{x})^{[p']},\mmu_{\ell^m}^{\otimes 2})$ for every $n,m\ge 1$ and every prime $\ell\neq p$. Passing to the colimit over $n\ge 1$ (resp.\ the limit over $m\ge 1$) for the evident transition maps, the assertion follows.\end{proof}

\subsection{A criterion for being unramified}\label{paragrapheexplicationcritère} Let $k$ be a field of exponential characteristic $p>0$ and let $G$ be a finite $k$-group of order coprime to $p$. Let $c\ge3$, fix a classifying pair $(U\subset \A^m_k,c)$, and denote $X:=U/G$. Let $L$ be a field of cohomological dimension $2$ containing $k$ and such that $\Gamma_L$ is a pro-$\ell$-group for some prime $\ell\neq p$. The choice of some $L((t))$-point $P$ of $X$ and the covariant functoriality of the étale fundamental group give rise, up to conjugation by an element of $\pi_1^\et(\overline{X},\overline{x})$, to a commutative square:
\[
    \begin{tikzcd} 
          \Gamma_{L((t))}\ar[d, "\phi_P"'] \ar[dr, phantom] \ar[r] & \Gamma_L\ar[d]  \\
        \pi_1^\text{ét}(X,\overline{x}) \ar[r, "\xi"] & \Gamma_k,
    \end{tikzcd}
\]
where $\xi:\pi_1^\et(X,\overline{x})\to\Gamma_k$ denotes the natural projection. Let $I\subset \Gamma_{L((t))}$ denote the inertia subgroup and $\Gamma_{\text{t}}:=\Gal(L((t))_\text{t}/L((t)))$ where $L((t))_\text{t}$ is the maximal tamely ramified extension of $L((t))$. We have a short exact sequence 
\[1\longrightarrow I_\text{t}\longrightarrow \Gamma_{\text{t}}\longrightarrow \Gamma_L\longrightarrow 1\]
where $I_\text{t}$ denotes the tame inertia subgroup, see e.g.\ \cite[§1.4.1]{CTS21}. The latter group is a quotient of $I$ and is canonically isomorphic to $\hat{\Z}'(1)$. Moreover, the above sequence naturally splits, and the action of $\Gamma_L$ on $I_\text{t}\cong \hat{\Z}'(1)$ is given by the cyclotomic character of $L$, cf.\ \emph{loc.\ cit., ibidem}. This hence provides a canonical isomorphism $\Gamma_{\text{t}}\cong \hat{\Z}'(1)\rtimes \Gamma_L$.

Therefore, the choice of an $L((t))$-point induces up to conjugation by an element of $\pi_1^\et(\overline{X},\overline{x})^{(p')}$ a homomorphism $I\to \pi_1^\et(\overline{X},\overline{x})^{(p')}$ such that we have a commutative diagram with exact rows:
\[
    \begin{tikzcd} 
        1 \ar[r] & I\ar[r] \ar[d, dashed] &  \Gamma_{L((t))} \ar[d, "\phi_P"'] \ar[dr, phantom] \ar[r] & \Gamma_L\ar[d] \ar[r]& 1 \\
        1 \ar[r] & \pi_1^\et(\overline{X},\overline{x})^{(p')} \ar[r]& \pi_1^\text{ét}(X,\overline{x})^{[p']} \ar[r, "\xi"] & \Gamma_k \ar[r]& 1.
    \end{tikzcd}
\]

Using Lemma \ref{cohomologieequivariante}, we see that asking whether a class $\alpha\in\rmH^3_\et(X,(\Q/\Z)'(2))$ evaluates trivially along a given $L((t))$-point $P:\Spec(L((t)))\to X$ amounts to asking whether the pullback map 
\[\rmH^3_\et(X,(\Q/\Z)'(2))\simeq\rmH^3(\pi_1^\et(X,\overline{x})^{[p']},(\Q/\Z)'(2))\overset{\phi_P^*}{\longrightarrow} \rmH^3(\Gamma_{L((t))},(\Q/\Z)'(2))\] sends $\alpha$ to zero. Moreover since $L$ is perfect of cohomological dimension at most $2$, the residue map $\rmH^3(L((t)),(\Q/\Z)'(2))\to \rmH^2(L,(\Q/\Z)'(1))\simeq\Br(L)$ coming from the Hochschild-Serre spectral sequence
\[E_2^{p,q}=\rmH^p(L,\rmH^q(I,(\Q/\Z)'(2)))\Rightarrow \rmH^{p+q}(L((t)),(\Q/\Z)'(2))\]
is an isomorphism, see \cite[Theorem 6.3.5]{gilleszamuely}. Up to passing to the tame quotient $\Gamma_{\text{t}}$ and by comparing the associated residue maps, we have a commutative diagram:
\[\begin{tikzcd}
	{\rmH^3(L((t)),(\Q/\Z)'(2))} & {\Br(L)} \\
	{\rmH^3(\Gamma_{\text{t}},(\Q/\Z)'(2))} & {\Br(L).}
	\arrow["\sim", from=1-1, to=1-2]
	\arrow[from=2-1, to=1-1]
	\arrow["\sim"', from=2-1, to=2-2]
	\arrow[equals, from=2-2, to=1-2]
\end{tikzcd}\]
Consequently, the pullback map $\rmH^3(\Gamma_{\text{t}},(\Q/\Z)'(2))\to \rmH^3(L((t)),(\Q/\Z)'(2))$ is an isomorphism. We hence infer that if the prescribed map $I\to \pi_1^\et(\overline{X},\overline{x})^{(p')}$ does not factor through a non-trivial map $I_\text{t}\to G(\overline{k})\subset \pi_1^\et(\overline{X},\overline{x})^{(p')}$, then $\varphi_P^*(\alpha)=0\in \rmH^3(L((t)),(\Q/\Z)'(2))$. 

Our claim is that a careful analysis of the map $I_\text{t}\rtimes\Gamma_L\to \pi_1^\et(X,\overline{x})^{[p']}$ for varying choices of $L((t))$-points suffices to determine all the classes $\alpha\in\rmH^3_\et(X,(\Q/\Z)'(2))$ whose image in $\rmH^3(k(X),(\Q/\Z)'(2))$, or equivalently, in $\rmH^0(X,\mathscr{H}^3_\et((\Q/\Z)'(2)))$, is unramified. 

\subsection{Inverse image of the unramified subgroup} Let us write $M_\nr(X)$ for the inverse image of the unramified subgroup $\rmH^3_\mathrm{nr}(k(X),(\Q/\Z)'(2))\subset \rmH^0(X,\mathscr{H}^3_\et((\Q/\Z)'(2)))$ along the map $\rmH^3(\pi_1^\et(X,\overline{x})^{[p']},(\Q/\Z)'(2))\to \rmH^0(X,\mathscr{H}^3_\et((\Q/\Z)'(2)))$ that appears in \eqref{suiteexactefondamentale}.
    \begin{lem}\label{fonctorialitegnr}
        Let $\pi:Y\to X$ be a good cover and suppose that $Y$ is irreducible. The pullback morphism $\pi^*:\rmH^3_\et(X,(\Q/\Z)'(2))\to \rmH^3_\et(Y,(\Q/\Z)'(2))$ naturally restricts to a map $M_\nr(X)\to M_\nr(Y).$
    \end{lem}
    \begin{proof} Let $K$ be the field of definition of $Y$. By virtue of the functoriality of unramified cohomology under dominant morphisms \cite[p.\ 143]{CTO89}, the image of any $\alpha\in\rmH^3_\mathrm{nr}(k(X),(\Q/\Z)'(2))$ in $\rmH^3(K(Y),(\Q/\Z)'(2))$, or equivalently in $\rmH^0(Y,\mathscr{H}^3_\et((\Q/\Z)'(2)))$, is unramified as well. In view of the functoriality of the Leray spectral sequence for the change-of-site maps $X_\et\to X_\mathrm{Zar}$ and $Y_\et\to Y_\mathrm{Zar}$ respectively, there is a commutative diagram:
\[\begin{tikzcd}
	{\rmH^3_\et(X,(\Q/\Z)'(2))} & {\rmH^0(X,\mathscr{H}^3_\et((\Q/\Z)'(2)))} \\
	{\rmH^3_\et(Y,(\Q/\Z)'(2))} & {\rmH^0(Y,\mathscr{H}^3_\et((\Q/\Z)'(2)))}.
	\arrow[from=1-1, to=1-2]
	\arrow[from=1-1, to=2-1]
	\arrow[from=1-2, to=2-2]
	\arrow[from=2-1, to=2-2]
\end{tikzcd}\]
Combining these observations with Lemma \ref{cohomologieequivariante} applied to $Y$ readily implies the result.
    \end{proof}

\begin{prop}\label{imageinverse}
    We have an equality
    \[M_\nr(X)=\rmH^3_{\mathrm{gnr}}(\pi_1^\et(X,\overline{x})^{[p']},(\Q/\Z)'(2))\]
    of subgroups of $\rmH^3(\pi_1^\et(X,\overline{x})^{[p']},(\Q/\Z)'(2))$.
\end{prop}
\begin{proof} First suppose that $\alpha\in \rmH^0(X,\mathscr{H}^3_\et((\Q/\Z)'(2)))$ does not belong to the subgroup $\rmH^3_\nr(k(X),(\Q/\Z)'(2))$. If we let $\widetilde{\alpha}\in\rmH^3_\et(X,(\Q/\Z)'(2))$ be a lift, we may suppose it to be of $\ell$-primary torsion for some prime~$\ell\neq p$. By Theorem \ref{WittenbergH3nr}, there exist a perfect field $L$ containing $k$ of cohomological dimension at most $2$ and some $P\in X(L((t)))$ such that $P^*(\widetilde{\alpha})\in\rmH^3(L((t)),(\Q/\Z)'(2))\simeq\Br(L)$ is non-zero. Up to passing to a pro-$\ell$-Sylow subgroup of $\Gamma_L$, we may further assume that the latter group is pro-$\ell$. By virtue of the discussion of §\ref{paragrapheexplicationcritère}, we may also assume that the induced map $\phi_P:\Gamma_{L((t))}\to \pi_1^\et(X,\overline{x})^{[p']}$ factors through a map $\Gamma_{\text{t}}\simeq\hat{\Z}'(1)\rtimes\Gamma_L\to \pi_1^\et(X,\overline{x})^{[p']}$. The latter fits into a commutative diagram of extensions:
\[
    \begin{tikzcd} 
        1 \ar[r] & \hat{\Z}'(1) \ar[r] \ar[d, dashed] &  \hat{\Z}'(1) \rtimes\Gamma_{L} \ar[d, "\phi_P"'] \ar[dr, phantom] \ar[r] & \Gamma_L\ar[d, "\tilde{\sigma}\mapsto\sigma"] \ar[r]& 1 \\
        1 \ar[r] & G(\overline{k}) \ar[r]& \pi_1^\text{ét}(X,\overline{x})^{[p']} \ar[r, "\xi"] & \Gamma_k \ar[r]& 1,
    \end{tikzcd}
\]
where the leftmost vertical map $\mu:\widehat{\Z}'(1)\to G(\overline{k})$ is induced by the commutativity of the right hand square and further lands into a cyclic subgroup. 

Let $\zeta\in \hat{\Z}'(1)$ be a topological generator and denote by $(\zeta,1)$ its image in $\hat{\Z}'(1)\rtimes\Gamma_L$. The group $\mathrm{Im}(\varphi_P)$ is generated by $\Im(\mu)\subset G(\overline{k})$ and by the image $D$ of $\Gamma_L\to \pi_1^\et(X,\overline{x})^{[p']}$. For $\tilde{\sigma}\in \Gamma_L$, write $(\gamma,\sigma):=\varphi_P(0,\tilde{\sigma})$, where we are using the isomorphism $\pi_1^\text{ét}(X,\overline{x})^{[p']}\simeq G(\overline{k})\rtimes \Gamma_k$. The semi-direct product structure on $\hat{\Z}'(1)\rtimes \Gamma_L$ is given by the formula $(0,\sigma)(\zeta,1)(0,\sigma)^{-1}=(\zeta^{\chi(\sigma)},1)$, where $\chi:\Gamma_L\to \mathrm{Aut}(\widehat{\Z}'(1))=(\widehat{\Z}')^*$ is the cyclotomic character. By the very definition of $\varphi_P$, we thus see that $(\gamma,\sigma)(\tau,1)(\gamma,\sigma)^{-1}=(\tau^{\chi(\sigma)},1)$ for every $\tau\in \Im(\mu)$ and every $(\gamma,\sigma)\in  D$. Note that $\widetilde{\alpha}$ does not die when pulled back to $\rmH^3(\langle\Im(\mu),D\rangle,(\Q/\Z)'(2))$, as it does not when pulled back to $\rmH^3(\Gamma_{\text{t}},(\Q/\Z)'(2))$. The morphism $\varphi_P$ may then be factored as 
\[\hat{\Z}'(1)\rtimes \Gamma_L\overset{\varphi_1}{\longrightarrow}\hat{\Z}'(1)\rtimes D\overset{\varphi_2}{\longrightarrow}\langle\Im(\mu),D\rangle\subset \pi_1^\et(X,\overline{x})^{[p']},\]
    such that $\phi_2$ is the composition of $(\mu,\mathrm{Id})$ and the multiplication $\Im(\mu)\rtimes D\to \langle \Im(\mu),D\rangle$, and such that $D$ acts on $\hat{\Z}'(1)$ through its image in $\Gamma_k$ and $\chi$. 
    By the functoriality of the Hochschild-Serre spectral sequence, we obtain a commutative diagram:
\[\begin{tikzcd}
	{\rmH^3(\pi_1^\et(X,\overline{x})^{[p']},(\Q/\Z)'(2))} & {\rmH^3(\langle\Im(\mu),D\rangle,(\Q/\Z)'(2))} & {\rmH^3(\Gamma_{\text{t}},(\Q/\Z)'(2))} \\
	& {\rmH^2(D,(\Q/\Z)'(1))} & {\rmH^2(L,(\Q/\Z)'(1))}
	\arrow[from=1-1, to=1-2]
	\arrow["{\varphi_P^*}", from=1-2, to=1-3]
	\arrow["{\partial_{\mu,D}}"', from=1-2, to=2-2]
	\arrow["\wr", from=1-3, to=2-3]
	\arrow[from=2-2, to=2-3]
\end{tikzcd}\]
which forces $\partial_{\mu,D}(\widetilde{\alpha})\neq 0$, i.e.\ $\widetilde{\alpha}\notin \rmH^3_\mathrm{any,{\mathrm{pro-}p'}}(\pi_1^\et(X,\overline{x})^{[p']},(\Q/\Z)'(2))$. This proves the inclusion $ \rmH^3_\mathrm{gnr}(\pi_1^\et(X,\overline{x})^{[p']},(\Q/\Z)'(2))\subset M_\nr(X)$.

    Conversely, assume that $\alpha\in M_\nr(X)$ and fix a pair $(\mu,D)\in {_\mathrm{any}\mathscr{D}_\mathrm{{\mathrm{pro-}p'}}}$, subject to the semi-commutation relation 
    \[(\gamma,\sigma)(\tau,1)(\gamma,\sigma)^{-1}=(\tau^{\chi(\sigma)},1)\in G(\overline{k})\rtimes \Gamma_k\simeq \pi_1^\et(X,\overline{x})^{[p']}\]
    for all $\tau\in \mathrm{Im}(\mu)$ and $(\gamma,\sigma)\in D$. By virtue of the general theory of étale covers \cite[Exposé~V, Théorème~4.1]{SGA1}, there exists an irreducible pro-étale cover $f:Y\to X$ together with a geometric point $\overline{y}$ such that the pair $(Y,\overline{y})$ is characterised by the fact that $f$ is an inverse limit of good covers and $\pi_1^\et(Y,\overline{y})^{[p']}=\langle\Im(\mu),D\rangle$; therefore, the restriction $\rmH^3(\pi_1^\et(X,\overline{x})^{[p']},(\Q/\Z)'(2))\to \rmH^3(\langle\Im(\mu),D\rangle,(\Q/\Z)'(2))$ along the inclusion of the closed subgroup $\langle \Im(\mu),D\rangle\subset \pi_1^\et(X,\overline{x})^{[p']}$ agrees with the pullback $f^*:\rmH^3_\et(X,(\Q/\Z)'(2))\to\rmH^3_\et(Y,(\Q/\Z)'(2))$. Applying Lemma~\ref{fonctorialitegnr}, this morphism restricts to a map $M_\nr(X)\to M_\nr(Y)$, so that we need only prove that the image of $\alpha$ in $\rmH^3(\pi_1^\et(Y,\overline{y})^{[p']},(\Q/\Z)'(2))$ vanishes along the residue map $\partial_{\mu,D}$. 
    
    Let $K\subset \overline{k}$ be the field of definition of $Y$. Again, \cite[Exposé V]{SGA1}, \emph{ibidem} shows that there exists an irreducible finite étale cover $g:Z\to Y$ together with a geometric point $\overline{z}$ such that $(Z,\overline{z})$ is characterised by the fact that $g$ is good and that $\pi_1^\et(Z,\overline{z})^{[p']}=D$. Thanks to the proof of \cite[Théorème]{ducrosdimensioncohomologique}, we may choose a field $L$ of cohomological dimension~$\le 2$ in which $K(Z)$ is algebraically closed and such that the restriction $\Br(K(Z))\to \Br(L)$ is injective. Since $K(Y)\subset L$, it follows that $Y(L)\neq \varnothing$ and we may choose an isomorphism $Y\otimes_K L\simeq (U\otimes_k L)/H$ for some adequate $L$-subgroup $H$ of $G\otimes_k L$ which in turn induces an isomorphism $\pi_1^\et(Y,\overline{y})^{[p']}\times_{\Gamma_K}\Gamma_{L}\simeq H\rtimes\Gamma_{L}$. As the $H$-torsor $U\otimes_k L\to (U\otimes_k L)/H$ is versal \cite[Example 5.4]{garibaldimerkurjevserre}, the map of pointed sets $Y(L((t))) \to \rmH^1(L((t)),H)$ is onto; see \emph{op.\ cit.}, Definition 5.1. Moreover, the natural map $\Gamma_{K(Z)}\to \pi_1^\et(Z,\overline{z})^{[p']}$ is surjective since $Z$ is normal \cite[Exposé I, Corollaire~10.2]{SGA1}. This yields a surjection $\Gamma_L\twoheadrightarrow D$, and the composite map
\[\phi_P:\Gamma_{L((t))}\longrightarrow \Gamma_\mathrm{t}\simeq\hat{\Z}'(1)\rtimes \Gamma_L\longrightarrow\hat{\Z}'(1)\rtimes D\longrightarrow \langle\Im(\mu),D\rangle=\pi_1^\et(Y,\overline{y})^{[p']}\]
defines a section of the projection map $\pi_1^\et(Y,\overline{y})^{[p']}\times_{\Gamma_K}\Gamma_{L((t))} \to \Gamma_{L((t))}$. This hence defines an element of $\rmH^1(L((t)),H)$ by virtue of \cite[§1.2, Proposition 8]{stixlivre}. Therefore, the map $\phi_P$ must come from an $L((t))$-point $P$ of $Y$ up to conjugation by an element of $H$. 

Since $\alpha\in M_\nr(Y)$, then Theorem~\ref{WittenbergH3nr} ensures that this class vanishes along the pullback $\phi_P^*:\rmH^3(\pi_1^\et(Y,\overline{y})^{[p']},(\Q/\Z)'(2))\to \rmH^3(\Gamma_{\text{t}},(\Q/\Z)'(2))$. By the functoriality of the Hochschild-Serre spectral sequence, we have a commutative diagram:
\[\begin{tikzcd}
	{\rmH^3(\pi_1^\et(Y,\overline{y})^{[p']},(\Q/\Z)'(2))} & {\rmH^3(\Gamma_{\text{t}},(\Q/\Z)'(2))} \\
	{\rmH^2(D,(\Q/\Z)'(1))} & {\rmH^2(L,(\Q/\Z)'(1)),}
	\arrow["\phi_P^*", from=1-1, to=1-2]
	\arrow["\partial_{\mu,D}"', from=1-1, to=2-1]
	\arrow[ from=1-2, to=2-2, "\wr"]
	\arrow[from=2-1, to=2-2]
\end{tikzcd}\]
where the bottom horizontal arrow coincides with the composite map
\[\rmH^2(D,(\Q/\Z)'(1))\overset{\sim}{\longrightarrow}\Br(Z)\{p'\}\longrightarrow\Br(K(Z))\longrightarrow\Br(L).\]
Note that the last map is injective by assumption. Since $Z$ is a smooth $K$-variety, the map $\Br(Z)\to \Br(K(Z))$ is injective as well; see for instance \cite[Theorem 3.5.5]{CTS21}. This forces $\partial_{\mu,D}(\alpha)=0$ and therefore implies the inclusion $M_\nr(X)\subset \rmH^3_\mathrm{gnr}(\pi_1^\et(X,\overline{x})^{[p']},(\Q/\Z)'(2))$, as desired.\end{proof}

\begin{rmk}\label{anyany} We can also show that the equality $M_\nr(X)=\rmH^3_{\mathrm{any,any}}(\pi_1^\et(X,\overline{x})^{[p']},(\Q/\Z)'(2))$ holds: indeed, the inclusion $\rmH^3_{\mathrm{any,any}}(\pi_1^\et(X,\overline{x})^{[p']},(\Q/\Z)'(2))\subset M_\nr(X)$ is evident in view of the inclusion $\rmH^3_{\mathrm{any,any}}(\pi_1^\et(X,\overline{x})^{[p']},(\Q/\Z)'(2))\subset \rmH^3_{\mathrm{any,pro-}p'}(\pi_1^\et(X,\overline{x})^{[p']},(\Q/\Z)'(2))$ and the equality in the statement of Proposition \ref{imageinverse}. As for the converse inclusion, it follows \emph{mutatis mutandis} from the second part of the proof of \emph{loc.\ cit.}\ applied to any pair $(\mu,D)\in {_\mathrm{any}\mathscr{D}_\mathrm{{\mathrm{any}}}}$.
\end{rmk}

\subsection{Cohomology classes trivial at the generic point} It now remains to precisely describe the kernel of the natural map $\rmH^3_\et(X,(\Q/\Z)'(2))\to \rmH^3(k(X),(\Q/\Z)'(2))$, that is, the subgroup of negligible classes in the sense of Serre; see \cite[§26, p.61]{garibaldimerkurjevserre}. Thanks to Proposition~\ref{suiteexactefondamentale}, we know that the latter group coincides with the image of the injective map $\CH^2(X)\to \rmH^3(\pi_1^\et(X,\overline{x})^{[p']},(\Q/\Z)'(2))$ which appears in \eqref{suitefondamentale}. We aim at precisely determining this image.

\begin{lem}\label{comparaisonclassesdepermutationCH2} Let $\pi:Y\to X$ be the finite good cover associated in
§\ref{definitionnegligible} with choices of a finite separable extension $k \subset K
\subset \overline{k}$ and of a subgroup $H \subset G \otimes_k K$. We have the following commutative diagram:
\[\begin{tikzcd}
	{\Pic(Y)^{\otimes 2}} & {\rmH^1(\pi_1^\et(Y,\overline{y})^{[p']},(\Q/\Z)'(1))^{\otimes 2}} \\
	{\mathrm{CH}^2(Y)} & {\rmH^3(\pi_1^\et(Y,\overline{y})^{[p']},(\Q/\Z)'(2)),}
	\arrow[from=1-1, to=1-2, hookrightarrow]
	\arrow[from=1-1, to=2-1]
	\arrow[from=1-2, to=2-2]
	\arrow[from=2-1, to=2-2, hookrightarrow]
\end{tikzcd}\]
where the top horizontal arrow is the inclusion coming from the choice of $y\in Y(K)$, the bottom horizontal one is the map appearing in \eqref{suitefondamentale}, and the leftmost (resp.\ rightmost) vertical map is  the intersection product (resp. the cup-product).
    \end{lem}
    \begin{proof}
In the derived category of complexes of abelian Zariski sheaves on $X$, the distinguished triangle $\Z(2)\to \Q(2)\to \Q/\Z(2)\to \Z(2)[1]$ induces an exact sequence:
\[\rmH^3_\mot(Y,\Q(2))\longrightarrow\rmH^3_\mot(Y,\Q/\Z(2))\longrightarrow \rmH^4_\mot(Y,\Z(2))\longrightarrow \rmH^4_\mot(Y,\Q(2)),\]
where $\rmH^4_\mot(Y,\Z(2))\simeq \CH^2(Y)$ because $Y$ is smooth and quasi-projective; see \cite[Théorème~7.5.2]{cthk}. As $\CH^2(Y)$ is torsion, the map $\CH^2(Y)\to \rmH^4_\mot(Y,\Q(2))$ in the above sequence is zero. By comparing the coniveau spectral sequences
    \[E_1^{p,q}(U_K)=\bigoplus_{u\in U_K^{(p)}}\rmH^{q-p}_\mot(\kappa(u),\Z(2-p))\Rightarrow\rmH^{p+q}_\mot(U_K,\Z(2))\]
    and 
    \[E_1^{p,q}(\A^m_K)=\bigoplus_{a\in {(\A^m_K)}^{(p)}}\rmH^{q-p}_\mot(\kappa(a),\Z(2-p))\Rightarrow\rmH^{p+q}_\mot(\A^m_K,\Z(2)),\]
we see that $\rmH^3_\mot(U_K,\Z(2))\simeq \rmH^3_\mot(\A^m_K,\Z(2))$, and the latter is isomorphic to $\rmH^3_\mot(K,\Z(2))$ by homotopy invariance \cite[Theorem 2.1]{bloch86}. Thanks to the norm-residue isomorphism \cite[Lemma 1.6, Theorem 1.7]{weibel}, we have $\rmH^3_\mot(K,\Z(2))=0$. In particular, a restriction-corestriction argument shows that $\rmH^3_\mot(Y,\Q(2))$ is killed by the order of $H$ and hence is trivial. Therefore, the connecting map $\rmH^3_\mot(Y,\Q/\Z(2))\to \CH^2(Y)$ is an isomorphism and we have $\rmH^3_\mot(Y,\Q/\Z(2))=\rmH^3_\mot(Y,(\Q/\Z)'(2))$. On the other hand, we have the following commutative diagram:
\[\begin{tikzcd}
	{\rmH^1_\mot(Y,(\Q/\Z)'(1))^{\otimes 2}} & {\rmH^1_\et(Y,(\Q/\Z)'(1))^{\otimes 2}} & {\rmH^1(\pi_1^\et(Y,\overline{y})^{[p']},(\Q/\Z)'(1))^{\otimes 2}} \\
	{\rmH^3_\mot(Y,(\Q/\Z)'(2))} & {\rmH^3_\et(Y,(\Q/\Z)'(2))} & {\rmH^3(\pi_1^\et(Y,\overline{y})^{[p']},(\Q/\Z)'(2)),}
	\arrow[from=1-1, to=1-2]
	\arrow[from=1-1, to=2-1]
	\arrow["\sim", from=1-2, to=1-3]
	\arrow[from=1-2, to=2-2]
	\arrow[from=1-3, to=2-3]
	\arrow[from=2-1, to=2-2]
	\arrow["\sim", from=2-2, to=2-3]
\end{tikzcd}\]
where the leftmost horizontal maps denote the étale realisation maps, the rightmost ones come from Lemma \ref{cohomologieequivariante}, and the vertical ones are induced by the cup-products in motivic and étale cohomology respectively. 

By applying \cite[Corollary 6.4]{bloch86}, we derive a short exact sequence
\[0\longrightarrow K^*\otimes_\Z(\Q/\Z)'\longrightarrow\rmH^1_\mot(Y,(\Q/\Z)'(1))\longrightarrow \Pic(Y)\longrightarrow 0\]
which, combined with the natural isomorphism $\rmH^1_\mot(Y,(\Q/\Z)'(1))\simeq \rmH^1_\et(Y,(\Q/\Z)'(1))$ due to \emph{op.\ cit., ibidem}, identifies with the short exact sequence coming from the Kummer sequence on $Y_\et$. Again, the choice of $y\in Y(K)$ provides a retraction of the leftmost map in the above sequence and hence an inclusion $\Pic(Y)\hookrightarrow \rmH^1_\mot(Y,(\Q/\Z)'(1))$, which is compatible with the isomorphism $\rmH^1_\et(Y,(\Q/\Z)'(1))\simeq \rmH^1(\pi_1^\et(Y,\overline{y})^{[p']},(\Q/\Z)'(1))$. The assertion hence follows from the commutativity of the following diagram:
\[\begin{tikzcd}
	 {\Pic(Y)^{\otimes 2}} & {\rmH^1_\mot(Y,(\Q/\Z)'(1))^{\otimes 2}} \\
	 {\CH^2(Y)} & {\rmH^3_\mot(Y,(\Q/\Z)'(2)),}
	\arrow[from=1-1, to=1-2, hookrightarrow]
	\arrow[from=1-1, to=2-1]
	\arrow[from=1-2, to=2-2]
	\arrow["\sim", from=2-1, to=2-2]
\end{tikzcd}\]
where the leftmost vertical map is induced by the intersection product and the bottom one is the inverse of the connecting map $\rmH^3_\mot(Y,\Q/\Z(2))\to \CH^2(Y)$. Indeed, to check the latter it suffices to note that the cup-product map $\rmH^1_\mot(Y,(\Q/\Z)'(1))^{\otimes 2}\to \rmH^3_\mot(Y,(\Q/\Z)'(2))$ and the intersection product $\Pic(Y)^{\otimes 2}\to \CH^2(Y)$ are by construction compatible with the evident connecting maps.\end{proof}

\begin{prop}\label{negligible} If $G$ is a constant $k$-group, we have an equality 
\[\mathrm{Im}\Big(\CH^2(X)\to \rmH^3(\pi_1^\et(X,\overline{x})^{[p']},(\Q/\Z)'(2))\Big)\{2'\}=\rmH^3_\mathrm{per}(\pi_1^\et(X,\overline{x})^{[p']},(\Q/\Z)'(2))\{2'\}\]
of subgroups of $\rmH^3(\pi_1^\et(X,\overline{x})^{[p']},(\Q/\Z)'(2))$.
\end{prop}

The proof of this assertion is a direct generalisation of that of \cite[Proposition~4.3]{peyreh3noether}. The main technical difference lies in the fact that the reduction argument involving monomial characters in \emph{loc.\ cit.}\ cannot work for arbitrary supersolvable groups over a non-closed field. However, it does work for $p$-groups of odd order over fields containing the $p$\textsuperscript{th}-roots of unity thanks to work of Kahn \cite{kahnrep}, which surprisingly suffices to describe the image of the map $CH^2(X)\to \rmH^3(\pi_1^\et(X,\overline{x})^{[p']},(\Q/\Z)'(2))$.

\begin{proof} Let $G_q$ be any $q$-Sylow subgroup of $G$ for some prime $q\neq 2$. If we let $Y:=U/G_q$, then the evident morphism $\pi:Y\to X$ is finite étale of degree coprime to $q$. Let $y\in Y(k)$ be the natural point lying above $x\in X(k)$; by the very definition of permutation negligible classes, we thus have an inclusion $\pi_*\rmH^3_\mathrm{per}(\pi_1^\et(Y,\overline{y})^{[p']},(\Q/\Z)'(2))\subset \rmH^3_\mathrm{per}(\pi_1^\et(X,\overline{x})^{[p']},(\Q/\Z)'(2))$. Moreover, we have a commutative diagram:
\[\small \begin{tikzcd}
	& {H^3(\pi_1^\et(X,\overline{x})^{[p']},(\Q/\Z)'(2))} & {\rmH^0(X,\mathscr{H}^3_\et((\Q/\Z)'(2)))} \\
	{\rmH^3_\mathrm{per}(\pi_1^\et(Y,\overline{y})^{[p']},(\Q/\Z)'(2))} & {H^3(\pi_1^\et(Y,\overline{y})^{[p']},(\Q/\Z)'(2))} & {\rmH^0(Y,\mathscr{H}^3_\et((\Q/\Z)'(2)))} \\
	{\rmH^3_\mathrm{per}(\pi_1^\et(X,\overline{x})^{[p']},(\Q/\Z)'(2))} & {H^3(\pi_1^\et(X,\overline{x})^{[p']},(\Q/\Z)'(2))} & {\rmH^0(X,\mathscr{H}^3_\et((\Q/\Z)'(2)))}.
	\arrow[from=1-2, to=1-3]
	\arrow["{\pi^*}"', from=1-2, to=2-2]
	\arrow["{\pi^*}"', from=1-3, to=2-3]
	\arrow[from=2-1, to=2-2]
	\arrow["{\pi_*}"', from=2-1, to=3-1]
	\arrow[from=2-2, to=2-3]
	\arrow["{\pi_*}"', from=2-2, to=3-2]
	\arrow["{\pi_*}"', from=2-3, to=3-3]
	\arrow[from=3-1, to=3-2]
	\arrow[from=3-2, to=3-3]
\end{tikzcd}\]
Since $\pi_*\pi^*$ is multiplication by $\deg(\pi)$, then assuming that the middle row is exact implies that the lower one is exact on its $q$-primary torsion subgroups. It therefore suffices to show that the kernel of $\rmH^3(\pi_1^\et(X,\overline{x})^{[p']},(\Q/\Z)'(2))\to\rmH^0(X,\mathscr{H}^3_\et((\Q/\Z)'(2)))$ is contained in $\rmH^3_\mathrm{per}(\pi_1^\et(X,\overline{x})^{[p']},(\Q/\Z)'(2))$ when $G$ is a $q$-group. Moreover since the degree of $k(\mmu_q)$ over $k$ is coprime to $q$, then by restriction-corestriction we may further assume that $k$ contains a primitive $q$\textsuperscript{th}-root of unity. 

As recalled earlier in this section, we know that $\mathrm{CH}^2(X)\simeq \CH^2(B_kG)$ consists of elements of the form $c_2(\rho)$ for $[\rho]\in \mathrm{Rep}_k(G)$. Moreover, for any two  $k$-representations $\rho$ and $\rho'$ of $G$, the Whitney formula \cite[Theorem 3.2]{fulton} reads
\[c_2(\rho+\rho')=c_2(\rho)+c_1(\rho)\smile c_1(\rho')+c_2(\rho')\in\rmH^3(\pi_1^\et(X,\overline{x})^{[p']},(\Q/\Z)'(2)).\]  
By Lemma \ref{comparaisonclassesdepermutationCH2} we have
 $c_1(\rho)\smile c_1(\rho')\in \rmH^3_\mathrm{per}(\pi_1^\et(X,\overline{x})^{[p']},(\Q/\Z)'(2))$,
so that the Chern class induces a group morphism:
\[\mathrm{Rep}_k(G)\longrightarrow \frac{\rmH^3(\pi_1^\et(X,\overline{x})^{[p']},(\Q/\Z)'(2))}{\rmH^3_{\mathrm{per}}(\pi_1^\et(X,\overline{x})^{[p']},(\Q/\Z)'(2))},\]
and it therefore remains to show that this map is identically zero under the assumption that $q\neq 2$. 

The group $\mathrm{Rep}_k(G)$ is by construction generated by classes of irreducible linear $k$-representations. As $G$ is a $q$-group and $\mmu_q\subset k$, then \cite[Propositions 7.4 and 9.1]{kahnrep} imply that each irreducible representation $\rho$ of $G$ is induced by a $1$-dimensional representation $\chi:H\to k^*$ for some adequate subgroup $H\subset G$ (i.e.\ a character of degree $1$). Let us write $Y:=U/H$ and denote by $\pi:Y\to X$ the associated étale cover. We prove by induction on the index $(G:H)$ that the Chern class $c_2(\pi_*\chi)$ is permutation negligible. If $(G:H)=1$, then $c_2(\pi_*\chi)=c_2(\chi)=0$ since $\chi$ is $1$-dimensional. If $(G:H)=q^r$ for some $r\ge 1$, assume that the desired property holds for subgroups of index strictly less than $q^r$. Thanks to \cite{fultonmacpherson}, see also \cite[§14.2]{NormsBachmannHoyois}, there exist higher transfer morphisms $\pi_*^{(n)}:\CH^r(Y)\to \CH^{nr}(X)$ for any $n\ge 1$ and any $r\ge0$, such that $\pi_*^{(1)}=\pi_*$ is the usual pushforward. By \cite[Corollary 5.3]{fultonmacpherson} we notably have, for any character $\chi$ of $H$ of degree~$1$, the following formula:
\[c_2(\pi_*\chi)=\pi_*c_2(\chi)+\pi_*^{(2)}c_1(\chi)+c_1(\pi_*\mathbf{1})\smile \pi_* c_1(\chi)+c_2(\pi_*\mathbf{1})\in \CH^2(X),\]
where $\mathbf{1}$ denotes the vector bundle associated with the trivial $1$-dimensional representation. Since $\chi$ is $1$-dimensional, then $c_2(\chi)=0$, and $c_1(\pi_*\mathbf{1})\smile \pi_* c_1(\chi)$ lies in the subgroup of permutation negligible classes in $\rmH^3(\pi_1^\et(X,\overline{x})^{[p']},(\Q/\Z)'(2))$. Moreover \cite[Corollary~10.9]{fultonmacpherson} yields an equality:
\[\pi_*\chi^2-(\pi_*\chi)^2+2\pi_*^{(2)}\chi=0.\]
This gives $\pi_*^{(2)}\chi=((\pi_*\chi)^2-\pi_*\chi^2)/2$, and we readily observe that the latter is an element of $\rmH^3_\mathrm{per}(\pi_1^\et(X,\overline{x})^{[p']},(\Q/\Z)'(2))$. 

We thus need only prove that $c_2(\pi_*\mathbf{1})$ is permutation negligible as well. By \cite[Chapter 2, Theorem 1.6 and Corollary]{suzuki}, there exists an intermediate subgroup $H\subset H'\subset G$ in which $H$ is normal and such that $(H':H)=q$. In particular, we can let $Y':=U/H'$ and $\pi':Y'\to X$, resp.\ $\varpi:Y\to Y'$ be the induced étale covers, so that $c_2(\pi_*\mathbf{1})=c_2(\pi'_*(\varpi_*\mathbf{1}))$. Since $\mmu_q\subset k$, we may choose a character $\lambda:H'\to \mmu_q$ such that $H=\Ker (\lambda)$. The class of $\varpi_*\mathbf{1}$ in $\mathrm{Rep}_k(H')$ then splits as $[\varpi_*\mathbf{1}]=[\mathbf{1}]+[\lambda]+\hdots+[\lambda]^{q-1}$. Applying the Whitney formula once again and using the induction hypothesis, we then see that $c_2(\pi_*\mathbf{1})=c_2(\pi'_*\mathbf{1}+\pi_*'\lambda+\hdots+\pi'_*\lambda^{q-1})$ is a sum of permutation negligible classes, and the desired assertion follows. 
\end{proof}

\begin{proof}[Proof of Theorem \ref{theoremeprincipal}] Combine Proposition \ref{suiteexactefondamentale}, Proposition \ref{imageinverse} and Proposition \ref{negligible}.
\end{proof}
\section{Algebraic classes}
Let $k$ be a field of exponential characteristic $p>0$, let $G$ be a finite $k$-group of order coprime to $p$, let $c\ge 3$, fix a classifying pair $(U\subset \A^m_k,c)$, and denote $X:=U/G$. We now move our attention to the subgroup of algebraic classes within $\rmH^3_\mathrm{nr}(k(X),(\Q/\Z)'(2))$, that is, the group: \[\rmH^3_{\nr,\mathrm{alg}}(k(X),(\Q/\Z)'(2)):=\Ker\Big(\rmH^3_\nr(k(X),(\Q/\Z)'(2))\to \rmH^3_\nr(\overline{k}(X),(\Q/\Z)'(2))\Big).\]

\subsection{An analogue of the algebraic unramified Brauer group} In this section, we define a subgroup of the kernel of the base-change map $\rmH^3_\et(X,(\Q/\Z)'(2))\to \rmH^3_\et(\overline{X},(\Q/\Z)'(2))$ whose image along the natural map $\rmH^3_\et(X,(\Q/\Z)'(2))\to\rmH^3(k(X),(\Q/\Z)'(2))$ is contained in the unramified subgroup $\rmH^3_\nr(k(X),(\Q/\Z)'(2))$. We then provide various refinements and vanishing results for this subgroup. 
\subsubsection{A distinguished map}
\label{flèchedistinguée}Consider the Hochschild-Serre spectral sequence 
\[E_2^{p,q}=\rmH^p(k,\rmH^q_\et(\overline{X},(\Q/\Z)'(2)))\Rightarrow\rmH^{p+q}_\et(X,(\Q/\Z)'(2)),\]
with associated filtration $F^\bullet\rmH^*_\et(X,(\Q/\Z)'(2))$. For $0\le r\le 3$, we denote \[F^r\rmH^3_\mathrm{gnr}(X,(\Q/\Z)'(2)):=F^r\rmH^3_\et(X,(\Q/\Z)'(2))\cap \rmH^3_\mathrm{gnr}(X,(\Q/\Z)'(2)),\]
where $\rmH^3_\mathrm{gnr}(X,(\Q/\Z)'(2))\subset \rmH^3_\et(X,(\Q/\Z)'(2))$ is the subgroup of globally unramified classes, see Definition \ref{definitionglobalementnr}.
Thanks to Lemma \ref{cohomologieequivariante}, there exists a natural isomorphism of $\Gamma_k$\nobreakdash-modules $\rmH^q_\et(\overline{X},(\Q/\Z)'(2))\simeq \rmH^q(G(\overline{k}),(\Q/\Z)'(2))$ for $0\le q\le 2c-2$. For simplicity, we will from now on denote $\widehat{G}^\ab(1):=\rmH^1(G(\overline{k}),(\Q/\Z)'(2))\simeq\Hom(G(\overline{k}),(\Q/\Z)'(2))$ and $\rmH^2(G,(\Q/\Z)'(2)):=\rmH^2(G(\overline{k}),(\Q/\Z)'(2))$, equipped with their $\Gamma_k$\nobreakdash-action. We also have an equality $\rmH^0(\overline{X},(\Q/\Z)'(2))=(\Q/\Z)'(2)$. In particular, the pullback along the structural morphism $X\to \Spec(k)$ coincides with the edge map $E_2^{3,0}=\rmH^3(k,(\Q/\Z)'(2))\to H^3_\et(X,(\Q/\Z)'(2))$. For any $x\in X(k)$, the associated section $\delta_x:\Gamma_k\to \pi_1^\et(X,\overline{x})$ yields a retraction of this map, and we therefore obtain an inclusion:
\[E_\infty^{2,1}=\Coker\Big(\rmH^0(k,\rmH^2(G,(\Q/\Z)'(2)))\to \rmH^2(k,\widehat{G}^\ab(1))\Big)\hookrightarrow \rmH^{3,\norm}_\et(X,(\Q/\Z)'(2)).\]
Composing this map with the projection $E_2^{2,1}\to E_{\infty}^{2,1}$ thus provides a map 
\begin{equation}\label{flècheE21versH3}\rmH^2(k,\widehat{G}^\ab(1))\longrightarrow \rmH^{3,\norm}_\et(X,(\Q/\Z)'(2)),\end{equation} which only depends on the choice of $x\in X(k)$.

\subsubsection{A criterion for algebraic classes}\label{definitionaccouplementcupproduit} Let $\mathrm{exp(G)}$ denote the exponent of $G(\overline{k})$ and consider the action of $(\Z/\exp(G))^*$ on the set $G(\overline{k})/(\text{conjugacy})$ of conjugacy classes of $G(\overline{k})$ via $[\tau]\mapsto [\tau^i]$ for $i\in (\Z/\mathrm{exp}(G))^*$. The \emph{outer exponent} of $G(\overline{k})$ is the smallest divisor $e$ of $\mathrm{exp}(G)$ such that this action factors through the quotient map $(\Z/\mathrm{exp}(G))^* \to (\Z/e)^*$. In particular, the cyclotomic action of $\Gamma_k$ on $G(\overline{k})/(\text{conjugacy})$ factors through the natural morphism $\chi: \Gamma_k \to(\Z/e)^*$.

Fix a pair $(\mu,D)\in  {_\mathrm{any}\mathscr{D}_\mathrm{{\mathrm{pro-}p'}}}$ and let $L\subset \overline{k}$ be the subfield fixed by the image of $D$ in $\Gamma_k$. Consider the composite map $\rho:D\hookrightarrow G(\overline{k})\rtimes \Gamma_k\to \Gamma_k$ and let $\sigma \in\Gamma_L$. Since $\mu:\widehat{\Z}'(1)\to G(\overline{k})$ satisfies the equality ${^\sigma[\tau]} = [\tau^{\chi(\sigma)}]$ in $G(\overline{k})/(\text{conjugacy})$ for any $\tau\in \Im(\mu)$, then its
image $\overline{\mu}$ in $G^\ab(-1):=\Hom(\widehat{\Z}'(1),G(\overline{k})^\ab)$ is $\sigma$-invariant. We may then view $\overline{\mu}$ as an element of $\rmH^0(L,G^\ab(-1))=\rmH^0(D,G^\ab(-1))$, and we obtain a composite map
\[\begin{tikzcd}
\rmH^2(k,\widehat{G}^{\mathrm{ab}}(1)) \ar[r, "\rho^*"] & \rmH^2(D,\widehat{G}^{\mathrm{ab}}(1))  \ar[r, "\smile\overline{\mu}"] & \rmH^2(D,(\Q/\Z)'(1)). \end{tikzcd}\]
\begin{prop}\label{nonramificationalgébrique}
Consider the subgroup $\Sha^2_\mathrm{nr}(k,\widehat{G}^\mathrm{ab}(1))$ of $\rmH^2(k,\widehat{G}^\mathrm{ab}(1))$ consisting of those classes $\alpha$ such that for every pair $(\mu,D)\in  {_\mathrm{any}\mathscr{D}_\mathrm{{\mathrm{pro-}p'}}}$ with associated morphism $\rho:D\to\Gamma_k$, we have $\rho^*\alpha\smile\overline{\mu}=0\in\rmH^2(D,(\Q/\Z)'(1))$. Consider the composite map 
\begin{equation}\label{flècheH2versH0H3}\rmH^2(k,\widehat{G}^\ab(1))\longrightarrow \rmH^3_\et(X,(\Q/\Z)'(2))\longrightarrow \rmH^0(X,\mathscr{H}^3_\et((\Q/\Z)'(2))),\end{equation}
where the first map is \eqref{flècheE21versH3}. Then a class $\alpha\in \rmH^2(k,\widehat{G}^\ab(1))$ lies in $\Sha^2_\nr(k,\widehat{G}^\ab(1))$ if and only if its image with respect to the above map lies in $\rmH^{3,\norm}_\mathrm{nr}(k(X),(\Q/\Z)'(2)).$ More precisely, the choice of $x\in X(k)$ provides the following splitting: 
   \[F^2\rmH^3_\mathrm{gnr}(X,(\Q/\Z)'(2))\simeq\rmH^3(k,(\Q/\Z)'(2))\oplus\Coker\Big(\rmH^0(k,\rmH^2(G,(\Q/\Z)'(2)))\to \Sha^2_\nr(k,\widehat{G}^\ab(1))\Big).\]

\end{prop}
\begin{proof} By virtue of Proposition \ref{imageinverse}, we need only prove that a class $\alpha\in \rmH^2(k,\widehat{G}^\ab(1))$ lies in the subgroup $\Sha^2_\nr(k,\widehat{G}^\ab(1))$ if and only if its image in $\rmH^3(\pi_1^\et(X,\overline{x}),(\Q/\Z)'(2))$ along \eqref{flècheE21versH3} lies in $\rmH^3_\mathrm{gnr}(\pi_1^\et(X,\overline{x}),(\Q/\Z)'(2))$. Let $\alpha\in \rmH^2(k,\widehat{G}^\ab(1))$ be a class whose image $\beta$ in $\rmH^3(\pi_1^\et(X,\overline{x}),(\Q/\Z)'(2))$ is not globally unramified. Consider a pair $(\mu,D)\in {_\mathrm{any}\mathscr{D}_\mathrm{{\mathrm{pro-}p'}}}$ such that $\partial_{\mu,D}(\beta)\neq 0\in \rmH^2(D,(\Q/\Z)'(1))$. By virtue of the Hochschild-Serre spectral sequences for the semi-direct product decompositions of $\widehat{\Z}'(1)\rtimes D$ and $\pi_1^\et(X,\overline{x})^{[p']}\simeq G(\overline{k})\rtimes\Gamma_k$ respectively, we derive a commutative diagram:

\begin{equation}\label{diagrammedefinitionsha2nr}\begin{tikzcd}
	{F^2\rmH^3(\pi_1^\et(X,\overline{x}),(\Q/\Z)'(2))} & {\rmH^3(\widehat{\Z}'(1)\rtimes D ,(\Q/\Z)'(2))} \\
	{\rmH^2(k,\widehat{G}^\ab(1))} & {\rmH^2(D,(\Q/\Z)'(1)).}
	\arrow[from=1-1, to=1-2]
	\arrow[from=1-1, to=2-1]
	\arrow["{\partial_{\mu,D}}"{description}, from=1-1, to=2-2]
	\arrow[from=1-2, to=2-2]
	\arrow["{\eqref{flècheE21versH3}}", shift left=3, from=2-1, to=1-1]
	\arrow["{\alpha\mapsto \rho^*\alpha\smile\overline{\mu}}"', from=2-1, to=2-2]
\end{tikzcd}\end{equation}
Note that the image of an element $\alpha\in\Coker\big(\rmH^0(k,\rmH^2(G,(\Q/\Z)'(2)))\to\rmH^2(k,\widehat{G}^\ab(1))\big)$ in $\rmH^2(D,(\Q/\Z)'(1))$ does not depend on any lift in $\rmH^2(k,\widehat{G}^\ab(1))$. Indeed, we have the following commutative diagram:
\[\begin{tikzcd}
	{\rmH^0(k,\rmH^2(G,(\Q/\Z)'(2)))} & {\rmH^0\Big(D,\rmH^2(\widehat{\Z}'(1),(\Q/\Z)'(2))\Big)} \\
	{\rmH^2(k,\widehat{G}^\ab(1))} & {\rmH^2(D,(\Q/\Z)'(1)),}
	\arrow[from=1-1, to=1-2]
	\arrow[from=1-1, to=2-1]
	\arrow[from=2-1, to=2-2]
	\arrow[from=1-2, to=2-2]
\end{tikzcd}\]
where the horizontal arrows are induced by $\phi_2^*$ and $\rmH^0\big(D,\rmH^2(\widehat{\Z}'(1),(\Q/\Z)'(2))\big)$ is trivial since the group $\rmH^2(\widehat{\Z}'(1),(\Q/\Z)'(2))$ is itself zero (as $\mathrm{cd}(\widehat{\Z}'(1))=1$). Since $\beta$ does not vanish in $\rmH^2(D,(\Q/\Z)'(1))$ by assumption, we must have $\rho^*\alpha\smile\overline{\mu}\neq 0\in \rmH^2(D,(\Q/\Z)'(1))$ as well. This proves the first assertion of the proposition. 

Conversely, let $\alpha\in \rmH^2(k,\widehat{G}^\ab(1))$ and consider its image $\beta$ in $\rmH^3(\pi_1^\et(X,\overline{x})^{[p']},(\Q/\Z)'(2))$ through \eqref{flècheE21versH3}. Assume that this class lies in $\rmH^3_\mathrm{gnr}(\pi_1^\et(X,\overline{x})^{[p']},(\Q/\Z)'(2))$. Fix a pair $(\mu,D)\in {_\mathrm{any}\mathscr{D}_\mathrm{{\mathrm{pro-}p'}}}$. By assumption, we have $\partial_{\mu,D}(\beta)=0$. The commutativity of \eqref{diagrammedefinitionsha2nr} readily implies that $\rho^*\alpha\smile \overline{\mu}=0\in \rmH^2(D,(\Q/\Z)'(1))$, as desired.\end{proof}

\begin{rmk}\label{remarquenonramificationalgébrique} (i) Consider the subgroup $\Sha^2_\mathrm{snr}(k,\widehat{G}^\mathrm{ab}(1))$ of $\rmH^2(k,\widehat{G}^\mathrm{ab}(1))$ consisting of those classes $\alpha$ such that for every  field $L\subset \overline{k}$ containing $k$ with $\Gamma_L$ pro-$\ell$ for some $\ell\neq p$ and with associated morphism $\rho:\Gamma_L\to\Gamma_k$ and every $\mu\in \Hom(\widehat{\Z}'(1),G(\overline{k}))$, if ${^\sigma[\tau]}=[\tau^{\chi(\sigma)}]$ in $G(\overline{k})/(\mathrm{conjugacy})$ for all $\sigma\in \Gamma_L$ and $\tau\in \Im(\mu)$, then $\rho^*\alpha\smile\overline{\mu}=0\in\rmH^2(L,(\Q/\Z)'(1))$. The first half of the proof of Proposition \ref{nonramificationalgébrique} shows that $\Sha^2_\mathrm{snr}(k,\widehat{G}^\mathrm{ab}(1))$ is always contained in the group $\Sha^2_\mathrm{nr}(k,\widehat{G}^\mathrm{ab}(1))$. However, these groups differ in general: see Lemma \ref{Sha2Q_16} and Remark \ref{remarqueSha2snrQ_16} for details.

(ii) Suppose that $G$ is an abelian $k$-group. Then the groups $\Sha^2_\mathrm{snr}(k,\widehat{G}^\mathrm{ab}(1))$ and $\Sha^2_{\mathrm{nr}}(k,\widehat{G}(1))$ coincide. Indeed, if $L\subset \overline{k}$ is a field containing $k$ with $\Gamma_L$ pro-$\ell$ and $\mu:\widehat{\Z}'(1)\to G(\overline{k})$ is a morphism, the condition that $^\sigma[\tau]=[\tau^{\chi(\sigma)}]$ in $G(\overline{k})/(\text{conjugacy})$ is equivalent to the equality $^\sigma\tau=\tau^{\chi(\sigma)}\in G(\overline{k})$  for all $\sigma\in \Gamma_L$ and $\tau\in \Im(\mu)$. If we let $D$ be the image of $\Gamma_L$ through $\rho$ and the section $\Gamma_k\to G(\overline{k})\rtimes \Gamma_k$, we thus observe that the datum of $\mu$ and $D$ defines a pair $(\mu,D)\in {_\mathrm{any}\mathscr{D}_\mathrm{{\mathrm{pro-}p'}}}$. We hence have the inclusion $\Sha^2_\mathrm{nr}(k,\widehat{G}^\mathrm{ab}(1))\subset \Sha^2_\mathrm{snr}(k,\widehat{G}^\mathrm{ab}(1))$.

(iii) If $G$ is an abelian $k$-group, then the group $\Sha^2_\mathrm{snr}(k,\widehat{G}^\mathrm{ab}(1))$ coincides with the group of those $\alpha\in \rmH^2(k,\widehat{G}^\ab(1))$ such that for every field $L\subset\overline{k}$ containing $k$ with $\Gamma_L$ pro-$\ell$ and with associated morphism $\rho:\Gamma_L\to\Gamma_k$, the pullback $\rho^*\alpha\in \rmH^2(L,\widehat{G}(1))$ lies in the left kernel of the cup-product pairing $\rmH^2(L,\widehat{G}(1))\times \rmH^0(L,G(-1))\to \rmH^2(L,(\Q/\Z)'(1))\subset \Q/\Z.$ Indeed, the quotient map $G(\overline{k})/(\mathrm{conjugacy})\to G(\overline{k})^\ab=G(\overline{k})$ is in this case bijective and $G(\overline{k})/(\mathrm{conjugacy})$ is an abelian group. Hence, for any such field $L\subset\overline{k}$ containing $k$, the map $\rmH^0(L,G(\overline{k})/(\mathrm{conjugacy})(-1))\to \rmH^0(L,G(-1))$ is bijective.
\end{rmk}

\begin{defn}
    Let $M$ be a torsion $\Gamma_k$-module. We let $\Sha^2_{\mathrm{alg}}(k,M)$ be the subgroup of $\rmH^2(k,M)$ consisting of those classes $\alpha$ such that for every field $L\subset\overline{k}$ containing $k$ with $\Gamma_L$ pro-$\ell$ for some $\ell\neq p$ and with associated morphism $\rho:\Gamma_L\to\Gamma_k$, the pullback $\rho^*\alpha\in \rmH^2(L,M)$ lies in the left kernel of the cup-product pairing $\rmH^2(L,M)\times \rmH^0(L,\widehat{M})\to \rmH^2(L,(\Q/\Z)'(1))\subset \Q/\Z.$

\end{defn}

In particular, Remark \ref{remarquenonramificationalgébrique} shows that if we let $M=\widehat{G}^\ab(1)$, we have an inclusion $\Sha^2_{\mathrm{alg}}(k,\widehat{G}^\ab(1))\subset \Sha^2_{\mathrm{nr}}(k,\widehat{G}^\ab(1))$.
\begin{prop}\label{annulation1}
    If $G$ is constant and $\mmu_{\mathrm{exp}(G)}\subset k$, then $\Sha^2_\mathrm{nr}(k,\widehat{G}^\mathrm{ab}(1))=0.$
\end{prop}
\begin{proof}
    Given the assumption, the cyclotomic character $\chi:\Gamma_k\to (\Z/\exp(G))^*$ is trivial. As $G$ is constant, we see that the condition ${^\sigma[\tau]}=[\tau^{\chi(\sigma)}]\in G(\overline{k})/\text{(conjugacy)}$ in the definition of $\Sha^2_\mathrm{nr}(k,\widehat{G}^\mathrm{ab}(1))$ is vacuous and always satisfied. Hence, the same line of arguments as in Remark \ref{remarquenonramificationalgébrique} shows that we have an equality $\Sha^2_\nr(k,\widehat{G}^\ab(1))=\Sha^2_\mathrm{snr}(k,\widehat{G}^\ab(1))$. On the other hand, we have isomorphisms of $\Gamma_k$\nobreakdash-modules $\widehat{G}^\mathrm{ab}(1)\simeq\Hom(G,\Z/\exp(G))\simeq G^\ab$. Because $\Gamma_k$ acts trivially on $G$, we may assume that the latter is an elementary abelian $\ell$\nobreakdash-group, i.e.\ $G=\Z/\ell^n$ for some integer $n\ge 0$. Coupled with the fact that $\Z/\ell^n\simeq \mmu_{\ell^n}$ as $\Gamma_k$-modules, it follows that for any field $L\subset\overline{k}$ containing $k$ with $\Gamma_L$ pro-$\ell$, we have $\rmH^2(L,\widehat{G}^\ab(1))=\rmH^2(L,\mmu_{\ell^n})=\Br(L)[\ell^n]$. Since $G$ is abelian, Remark \ref{remarquenonramificationalgébrique}, (iii) further shows that $\Sha^2_\mathrm{nr}(k,\widehat{G}^\mathrm{ab}(1))=\Sha^2_\mathrm{alg}(k,\mmu_{\ell^n})$. 
    
 As the restriction map $\rmH^2(k,\mmu_{\ell^n})\to \rmH^2(L,\mmu_{\ell^n})$ coincides with $\Br(k)[\ell^n]\to \Br(L)[\ell^n]$ for any field $L$ containing $k$ and the cup-product pairing $\Br(L)[p^n]\times \rmH^0(L,\mmu_{\ell^n})\to \Br(L)[\ell^n]$ is non-degenerate on the left, we infer that an element of $\Sha^2_\mathrm{alg}(k,\mmu_{\ell^n})$ must vanish in $\rmH^2(L,\mmu_{\ell^n})$ for any such $L$. If a class $\alpha\in \Br(k)[\ell^n]$ is non-zero, then \cite[Théorème]{ducrosdimensioncohomologique} implies that there exists a field $L\subset\overline{k}$ containing $k$ with $\Gamma_L$ of cohomological dimension $2$ and such that the image of $\alpha$ in $\Br(L)$ is non-zero as well. Up to replacing $\Gamma_L$ by a pro-$\ell$-Sylow subgroup, this hence forces $\Sha^2_{\mathrm{alg}}(k,\mmu_{\ell^n})=0$, as desired.
\end{proof}

Let us recall that a finite group is said to be \emph{metacyclic} if all of its Sylow subgroups are cyclic. Similarly, we say that a finite Galois extension is metacyclic if its Galois group is. A $k$-torus is said to be \emph{quasi-trivial} if it is isomorphic to $\mathrm{Res}_{A/k}(\Gm)$ for some étale $k$-algebra $A$, where $\mathrm{Res}_{A/k}$ denotes Weil's restriction of scalars.

\begin{prop}\label{annulation2} If $M$ is a finite $\Gamma_k$-module of order coprime to $p$ which is split by a metacyclic extension, then  $\Sha^2_\mathrm{alg}(k,M)=~0$.
\end{prop}
\begin{proof} Let $k\subset K$ be such a metacyclic extension. Taking a coflasque resolution of $M$ in the sense of \cite[Lemma 0.6]{colliotthelenesansucflasquetori} provides an exact sequence of finitely generated $\Gal(K/k)$-modules
\begin{equation}\label{resolutioncoflasque}0\longrightarrow \widehat{T}\longrightarrow \widehat{Q}\longrightarrow \widehat{M}\longrightarrow 0\end{equation}
where $\widehat{T}$ and $\widehat{Q}$ are free as abstract groups, $\widehat{Q}$ is a permutation $\Gal(K/k)$-module, and $\widehat{T}$ is coflasque, i.e.\ $\rmH^1(H,\widehat{T})=0$ for any subgroup $H\subset \Gal(K/k)$. As $\Gal(K/k)$ is a metacyclic group, a theorem of Endo-Miyata \cite[Theorem 1.5]{endomiyata} asserts that $\widehat{T}$ is a direct factor of a permutation $\Gal(K/k)$-module. Dualizing \eqref{resolutioncoflasque} provides a short exact sequence of groups of multiplicative type
\begin{equation}\label{resolutionflasque}1\longrightarrow M\longrightarrow Q\longrightarrow T\longrightarrow 1,\end{equation}
where $Q$ is a quasi-trivial torus and $T$ is a direct factor of a quasi-trivial torus. As a consequence, we have $\rmH^1(L,T)=0$ for any field extension $k\subset L$. Taking the long exact cohomology sequence associated to \eqref{resolutionflasque} thus yields an inclusion $\rmH^2(L,M)\hookrightarrow \rmH^2(L,Q)$ for any field extension $k\subset L$. 

It remains to show that the image of $\Sha^2_\mathrm{alg}(k,M)\subset \rmH^2(k,M)$ in $\rmH^2(k,Q)$ is trivial. Since $Q$ is quasi-trivial, then by Shapiro's lemma we need only treat the case $Q=\Gm$, for which $\rmH^2(k,Q)=\Br(k)$, and similarly for any field extension $k\subset L$. By the functoriality of the cup-product pairing, the image of $\Sha^2_\mathrm{alg}(k,M)$ in $\rmH^2(k,\Gm)$ must be contained in $\Sha^2_\mathrm{alg}(k,(\Q/\Z)'(1))$. The conclusion follows \emph{mutatis mutandis} as in the proof of Proposition \ref{annulation1}.
\end{proof}

\begin{prop} Let $G$ be a finite constant abelian group whose exponent is not divisible by $8$. Then $\Sha^2_\mathrm{nr}(k,\widehat{G}(1))=0$.
\end{prop}
\begin{proof} We may assume without loss of generality that $G=\Z/m$ where $m=2,4$ or $p^n$ with $p\neq 2$ and $n\ge 1$. In all of these cases, $k(\mmu_m)/k$ is cyclic and splits $M:=\widehat{G}(1)$, so that Remark \ref{remarquenonramificationalgébrique}, (ii) and (iii) and Proposition \ref{annulation2} directly imply the desired vanishing.
\end{proof}

\begin{rmk} This vanishing is quite sharp, as Theorem \ref{contrexempleZ/8} suggests.
\end{rmk}

\section{Universal triviality of unramified cohomology} For smooth and projective varieties $X$ over an algebraically closed field $k$ of exponential characteristic $p>0$, a question of prime interest is whether the unramified cohomology group $\rmH^3_\nr(k(X),(\Q/\Z)'(2))$ is \emph{universally trivial} as soon as it is trivial, i.e., for any field $L$ containing $k$, we have an isomorphism $\rmH^3(L,(\Q/\Z)'(2))\simeq \rmH^3_\nr(L(X),(\Q/\Z)'(2))$. This indeed leads to finer obstructions to rationality for such varieties, as there exist complex varieties with trivial unramified cohomology in degree three which is not universally trivial, as was highlighted by \cite{voisindegen}; see also \cite{voisinuniversel}. On the other hand, other varieties, e.g.\ cubic fourfolds containing a plane \cite{auelcolliotparimala}, have been shown to have universally trivial unramified cohomology in degree at most $3$. In view of the counter-examples provided by Peyre in \cite{peyreh3noether}, it is natural to ask whether the unramified cohomology of homogeneous spaces with finite stabilisers is universally trivial as soon as it is trivial, at least in degrees at most $3$. 

Let $G$ be a finite group of order coprime to $p$ viewed as a constant $k$-group. We recall that the \emph{Bogomolov multiplier} of $G$ is defined as the subgroup $B_0(G)$ of $\rmH^2(G,\Q/\Z)$ consisting of those classes which vanish along the restriction $\rmH^2(G,\Q/\Z)\to \rmH^2(B,\Q/\Z)$ for any bicyclic (or equivalently, abelian) subgroup $B\subset G$. Let $c\ge3$ and fix a classifying pair $(U\subset \A^m_k,c)$. If we denote $X:=U/G$, then we know from the classical work of Bogomolov that $\Brnr(\overline{X})\simeq B_0(G)$, see \cite[Lemma 3.10]{bogomolovbrauer}. In this section, we show that under the assumption that $B_0(G)$ vanishes, the algebraic subgroup of $\rmH^3_\nr(k(X),(\Q/\Z)'(2))$ eventually stabilises over an adequate finite field extension of $k$, i.e.\ it remains unchanged under further field extensions.

\subsection{Algebraic classes over a field with many roots of unity} We actually provide a much more general statement, which is the precise determination of the algebraic subgroup of the unramified cohomology in degree $3$ over any field which contains all the roots of unity that divide the order of the group $G$.
\begin{thm}\label{trivialitépartiealgébrique}
    Assume that $\mmu_{|G|}\subset k$. Then there exists a subgroup $H_{k,G}\subset \Hom(\Gamma_k,B_0(G))$ and a short exact sequence:
    \[ 0\longrightarrow \Ker\Big(\CH^2(X)\to\CH^2(\overline{X})\Big)\longrightarrow H_{k,G}\longrightarrow\rmH^{3,\norm}_{\nr,\mathrm{alg}}(k(X),(\Q/\Z)'(2))\longrightarrow 0.\]
    If we further assume that $\mmu_{|G|^2}\subset k$, then $H_{k,G}=\Hom(\Gamma_k,B_0(G))$.
\end{thm}
\begin{rmk} Here are two cases where the exact sequence of Theorem \ref{trivialitépartiealgébrique} simplifies.

(i) If $k$  admits no finite field extension of degree dividing the order of $B_0(G)$, then the group $\Hom(\Gamma_k,B_0(G))$ is zero. Theorem \ref{trivialitépartiealgébrique} hence implies that $\rmH^{3,\norm}_{\nr,\mathrm{alg}}(k(X),(\Q/\Z)'(2))$ is trivial under the further assumption that $\mmu_{|G|}\subset k$.

(ii) If $k$ contains the algebraic closure of its prime subfield, then the restriction map $\CH^2(X)\to \CH^2(\overline{X})$ is an isomorphism; see the proof of \cite[Lemma~17.2]{totaro}. Thus, the group $\Ker\big(\CH^2(X)\to\CH^2(\overline{X})\big)$ is trivial.

\end{rmk}
\begin{cor}\label{trivialitéuniverselle} If $\mmu_{|G|}\subset k$ and the groups $B_0(G)$ and $\rmH^3_\nr(\overline{k}(X),(\Q/\Z)'(2))$ both vanish, then the group $\rmH^3_\nr(k(X),(\Q/\Z)'(2))$ is universally trivial, i.e.\ for every field $L$ containing $k$, the natural map \[\rmH^3(L,(\Q/\Z)'(2))\to \rmH^3_\nr(L(X),(\Q/\Z)'(2))\]
is an isomorphism.
\end{cor}
\begin{proof} Combine Theorem \ref{trivialitépartiealgébrique} with the invariance of unramified cohomology under extensions of separably closed fields \cite[Theorem 4.4.1]{ctgersten}.
\end{proof}

\begin{rmk}\label{remarquebogomolov} The assumption that $B_0(G)$ be trivial is quite sharp in the following sense. For simplicity, let us assume that $k$ has characteristic $0$. If we indeed further assume that $\mmu_{|G|^2}\subset k $, then by \cite[Proposition 5.3]{lagardeunramifiedbrauer} the natural restriction map $\Br^\norm_\nr(X)\to \Br_\nr(\overline{X})\simeq B_0(G)$ is an isomorphism. In particular, if we let $X\hookrightarrow X^c$ be any smooth compactification, we have an equality $\Br(\overline{X}^c)=\rmH^0(k,\Br(\overline{X}^c))$. We may thus take the cup-product of a class $\alpha\in \Br(\overline{X}^c)$ with an element of $\rmH^1(k,\widehat{\Z}(1))$ and obtain a class $\beta\in \rmH^1(k,\Br(\overline{X}^c)\otimes\widehat{\Z}(1))$. 

Alternatively, if we descend this class $\alpha$ to an element of $\Br(X^c)$, lift it to $\rmH^2_\et(X^c,\Q/\Z(1))$ and take the cup-product with the same element of $\rmH^1(k,\widehat{\Z}(1))$, then we obtain a cohomology class $\widetilde{\beta}\in\rmH^3_\et(X^c,\Q/\Z(2))$ which dies in $\rmH^3_\et(\overline{X}^c,\Q/\Z(2))$ and whose image in $\rmH^1(k,\Br(\overline{X}^c)\otimes_\Z\widehat{\Z}(1))$ is $\beta$. In the case where $\overline{\Q}\subset k$, the map $\CH^2(X)\to \CH^2(\overline{X})$ is an isomorphism \cite[Lemma~17.2]{totaro} and Proposition \ref{suiteexactefondamentale} implies that the image of $\widetilde{\beta}$ in $\rmH^3_\nr(k(X),(\Q/\Z)'(2))$ cannot be zero as soon as $\widetilde{\beta}$ itself is non-zero. For instance, if $k=\C((t))$, then any non-zero element of $\rmH^1(k,\Br_\nr(\overline{X})\otimes_\Z\widehat{\Z}(1))\simeq \Br_\nr(\overline{X})$ provides a non-zero class in $\rmH^3_\nr(k(X),(\Q/\Z)'(2))$.
\end{rmk}

\subsection{Proof of Theorem \ref{trivialitépartiealgébrique}} We now give a proof of the theorem announced at the beginning of this section:
\begin{proof}[Proof of Theorem \ref{trivialitépartiealgébrique}] Since $G$ is constant, the group $\pi_1^\et(X,\overline{x})^{[p']}$ is isomorphic to the direct product $G\times \Gamma_k$; given the assumption on $k$, the cyclotomic character $\chi:\Gamma_k\to (\Z/\exp(G))^*$ is trivial. In particular, in the criterion from Proposition~\ref{imageinverse}, we need only consider pairs $(\mu,D)$ where $\mu:\widehat{\Z}'(1)\to G$ is a morphism, $D\subset \pi_1^\et(X,\overline{x})^{[p']}$, and the image of $\mu$ commutes with any element of $D$. We have an inclusion $D\subset \langle D,\Gamma_k\rangle$, and the latter group is of the form $H\times\Gamma_k$ for some subgroup $H\subset G$. If we let $C:=\Im(\mu)$, then since $C$ commutes with both $D$ and $\Gamma_k$, it commutes with $\langle D,\Gamma_k\rangle=H\times \Gamma_k$ as well. On the other hand, the residue $\partial_{\mu,D}$ factors through $\partial_{\mu,\langle D,\Gamma_k\rangle}$. By virtue of Remark \ref{anyany}, we may hence further restrict to those pairs $(\mu,D)$ such that $D=H\times \Gamma_k$ is a direct product. 

By the very construction of the residue $\partial_{\mu,D}$ (see~§\ref{paragraphegéometriquementnonramifié} for details), we may factor the associated map $\phi:\widehat{\Z}'(1)\rtimes D\to G\times \Gamma_k$ as follows:

\begin{equation}\label{diagrammechangementfiltration}
\begin{tikzcd}
	1 & {\widehat{\Z}'(1)} & {\widehat{\Z}'(1)\rtimes D} & D & 1 \\
	1 & C & {C\times D} & D & 1 \\
	1 & {C\times H} & {C\times D} & {\Gamma_k} & 1 \\
	1 & G & {G\times \Gamma_k} & {\Gamma_k} & 1
	\arrow[from=1-1, to=1-2]
	\arrow[from=1-2, to=1-3]
	\arrow[two heads, from=1-2, to=2-2]
	\arrow[from=1-3, to=1-4]
	\arrow[two heads, from=1-3, to=2-3]
	\arrow["\phi"{description}, curve={height=24pt}, from=1-3, to=4-3]
	\arrow[from=1-4, to=1-5]
	\arrow[equals, from=1-4, to=2-4]
	\arrow[from=2-1, to=2-2]
	\arrow[from=2-2, to=2-3]
	\arrow[hook, from=2-2, to=3-2]
	\arrow[from=2-3, to=2-4]
	\arrow[equals, from=2-3, to=3-3]
	\arrow[from=2-4, to=2-5]
	\arrow[two heads, from=2-4, to=3-4]
	\arrow[from=3-1, to=3-2]
	\arrow[from=3-2, to=3-3]
	\arrow["m"', from=3-2, to=4-2]
	\arrow[from=3-3, to=3-4]
	\arrow[from=3-3, to=4-3]
	\arrow[from=3-4, to=3-5]
	\arrow[equals, from=3-4, to=4-4]
	\arrow[from=4-1, to=4-2]
	\arrow[from=4-2, to=4-3]
	\arrow[from=4-3, to=4-4]
	\arrow[from=4-4, to=4-5]
\end{tikzcd}\end{equation}
where the horizontal maps are the obvious inclusions and projections, the map $C\to C\times H$ is the inclusion of the first factor, $D=H\times \Gamma_k\to \Gamma_k$ is the projection, and $m:C\times H\to G$ is the multiplication map $(\tau,\gamma)\mapsto \tau\gamma$. Thanks to Lemma~\ref{cohomologieequivariante}, the Hochschild-Serre spectral sequence $E_2^{p,q}=\rmH^p(k,\rmH^q_\et(\overline{X},(\Q/\Z)'(2)))\Rightarrow\rmH^{p+q}_\et(X,(\Q/\Z)'(2))$ identifies with $E_2^{p,q}=\rmH^p(k,\rmH^q(G,(\Q/\Z)'(2)))\Rightarrow\rmH^{p+q}(\pi_1^\et(X,\overline{x})^{[p']},(\Q/\Z)'(2))$ in bi-degrees $p+q\le 4$. 

A standard restriction-corestriction argument shows that $\rmH^{3,\norm}(\pi_1^\et(X,\overline{x})^{[p']},(\Q/\Z)'(2))$ is killed by $|G|$. Hence, the long exact cohomology sequence associated to the short exact sequence of trivial $G$-modules $0 \to {\mmu_{|G|}^{\otimes 2}} \to {(\Q/\Z)'(2)} \to {(\Q/\Z)'(2)} \to 0$ shows that the map $\rmH^{3,\norm}(\pi_1^\et(X,\overline{x})^{[p']},{\mmu_{|G|}^{\otimes 2}}) \to\rmH^{3,\norm}(\pi_1^\et(X,\overline{x})^{[p']},(\Q/\Z)'(2))$ induced by the inclusion ${\mmu_{|G|}^{\otimes 2}} \hookrightarrow {(\Q/\Z)'(2)}$ is surjective. Similarly, the map $\rmH^{i}(G,{\mmu_{|G|}^{\otimes 2}}) \to\rmH^{i}(G,(\Q/\Z)'(2))$ is surjective for any $i\ge 1$, so that $\Gamma_k$ acts trivially on both groups.  Since $\pi_1^\et(X,\overline{x})^{[p']}\simeq G\times \Gamma_k$ is a direct product and $\Gamma_k$ acts trivially on $\mmu_{|G|}^{\otimes2}$, it follows from \cite[Theorem 2.4.6]{neukirchschmidtwingberg} that the spectral sequence $E_2^{p,q}=\rmH^p(k,\rmH^q(G,\mmu_{|G|}^{\otimes2}))\Rightarrow\rmH^{p+q}(\pi_1^\et(X,\overline{x})^{[p']},\mmu_{|G|}^{\otimes2})$ degenerates on the second page. 

Denote by $F^\bullet\rmH^3(\pi_1^\et(X,\overline{x})^{[p']},(\Q/\Z)'(2))$ and $F^\bullet\rmH^3(\pi_1^\et(X,\overline{x})^{[p']},\mmu_{|G|}^{\otimes2})$ the filtrations associated to the aforementioned spectral sequences and consider the following commutative diagram:

\begin{equation}\label{diagrammescindageF1}\begin{tikzcd}
	{F^1\rmH^3(\pi_1^\et(X,\overline{x})^{[p']},{\mmu_{|G|}^{\otimes 2}}) } & {\rmH^1(k,\rmH^2(G,{\mmu_{|G|}^{\otimes 2}}))} \\
	{F^1\rmH^3(\pi_1^\et(X,\overline{x})^{[p']},(\Q/\Z)'(2))} & {\rmH^1(k,\rmH^2(G,(\Q/\Z)'(2)))).}
	\arrow[two heads, from=1-1, to=1-2]
	\arrow[from=1-1, to=2-1]
	\arrow[from=1-2, to=2-2]
	\arrow[from=2-1, to=2-2]
\end{tikzcd}\end{equation} 
The short exact sequence $0\to \rmH^1(G,(\Q/\Z)'(2))\to \rmH^2(G,\mmu_{|G|}^{\otimes2})\to \rmH^2(G,(\Q/\Z)'(2))\to 0$ of trivial $\Gamma_k$-modules induces an exact sequence
\[\rmH^1(k,\rmH^2(G,{\mmu_{|G|}^{\otimes 2}}))\to \rmH^1(k,\rmH^2(G,(\Q/\Z)'(2)))\to \rmH^2(k,\rmH^1(G,(\Q/\Z)'(2)))\to \rmH^2(k,\rmH^2(G,{\mmu_{|G|}^{\otimes 2}})).\]
Up to fixing isomorphisms $\rmH^1(G,(\Q/\Z)'(2))\simeq \bigoplus_{i\in I}\mmu_{n_i}$ and $\rmH^2(G,{\mmu_{|G|}^{\otimes 2}})\simeq \bigoplus_{j\in J} \mmu_{n_j}$ for some finite sets $I$ and $J$ and integers $n_i$ for $i\in I$ (resp.\ $n_j$ for $j\in J$), we see that the last map in the above sequence identifies with an inclusion $\bigoplus_{i\in I}\Br(k)[n_i] \hookrightarrow \bigoplus_{j\in J}\Br(k)[n_j]$. This implies that the rightmost vertical arrow in \eqref{diagrammescindageF1} is surjective and hence so is the map $F^1\rmH^3(\pi_1^\et(X,\overline{x})^{[p']},(\Q/\Z)'(2))\to\rmH^1(k,\rmH^2(G,(\Q/\Z)'(2)))$. 

Denote by $F^\bullet_a\rmH^3(C\times D,(\Q/\Z)'(2))$ and $F^\bullet_b\rmH^3(C\times D,(\Q/\Z)'(2))$ the filtrations given by the Hochschild-Serre spectral sequences associated to the third and the second row of \eqref{diagrammechangementfiltration} respectively. The map of filtrations $F^\bullet_a\rmH^3(C\times D,(\Q/\Z)'(2))\to F^\bullet_b\rmH^3(C\times D,(\Q/\Z)'(2))$ viewed as subgroups of $\rmH^3(C\times D,(\Q/\Z)'(2))$ is the identity (as it is induced by the identity on $C\times D$), even though the filtrations are themselves distinct. The group $F^1_b\rmH^3(C\times D,(\Q/\Z)'(2))$ fits into a short exact sequence
\[0\longrightarrow F^2_b\rmH^3(C\times D,(\Q/\Z)'(2))\longrightarrow F^1_b\rmH^3(C\times D,(\Q/\Z)'(2))\longrightarrow \rmH^1(D,\rmH^2(C,(\Q/\Z)'(2))),\]
where the rightmost group is zero since $\rmH^2(C,(\Q/\Z)'(2))$ is ($C$ is cyclic and $(\Q/\Z)'(2)$ is divisible). We thus obtain a map $F^1_a\rmH^3(C\times D,(\Q/\Z)'(2))\to F^2_b\rmH^3(C\times D,(\Q/\Z)'(2))$. Denote by $F^\bullet\rmH^2(D,\widehat{C}(1))$ the filtration given by the Hochschild-Serre spectral sequence $E_2^{p,q}=\rmH^p(k,\rmH^q(H,\widehat{C}(1)))\Rightarrow\rmH^{p+q}(D,\widehat{C}(1))$. We have a commutative diagram 
\begin{equation}\label{diagrammeuniversel1}\begin{tikzcd}
	{F^2_a\rmH^3(C\times D,(\Q/\Z)'(2))} & {F^2\rmH^2(D,\widehat{C}(1))} \\
	{F^1_a\rmH^3(C\times D,(\Q/\Z)'(2))} & {F^1\rmH^2(D,\widehat{C}(1))} \\
	{F^2_b\rmH^3(C\times D,(\Q/\Z)'(2))} & {\rmH^2(D,\widehat{C}(1)),}
	\arrow[from=1-1, to=1-2]
	\arrow[hook, from=1-1, to=2-1]
	\arrow[hook, from=1-2, to=2-2]
	\arrow[from=2-1, to=2-2]
	\arrow[from=2-1, to=3-1]
	\arrow[hook, from=2-2, to=3-2]
	\arrow[from=3-1, to=3-2]
\end{tikzcd}\end{equation}
where the top horizontal arrow is the composite map 
\[F^2_a\rmH^3(C\times D,(\Q/\Z)'(2))\longrightarrow \rmH^2(k,\widehat{C}(1)\times \widehat{H}(1)) \longrightarrow \rmH^2(k,\widehat{C}(1))=F^2\rmH^2(D,\widehat{C}(1)).\]
By analysing the Hochschild-Serre spectral sequence $E_2^{p,q}=\rmH^p(H,\rmH^q(C,(\Q/\Z)'(2)))\Rightarrow \rmH^{p+q}(C\times H,(\Q/\Z)'(2))$, we see that the group $\rmH^2(C\times H,(\Q/\Z)'(2))$ fits into the following split short exact sequence of trivial $\Gamma_k$-modules:
\[0\longrightarrow \rmH^2(H,(\Q/\Z)'(2))\longrightarrow\rmH^2(C\times H,(\Q/\Z)'(2))\longrightarrow \rmH^1(H,\widehat{C}(1))\longrightarrow 0.\]
Taking Galois cohomology, we derive a split short exact sequence:
\begin{equation}\label{definitionflècheH1H2versH2résidu}0\longrightarrow \rmH^1(k,\rmH^2(H,(\Q/\Z)'(2)))\longrightarrow \rmH^1(k,\rmH^2(C\times H,(\Q/\Z)'(2)))\longrightarrow \rmH^1(k,\rmH^1(H,\widehat{C}(1)))\longrightarrow 0.\end{equation}
We then have a commutative diagram with exact rows:

\begin{equation}\label{diagrammeuniversel2}\begin{tikzcd}
	{F^2_a\rmH^3(C\times D,(\Q/\Z)'(2))} & {F^1_a\rmH^3(C\times D,(\Q/\Z)'(2))} & {\rmH^1(k,\rmH^2(C\times H,(\Q/\Z)'(2)))} \\
	{F^2\rmH^2(D,\widehat{C}(1))} & {F^1\rmH^2(D,\widehat{C}(1))} & {\rmH^1(k,\rmH^1(H,\widehat{C}(1))),}
	\arrow[from=1-1, to=1-2]
	\arrow[from=1-1, to=2-1]
	\arrow[from=1-2, to=1-3]
	\arrow[from=1-2, to=2-2]
	\arrow[from=1-3, to=2-3]
	\arrow[from=2-1, to=2-2]
	\arrow[from=2-2, to=2-3]
\end{tikzcd}\end{equation}
where the rightmost arrow is given by \eqref{definitionflècheH1H2versH2résidu}. 

The map $\rmH^2(D,\widehat{C}(1))\to \rmH^2(D,(\Q/\Z)'(1))$ associated to $\widehat{\Z}'(1)\rtimes D \to C\times D$ in \eqref{diagrammechangementfiltration} is induced by the cup-product $\alpha\mapsto \alpha\smile \overline{\mu}$, where we view $\overline{\mu}:\widehat{\Z}'(1)\to G$ as an element of $\rmH^0(D,C(-1))=\rmH^0(k,C(-1))$ (see §\ref{definitionaccouplementcupproduit} for details). Since $\Gamma_k$ acts trivially on $\rmH^1(H,\widehat{C}(1))=\Hom(H,\widehat{C}(1))$, then $\rmH^1(k,\Hom(H,\widehat{C}(1)))=\Hom(\Gamma_k,\Hom(H,\widehat{C}(1)))$ and the cup-product pairing $\rmH^1(k,\rmH^1(H,\widehat{C}(1)))\times \rmH^0(k,C(-1))\to \rmH^1(k,\rmH^1(H,(\Q/\Z)'(1)))$ is non-degenerate on the left. We infer from the combination of \eqref{diagrammeuniversel1} and \eqref{diagrammeuniversel2} that for any $\alpha\in F^1\rmH^3(\pi_1^\et(X,\overline{x})^{[p']},(\Q/\Z)'(2))$ with image $\beta\in \rmH^1(k,\rmH^2(G,(\Q/\Z)'(2)))$, the vanishing of $\partial_{\mu,D}(\alpha)\in \rmH^2(D,(\Q/\Z)'(1))$ implies that of $\beta$ along the composite map 
\[\rmH^1(k,\rmH^2(G,(\Q/\Z)'(2)))\longrightarrow \rmH^1(k,\rmH^2(C\times H,(\Q/\Z)'(2)))\longrightarrow \rmH^1(k,\rmH^1(H,\widehat{C}(1))).\]
Moreover, we see that an element of $\rmH^1(k,\rmH^1(H,\widehat{C}(1)))\simeq\Hom(\Gamma_k,\rmH^1(H,\widehat{C}(1)))$ is trivial if and only if its restriction to $\Hom(\Gamma_k,\rmH^1(Q,\widehat{C}(1)))$ is zero for every cyclic subgroup $Q\subset H$. Let us consider the following commutative diagram:
\[\begin{tikzcd}
	{\Hom(\Gamma_k,\rmH^2(G,(\Q/\Z)'(2)))} & {\Hom(\Gamma_k,\rmH^1(H,\widehat{C}(1)))} \\
	{\Hom(\Gamma_k,\rmH^2(C\times Q,(\Q/\Z)'(2)))} & {\Hom(\Gamma_k,\rmH^1(Q,\widehat{C}(1))),}
	\arrow[from=1-1, to=1-2]
	\arrow[from=1-1, to=2-1]
	\arrow[from=1-2, to=2-2]
	\arrow["\sim", from=2-1, to=2-2]
\end{tikzcd}\]
where the leftmost vertical arrow is induced by the inclusion of a bicyclic subgroup $C\times Q$ in $G$. Varying the choice of $Q$ and $C$, this shows that the vanishing of $\partial_{\mu,D}(\alpha)$ for any pair $(\mu,D)\in {_\mathrm{any}}\mathscr{D}_\mathrm{{\mathrm{any}}}$ implies that $\beta$ dies along the pushforward $\Hom(\Gamma_k,\rmH^2(G,(\Q/\Z)'(2)))\to \Hom(\Gamma_k,\rmH^2(B,(\Q/\Z)'(2)))$ attached to the map $\rmH^2(G,(\Q/\Z)'(2))\to \rmH^2(B,(\Q/\Z)'(2))$ for any bicyclic subgroup $B\subset G$. Therefore, if a class $\alpha\in F^1\rmH^3(\pi_1^\et(X,\overline{x})^{[p']},(\Q/\Z)'(2))$ is globally unramified, then its image $\beta\in \rmH^1(k,\rmH^2(G,(\Q/\Z)'(2)))$ belongs to the subgroup $\Hom(\Gamma_k,B_0(G))$. 

Applying Proposition \ref{suiteexactefondamentale} and Proposition~\ref{imageinverse} to both $X$ and $\overline{X}$ and taking Galois invariants, we obtain a commutative diagram with exact rows:
\[\small\begin{tikzcd}
	0 & {\CH^2(\overline{X})^{\Gamma_k}} & {\rmH^3_\mathrm{gnr}(\pi_1^\et(\overline{X},\overline{x})^{[p']},(\Q/\Z)'(2))^{\Gamma_k}} & {\rmH^3_\nr(\overline{k}(X),(\Q/\Z)'(2))^{\Gamma_k}} \\
	0 & {\CH^2(X)} & {\rmH^3_\mathrm{gnr}(\pi_1^\et(X,\overline{x})^{[p']},(\Q/\Z)'(2))} & {\rmH^3_\nr(k(X),(\Q/\Z)'(2))} & 0.
	\arrow[from=1-1, to=1-2]
	\arrow[from=1-2, to=1-3]
	\arrow[from=1-3, to=1-4]
	\arrow[from=2-1, to=2-2]
	\arrow[from=2-2, to=1-2]
	\arrow[from=2-2, to=2-3]
	\arrow[from=2-3, to=1-3]
	\arrow[from=2-3, to=2-4]
	\arrow[from=2-4, to=1-4]
	\arrow[from=2-4, to=2-5]
\end{tikzcd}\]
The previous observations imply that the group $F^1\rmH^3_\mathrm{gnr}(\pi_1^\et(X,\overline{x}),(\Q/\Z)'(2))$ fits into a short exact sequence
\begin{equation}\label{suitepreuvetrivialitéuniverselle}0\to F^2\rmH^3_\mathrm{gnr}(\pi_1^\et(X,\overline{x})^{[p']},(\Q/\Z)'(2))\to F^1\rmH^3_\mathrm{gnr}(\pi_1^\et(X,\overline{x})^{[p']},(\Q/\Z)'(2)) \to \Hom(\Gamma_k,B_0(G)),\end{equation}
where the group $F^2\rmH^3_\mathrm{gnr}(\pi_1^\et(X,\overline{x})^{[p']},(\Q/\Z)'(2))$ is reduced to $\rmH^3(k,(\Q/\Z)'(2))$ by virtue of Proposition \ref{nonramificationalgébrique} and Proposition \ref{annulation1}. On the other hand, Lemma \ref{surjectivitéCH2} shows that the cokernel of the restriction map $\CH^2(X)\to \CH^2(\overline{X})^{\Gamma_k}$ vanishes. Applying the snake lemma to the previous diagram thus yields the first part of the statement of the theorem. 

Assume now that $\mmu_{|G|^2}\subset k$. As explained in Remark \ref{remarquebogomolov}, we can show that any element of $\Brnr(\overline{X})$ is $\Gamma_k$-invariant and descends to an element of $\Brnr(X)$. On the other hand, up to decomposing $\rmH^0(k,\Brnr(\overline{X}))=\Brnr(\overline{X})$ into a direct sum of elementary abelian subgroups and using the equality $\rmH^1(k,\mmu_{|G|^2})=\Hom(\Gamma_k,\mmu_{|G|^2})$, we see that the cup-product map $\rmH^0(k,\Brnr(\overline{X}))\otimes \rmH^1(k,\mmu_{|G|^2})\to \rmH^1(k,\Brnr(\overline{X})\otimes \mmu_{|G|^2})\simeq \Hom(\Gamma_k,B_0(G))$ is surjective. Remark \ref{remarquebogomolov} further shows that any element in the image of this map comes by construction from a class in $F^1\rmH^3_\mathrm{gnr}(\pi_1^\et(X,\overline{x}),(\Q/\Z)'(2))$. This implies that the last arrow in \eqref{suitepreuvetrivialitéuniverselle} is surjective, hence the second claim in the statement of the theorem.\end{proof}

In the proof, we made use of the following lemma:

\begin{lem}\label{surjectivitéCH2} Assume that $G$ is constant and that $\mmu_{\mathrm{exp}(G)}\subset k$. Then for any field $L$ containing $k$, the base change map $\CH^2(X)\to \CH^2(X\otimes_k L)$ is surjective.
\end{lem}
\begin{proof} Since $\CH^2(X)$ (resp.\ $\CH^2(X\otimes_k L)$) is generated by Chern classes of $k$-representations (resp.\ $L$-representations) of $G$, it suffices to remark that every such $L$-representation descends to a $k$-representation. Given the assumptions, this readily follows from Brauer's theorem on induced characters \cite[§12.3, Théorème 24]{serrereplin} applied to both $k$ and $L$.
\end{proof}

\begin{rmk} In \cite[Chapter 17]{totaro}, Totaro asks whether one can show that the restriction map $\CH^2(B_kG)\to\CH^2(B_{L}G)$ is an isomorphism for any field $L$ containing $k$, as soon as $\mmu_{\exp(G)}\subset k$. Granted this hypothesis, Lemma \ref{surjectivitéCH2} ensures that this map is surjective. We claim that it is an isomorphism under the further assumption that $\mmu_{|G|}\subset k$ and $B_0(G)$ is trivial. 

Let indeed $c\ge 3$ and $(U\subset \A^m_k,c)$ be a classifying pair and $X:=U/G$; Theorem \ref{trivialitépartiealgébrique} implies in this case that in the following commutative diagram:
\[\begin{tikzcd}
	{\CH^2(B_{\overline{k}}G)} & {\rmH^3_\mathrm{gnr}(\overline{X},(\Q/\Z)'(2))} \\
	{\CH^2(B_kG)} & {\rmH^3_\mathrm{gnr}(X,(\Q/\Z)'(2)),}
	\arrow[from=1-1, to=1-2]
	\arrow[from=2-1, to=1-1]
	\arrow[from=2-1, to=2-2]
	\arrow[from=2-2, to=1-2]
\end{tikzcd}\]
the image of $\Ker\big(\CH^2(B_kG)\to \CH^2(B_{\overline{k}}G)\big)$ in $\rmH^3_\mathrm{gnr}(X,(\Q/\Z)'(2))$ lies in $\rmH^3(k,(\Q/\Z)'(2))$, hence is zero since the image of $\CH^2(B_kG)$ meets $\rmH^3(k,(\Q/\Z)'(2))$ trivially. As the bottom horizontal map in the above diagram is injective thanks to Proposition \ref{suiteexactefondamentale}, this implies that the kernel of the restriction map $\CH^2(B_kG)\to \CH^2(B_{\overline{k}}G)$ is trivial. Thus, $\CH^2(B_kG)\to \CH^2(B_{\overline{k}}G)$ is an isomorphism, and similarly when replacing $k$ by $L$. 

Since $\CH^2(B_{\overline{k}}G)\to \CH^2(B_{\overline{L}}G)$ is an isomorphism as well \cite[Corollary 17.3]{totaro}, the desired claim follows. It would be interesting to determine whether the group $\Ker\big(\CH^2(B_kG)\to \CH^2(B_{\overline{k}}G)\big)$ is always trivial, that is, without assuming the vanishing of $B_0(G)$.

\end{rmk}

\section{Fields of cohomological dimension $\le 1$} 

\subsection{A formula for algebraic classes} Let $k$ be a field of cohomological dimension $\le 1$ and exponential characteristic $p>0$. Let $G$ be a finite $k$-group of order coprime to $p$, let $c\ge 3$ and fix a classifying pair $(U\subset\A^m_k,c)$. Denote $X:=U/G$. Since $\mathrm{cd}(k)\le 1$, we see that in the Hochschild-Serre spectral sequence \[E_2^{p,q}=\rmH^p(k,\rmH^q_\et(\overline{X},(\Q/\Z)'(2)))\Rightarrow\rmH^{p+q}(X,(\Q/\Z)'(2)),\] 
we have $F^1\rmH^3_\et(X,(\Q/\Z)'(2))=H^1(k,\rmH^2_\et(\overline{X},(\Q/\Z)'(2)))$ and $F^r\rmH^3_\et(X,(\Q/\Z)'(2))=0$ for $r\ge 2$, as well as an isomorphism $\rmH^2(\overline{X},(\Q/\Z)'(2))\simeq \rmH^2(G(\overline{k}),(\Q/\Z)'(2))$ thanks to Lemma \ref{cohomologieequivariante}.

Fix a pair $(\mu,D)\in{_\mathrm{any}}\mathscr{D}_\mathrm{{\mathrm{pro-}p'}}$, denote $C:=\Im(\mu)$, $H:=D\cap G(\overline{k})$ and let $K$ be the subfield of $\overline{k}$ fixed by the image of $D\to \Gamma_k$. Since $\Gamma_K$ has cohomological dimension $\le 1$, we can fix a continuous section $s:\Gamma_K\to D$ and view $D$ as a semi-direct product $H\rtimes_s\Gamma_K$ endowed with a morphism $H\rtimes_s\Gamma_K\hookrightarrow G(\overline{k})\rtimes\Gamma_k\simeq \pi_1^\et(X,\overline{x})^{[p']}$. If $\sigma\in \Gamma_K$, we let $(\gamma_\sigma,\sigma):=s(\sigma)\in G(\overline{k})\rtimes\Gamma_k$. The $\Gamma_K$-action on $H$ is then given by ${^\sigma}h:=\gamma_\sigma {^\sigma h} \gamma_\sigma^{-1}$, where the action specified on the right is induced by the prescribed action of $\Gamma_K$ on $G(\overline{k})$. By \cite[§2.2, Proposition 8]{stixlivre}, the classes of sections of the tautological short exact sequence
\[1\longrightarrow H\longrightarrow D\longrightarrow \Gamma_K\longrightarrow 1\]
modulo conjugation by an element of $H$ are in bijection with the pointed set $\rmH^1(K,H)$ (for the aforementioned action of $\Gamma_K$ on $H$). Let us form the semi-direct product $C\rtimes D$, where $D$ acts on $C$ through the cyclotomic character $\chi:\Gamma_K\to \mathrm{Aut}(\widehat{\Z}'(1))$. 

The choice of $K$, $\mu$, $H$, action of $\Gamma_K$ on $H$ and class of a section $s:\Gamma_K\to D$ thus determines a split short exact sequence
\begin{equation}\label{produitsemidirectcritèrecorpsdedimension1} 1\longrightarrow C\times H\longrightarrow C\rtimes D\longrightarrow \Gamma_K\longrightarrow 1\end{equation}
and hence an action of $\Gamma_K$ on $C\times H$ given by ${^\sigma}(\tau,h):=(\tau^{\chi(\sigma)},\gamma_\sigma {^\sigma}h\gamma_\sigma^{-1})$.
The Hoschschild-Serre spectral sequence $E_2^{p,q}=\rmH^p(H,\rmH^q(C,(\Q/\Z)'(2)))\Rightarrow \rmH^{p+q}(C\times H,(\Q/\Z)'(2))$ yields a split short exact sequence of $\Gamma_K$-modules:
\[0\longrightarrow \rmH^2(H,(\Q/\Z)'(2))\longrightarrow \rmH^2(C\times H,(\Q/\Z)'(2))\longrightarrow \Hom(H,\widehat{C}(1))\longrightarrow 0.\]
We thus obtain a map $\rmH^1(K,\rmH^2(C\times H,(\Q/\Z)'(2)))\to \rmH^1(K,\Hom(H,\widehat{C}(1)))$.

Since we have the equality ${^\sigma[\tau]} = [\tau^{\chi(\sigma)}]$ in $G(\overline{k})/(\text{conjugacy})$ for every $\tau\in C$ and every $\sigma\in \Gamma_K$, we see that the
image $\overline{\mu}\in \Hom(\widehat{\Z}'(1),G(\overline{k})^\ab)=G^\ab(-1)$ of $\mu$ is $\Gamma_K$-invariant. In particular, we may regard $\overline{\mu}$ as an element of $\rmH^0(D,C(-1))=\rmH^0(K,C(-1))$ and consider the cup-product map $\rmH^1(K,\Hom(H,\widehat{C}(1)))\to \rmH^1(K,\widehat{H}^\ab)$ given by $\beta\mapsto\beta\smile \overline{\mu}$. 

On the other hand, the multiplication $C\times H\to G(\overline{k}),~(\tau,\gamma)\mapsto \tau\gamma$ is $\Gamma_K$-equivariant when twisting the natural $\Gamma_k$-action on $G(\overline{k})$ by conjugation by $\gamma_\sigma\in H$. Since the morphism $G(\overline{k})\to G(\overline{k})$, $g\mapsto \gamma_{\sigma}g\gamma_\sigma^{-1}$ induces the identity on $\rmH^1(k,\rmH^2(G(\overline{k}),(\Q/\Z)'(2)))$, we have a well-defined map $\rmH^1(k,\rmH^2(G(\overline{k}),(\Q/\Z)'(2)))\to \rmH^1(K,\rmH^2(C\times H,(\Q/\Z)'(2)))$. We hence have a composite map
\[\lambda_{\mu,D}:\rmH^1(k, \rmH^2(G(\overline{k}),(\Q/\Z)'(2)))\longrightarrow \rmH^1(K,\rmH^2(C\times H,(\Q/\Z)'(2)))\longrightarrow \rmH^1(K,\Hom(H,\widehat{C}(1))).\]
\begin{prop}\label{nonramificationalgébriquecorpsdedimension1}
Consider the subgroup \[\Sha^1_\mathrm{nr}(k, \rmH^2(G(\overline{k}),(\Q/\Z)'(2))\subset\rmH^1(k, \rmH^2(G(\overline{k}),(\Q/\Z)'(2)))\] 
consisting of those classes $\alpha$ such that for every pair $(\mu,D)\in{_\mathrm{any}}\mathscr{D}_\mathrm{{\mathrm{pro-}p'}}$, if we denote $H:= D\cap G(\overline{k})$ and $K\subset \overline{k}$ the field fixed by the image of $D\to \Gamma_k$, we have $\lambda_{\mu,D}(\alpha)\smile\overline{\mu}=0\in\rmH^1(K,\widehat{H}^\ab)$. Consider the composite map
\begin{equation}\label{flècheH1H2versH0H3}\rmH^1(k, \rmH^2(G(\overline{k}),(\Q/\Z)'(2)))\longrightarrow \rmH^3_\et(X,(\Q/\Z)'(2))\longrightarrow \rmH^0(X,\mathscr{H}^3_\et((\Q/\Z)'(2))),\end{equation}
where the first arrow is the inclusion. Then a class $\alpha\in \rmH^1(k,\rmH^2(G(\overline{k}),(\Q/\Z)'(2)))$ lies in $\Sha^1_\mathrm{nr}(k, \rmH^2(G(\overline{k}),(\Q/\Z)'(2)))$ if and only if its image along the above map lies in $\rmH^3_\mathrm{nr}(k(X),(\Q/\Z)'(2)).$ More precisely, we have an equality:
   \[F^1\rmH^3_\mathrm{gnr}(X,(\Q/\Z)'(2))=\Sha^1_\mathrm{nr}(k, \rmH^2(G(\overline{k}),(\Q/\Z)'(2))).\]
\end{prop}
\begin{proof} Let $\alpha\in \rmH^1(k, \rmH^2(G(\overline{k}),(\Q/\Z)'(2)))$ and view it as an element of $\rmH^3_\et(X,(\Q/\Z)'(2))$. Proposition \ref{imageinverse} implies that the image of $\alpha$ in $\rmH^0(X,\mathscr{H}^3_\et((\Q/\Z)'(2)))$ along \eqref{flècheH1H2versH0H3} is unramified if and only if for every pair $(\mu,D)\in {_\mathrm{any}}\mathscr{D}_\mathrm{{\mathrm{pro-}p'}}$, the residue $\partial_{\mu,D}(\alpha)\in \rmH^2(D,(\Q/\Z)'(1))$ is zero. Let us fix such a pair $(\mu,D)$,  let $K$ be the subfield of $\overline{k}$ fixed by the image of $D\to\Gamma_k$, denote $C:=\Im(\mu)$, $H:=D\cap G(\overline{k})$ and fix a continuous section $s:\Gamma_K\to D$. By virtue of the functoriality of the Hochschild-Serre spectral sequence for the semi-direct product decompositions of $\widehat{\Z}'(1)\rtimes D$, $C\rtimes D$ and $\pi_1^\et(X,\overline{x})^{[p']}\simeq G(\overline{k})\rtimes\Gamma_k$ respectively, we derive a commutative diagram:
\[\begin{tikzcd}
	{F^1\rmH^3(\pi_1^\et(X,\overline{x})^{[p']},(\Q/\Z)'(2))} & {\rmH^2(D,\widehat{C}(1))} & {\rmH^2(D,(\Q/\Z)'(1))} \\
	{\rmH^1(k, \rmH^2(G(\overline{k}),(\Q/\Z)'(2)))} & {\rmH^1(K,\Hom(H,\widehat{C}(1)))} & {\rmH^1(K,\widehat{H}^\ab),}
	\arrow[from=1-1, to=1-2]
	\arrow["{\partial_{\mu,D}}", curve={height=-18pt}, from=1-1, to=1-3]
	\arrow[from=1-2, to=1-3]
	\arrow[hook, from=2-1, to=1-1]
	\arrow["{\lambda_{\mu,D}}"', from=2-1, to=2-2]
	\arrow[hook, from=2-2, to=1-2]
	\arrow[from=2-2, to=2-3]
	\arrow[hook, from=2-3, to=1-3]
\end{tikzcd}\]
where the second (resp.\ third) injective vertical arrow is constructed by considering the Hochschild-Serre spectral sequence $E_2^{p,q}=\rmH^p(K,\rmH^q(H,\widehat{C}(1)))\Rightarrow\rmH^{p+q}(D,\widehat{C}(1))$ (resp.\ the spectral sequence $E_2^{p,q}=\rmH^p(K,\rmH^q(H,(\Q/\Z)'(1)))\Rightarrow\rmH^{p+q}(D,(\Q/\Z)'(1))$) and the equality $\rmH^1(K,\Hom(H,\widehat{C}(1))=F^1\rmH^2(D,\widehat{C}(1))$  (resp.\ the equality $\rmH^1(K,\widehat{H}^\ab)=F^1\rmH^2(D,(\Q/\Z)'(1))$) coming from the fact that $K$ has cohomological dimension $\le 1$. 

On the other hand, the map $\rmH^1(K,\Hom(H,\widehat{C}(1)))\to \rmH^1(K,\widehat{H}^\ab)$ is given by $\beta\mapsto\beta\smile\overline{\mu}$. Since the rightmost vertical arrow in the above diagram is injective, we infer that $\partial_{\mu,D}(\alpha)=0\in \rmH^2(D,(\Q/\Z)'(1))$ if and only if $\lambda_{\mu,D}(\alpha)\smile\overline{\mu}=0\in\rmH^1(K,\widehat{H}^\ab)$, as desired.\end{proof}

\begin{thm}\label{cohomologiecorpsdedimension1}  There is a short exact sequence:
\begin{equation}\begin{tikzcd}
 0 \rar & \Ker\Big(\CH^2(X)\to \CH^2(\overline{X})\Big) \rar & \Sha^1_\nr(k, \rmH^2(G(\overline{k}),(\Q/\Z)'(2))) \ar[rounded corners,
   to path={ -- ([xshift=2ex]\tikztostart.east)
             |- ($(\tikztostart)!.5!(\tikztotarget)$) \tikztonodes
             -| ([xshift=-2ex]\tikztotarget.west)
             -- (\tikztotarget)}]{dl} \\
 & \rmH^3_{\nr,\mathrm{alg}}(k(X),(\Q/\Z)'(2)) \rar & \Coker\Big(\CH^2(X)\to \CH^2(\overline{X})^{\Gamma_k}\Big).
\end{tikzcd}\end{equation}
\end{thm}
\begin{rmk} Let us discuss a sufficient condition under which the exact sequence from Theorem \ref{cohomologiecorpsdedimension1} simplifies.

(i) If $G$ is constant and the Chern class $c_2:\mathrm{Fil}^2K_0(\overline{X})^{\Gamma_k}\to \CH^2(\overline{X})^{\Gamma_k}$ is surjective, then the group  $\Coker\big(\CH^2(X)\to \CH^2(\overline{X})^{\Gamma_k}\big)$ is trivial. Indeed, since $\mathrm{gr}^2K_0(X)\simeq \mathrm{gr}^2\mathrm{Rep}_k(G)$ and $\mathrm{gr}^2K_0(\overline{X})\simeq \mathrm{gr}^2\mathrm{Rep}_{\overline{k}}(G)$, then the claim follows from the fact that the natural restriction map $\mathrm{Rep}_k(G)\to \mathrm{Rep}_{\overline{k}}(G)^{\Gamma_k}$ is surjective. This holds because the Brauer group of $k$ is trivial, see for instance \cite[§12.2, Corollaire]{serrereplin}.

(ii) The condition that the Chern class $c_2:\mathrm{Fil}^2K_0(\overline{X})^{\Gamma_k}\to \CH^2(\overline{X})^{\Gamma_k}$ be surjective is for instance satisfied when all the irreducible characters of $G$ over $\overline{k}$ are $\Gamma_k$-invariant (this holds for instance for symmetric groups). One can however find counter-examples to this property when considering e.g.\ generalised quaternion groups: see the proof of Lemma~\ref{CH2Q_16} applied to $k=\F_p$.

\end{rmk}
\begin{proof}[Proof of Theorem \ref{cohomologiecorpsdedimension1}] Applying Proposition \ref{imageinverse} to both $X$ and $\overline{X}$ and taking Galois invariants, we obtain a commutative diagram with exact rows:
\[\small\begin{tikzcd}
	0 & {\CH^2(\overline{X})^{\Gamma_k}} & {\rmH^3_\mathrm{gnr}(\pi_1^\et(\overline{X},\overline{x})^{[p']},(\Q/\Z)'(2))^{\Gamma_k}} & {\rmH^3_\nr(\overline{k}(X),(\Q/\Z)'(2))^{\Gamma_k}} \\
	0 & {\CH^2(X)} & {\rmH^3_\mathrm{gnr}(\pi_1^\et(X,\overline{x})^{[p']},(\Q/\Z)'(2))} & {\rmH^3_\nr(k(X),(\Q/\Z)'(2))} & 0.
	\arrow[from=1-1, to=1-2]
	\arrow[from=1-2, to=1-3]
	\arrow[from=1-3, to=1-4]
	\arrow[from=2-1, to=2-2]
	\arrow[from=2-2, to=1-2]
	\arrow[from=2-2, to=2-3]
	\arrow[from=2-3, to=1-3]
	\arrow[from=2-3, to=2-4]
	\arrow[from=2-4, to=1-4]
	\arrow[from=2-4, to=2-5]
\end{tikzcd}\] 
by virtue of Proposition \ref{nonramificationalgébriquecorpsdedimension1}, we have $F^1\rmH^3_\mathrm{gnr}(X,(\Q/\Z)'(2))=\Sha^1_\nr(k,\ \rmH^2(G(\overline{k}),(\Q/\Z)'(2)))$. A straightforward application of the snake lemma thus provides the short exact sequence in the statement of the theorem.\end{proof}

\subsection{An algorithm over a finite field} Let now $k$ be a finite field of characteristic $p$. The computation of unramified cohomology in degree $3$ for various classes of smooth and projective varieties over $k$ is a delicate and open task which is of prime interest, due to its applications to stable-rationality and its connections to questions relative to algebraic cycles in codimension $2$; a wide variety of sources in the literature discuss such questions, see e.g.\ \cite{colliotkahn}, \cite{colliotscavia}, \cite{colliotvoisin} and \cite{pirutka}. 

However, the existence of an algorithm that effectively computes the unramified cohomology in degree $3$ for a large class of varieties does not seem to be achievable in general. Let $G$ be a finite $k$-group of order coprime to $p$, let $c\ge 3$, let $(U\subset \A^m_k,c)$ be a classifying pair, and denote $X:=U/G$. In this section, we present an effective method to partially compute the unramified cohomology group $\rmH^3_\nr(k(X),(\Q/\Z)'(2))$. When $G$ is constant of odd order, we explain how this method can be refined to compute the group $\rmH^3_\nr(k(X),(\Q/\Z)'(2))$.
\begin{thm}\label{algorithmecorpsfinis} There exists an algorithm that, given the datum of a finite field $k$ of characteristic $p$ and a finite $k$-group $G$ of odd order coprime to $p$, computes the unramified cohomology group $\rmH^3_\nr(k(U/G),(\Q/\Z)'(2))$ for any classifying pair $(U\subset \A^m_k,c)$ with $c\ge 3$.
\end{thm}

\subsubsection{Proof of Theorem \ref{algorithmecorpsfinis}} The computation of the unramified cohomology group can be conducted as follows: we will prove the existence of the following sequence of algorithms, which combine to provide the announced algorithm in the statement of the theorem.

\begin{lem}\label{algorithmeH3} There exists an algorithm that, given the datum of a finite field $k$ of characteristic $p$ and a finite $k$-group $G$ of order coprime to $p$, returns a finite group that is isomorphic to the cohomology group $\rmH^3_\et(U/G,(\Q/\Z)'(2))$ together with an explicit identification of its cohomology classes for any classifying pair $(U\subset \A^m_k,c)$ with $c\ge 3$.
\end{lem}
\begin{proof}  It suffices to note that $\rmH^3_\et(U/G,(\Q/\Z)'(2))$ can be computed thanks to Lemma~\ref{cohomologieequivariante} as the cohomology of the profinite group $\pi_1^\et(U/G,\overline{x})^{[p']}\simeq G(\overline{k})\rtimes\Gamma_k$. Let us choose a topological generator $\sigma$ of $\Gamma_k$. We claim that one can specify an integer $N\ge 1$ such that the pullback $\rmH^3(G(\overline{k})\rtimes (\Gamma_k/\langle \sigma^N\rangle),(\Q/\Z)'(2))\to \rmH^3(G(\overline{k})\rtimes \Gamma_k,(\Q/\Z)'(2))$ is an isomorphism. Let $d$ be the degree of the minimal Galois extension that splits $G$ and contains the roots of unity of order dividing $|G|$. The Hochschild-Serre spectral sequence $E_2^{p,q}=\rmH^p(G(\overline{k})\rtimes (\Gamma_k/\langle\sigma^d\rangle),\rmH^q(\langle\sigma^d\rangle,\mmu_{d}^{\otimes 2}))\Rightarrow \rmH^{p+q}(G(\overline{k})\rtimes \Gamma_k,\mmu_d^{\otimes2})$ provides an injection $\rmH^3(G(\overline{k})\rtimes (\Gamma_k/\langle \sigma^d\rangle),(\Q/\Z)'(2))\hookrightarrow \rmH^3(G(\overline{k})\rtimes \Gamma_k,\mmu_d^{\otimes2})$ whose cokernel is contained in $\rmH^2\big(G(\overline{k})\rtimes (\Gamma_k/\langle \sigma^d\rangle),\rmH^1(\langle \sigma^N\rangle,\mmu_d^{\otimes2}))\big)$. 

On the other hand, for any $n\ge1$, the map $\rmH^1(\langle \sigma^d\rangle,\mmu_d^{\otimes2})\to \rmH^1(\langle \sigma^{nd}\rangle,\mmu_d^{\otimes2})$ is the multiplication by $n$. If we thus let $n=d$, we see that the latter map is zero. Denote $N:=d^2$. By pushing forward along the inclusion $G(\overline{k})\rtimes(\Gamma/\langle\sigma^d\rangle)\hookrightarrow G(\overline{k})\rtimes(\Gamma/\langle\sigma^{N}\rangle)$, we infer that the map $\rmH^3(G(\overline{k})\rtimes (\Gamma_k/\langle \sigma^{N}\rangle),\mmu_d^{\otimes2})\hookrightarrow \rmH^3(G(\overline{k})\rtimes \Gamma_k,\mmu_d^{\otimes2})$ is surjective and hence an isomorphism. Since $\rmH^3(G(\overline{k})\rtimes \Gamma_k,(\Q/\Z)'(2))$ is killed by $|G|$ and hence by $d$, the Kummer sequence $0 \to {\mmu_{d}^{\otimes 2}} \to {(\Q/\Z)'(2)} \to {(\Q/\Z)'(2)} \to 0$ provides a commutative diagram
\[\begin{tikzcd}
	{\rmH^3(G(\overline{k})\rtimes (\Gamma_k/\langle \sigma^{N}\rangle),\mmu_d^{\otimes2})} & {\rmH^3(G(\overline{k})\rtimes \Gamma_k,\mmu_d^{\otimes2})} \\
	{\rmH^3(G(\overline{k})\rtimes (\Gamma_k/\langle \sigma^{N}\rangle),(\Q/\Z)'(2))} & {\rmH^3(G(\overline{k})\rtimes \Gamma_k,(\Q/\Z)'(2)),}
	\arrow["\sim", from=1-1, to=1-2]
	\arrow[two heads, from=1-1, to=2-1]
	\arrow[two heads, from=1-2, to=2-2]
	\arrow[hookrightarrow, from=2-1, to=2-2]
\end{tikzcd}\]
with surjective vertical arrows, so that the bottom horizontal map is an isomorphism as well. It therefore suffices to compute the third cohomology group for the finite group $G(\overline{k})\rtimes (\Gamma_{k}/\langle \sigma^N\rangle)$ acting on $(\Q/\Z)'(2)$ through the cyclotomic action of $\Gamma_{k}/\langle \sigma^N\rangle$; this may be carried out algorithmically, for instance by means of the HAP package for the computer algebra system GAP.
\end{proof}
\begin{lem}\label{algorithmeresidus} Let $k$ be a finite field of characteristic $p$, let $G$ be a finite $k$-group of order coprime to $p$, let $c\ge 3$ and fix a classifying pair $(U\subset \A^m_k,c)$. For any class $\alpha\in\rmH^3_\et(U/G,(\Q/\Z)'(2))$, there exists an algorithm that computes the residues $\partial_{\mu,D}(\alpha)$ for any pair $(\mu,D)\in{_\mathrm{any}}\mathscr{D}_\mathrm{{\mathrm{pro-}p'}}$.
\end{lem}
\begin{proof} First assume that $G$ is constant and $\mmu_{\exp(G)} \subset k$. We then see that the group $\pi_1^\et(X,\overline{x})^{[p']}$ is isomorphic to the direct product $G\times \Gamma_k$ and the cyclotomic character $\chi:\Gamma_k\to (\Z/\mathrm{exp}(G))^*$ is trivial. In particular, in the criterion from Proposition~\ref{imageinverse}, we may only consider pairs $(\mu,D)$ where $\mu:\widehat{\Z}'(1)\to G$ is a morphism, $D\subset \pi_1^\et(X,\overline{x})^{[p']}$ is a closed subgroup, and the image of $\mu$ commutes with any element of $D\cap G$. Let $H\subset G$ be the subgroup such that $\langle D,\Gamma_k\rangle=H\times \Gamma_k$ and denote $C:=\Im(\mu)$. Given the inclusion $\langle C,D\rangle\subset \langle C,H\rangle\times\Gamma_k$, the residue $\partial_{\mu,H\times \Gamma_k}$ factors $\partial_{\mu,D}$. We hence need only consider those pairs $(\mu,D)$ such that $D=H\times \Gamma_k$ is a direct product. These subgroups are finitely many, so that a class $\alpha\in \rmH^3(\pi_1^\et(X,\overline{x}),(\Q/\Z)'(2))$ lies in $\rmH^3_{\mathrm{gnr}}(\pi_1^\et(X,\overline{x}),(\Q/\Z)'(2)$ if and only if it vanishes along the finitely many residues defined by these subgroups.

Let now $G$ and $k$ be arbitrary, and consider the minimal finite Galois sub-extension $k\subset L$ of $\overline{k}$ which splits $G$ and contains $\mmu_{\exp(G)}$. Let $(\mu,D)\in{_\mathrm{any}}\mathscr{D}_\mathrm{{\mathrm{pro-}p'}}$; since $k$ is finite, the subfield $K$ of $\overline{k}$ associated to this residue is either contained in $L$ or contains this field.
\begin{enumerate}
\item[(i)] In the first case, there are finitely many intermediate field extensions $k\subset K\subset L$. For a fixed $K$, we may enumerate all the finitely many pairs of morphisms $\mu:\widehat{\Z}'(1)\to G(\overline{k})$ and subgroups $H\subset G(\overline{k})$ which satisfy the semi-commutation relation $\gamma{^\sigma}\tau\gamma^{-1}=\tau^{\chi(\sigma)}$ for any $\tau\in\Im(\mu)$, any $\gamma\in H$ and any $\sigma\in \Gal(K/k)$. 
\item[(ii)] In the second case, the beginning of the proof ensures that we may only consider finitely many residues. We may thus enumerate all the finitely many pairs of morphisms $\mu:\widehat{\Z}'(1)\to G(\overline{k})$ and subgroups $H\subset G(\overline{k})$ which satisfy the commutation relation $\gamma\tau\gamma^{-1}=\tau$ for any $\tau\in\Im(\mu)$ and any $\gamma\in H$.
\end{enumerate}
The datum of each pair $(\mu,D)$ in case (i) or (ii) allows us to construct a residue map $\partial_{\mu,D}:\rmH^3(\pi_1^\et(X,\overline{x})^{[p']},(\Q/\Z)'(2))\to \rmH^2(D,(\Q/\Z)'(1))$, and the computation of $\partial_{\mu,D}(\alpha)$ can be carried out using e.g.\ the HAP package for the computer algebra system GAP. We may thus form a list $R$ of the values $\partial_{\mu,D}(\alpha)$ in case (i) as well as a list $R'$ of the values $\partial_{\mu,D}(\alpha)$ in case (ii). The concatenation of $R$ and $R'$ accounts for all the possible residues values of $\alpha$.
\end{proof}
\begin{lem}\label{algorithmenegligeables} There exists an algorithm that, given the datum of a finite field $k$ of characteristic $p$ and a finite $k$-group $G$ of order coprime to $p$, computes the group of permutation negligible classes $\rmH^3_\mathrm{per}(U/G,(\Q/\Z)'(2))$ for any classifying pair $(U\subset \A^m_k,c)$ with $c\ge 3$.
\end{lem}
\begin{proof} We first find a list $Q:=\{q_1,\hdots,q_n\}$ of the finitely many prime numbers that divide the order of $G$. For any such $q\in Q$, we can then exhibit an associated $q$-Sylow subgroup $G_q$ of $G$, which we fix once and for all. If we denote $K:=k(\mmu_q)$, we then find a list $S_q:=\{H_1,\hdots, H_r\}$ of all the subgroups of $G_q$. For each such $H\in S_q$, we provide a list $M_{H,q}:=\{\phi_1,\hdots,\phi_s\}$ of all the $\Gamma_K$-equivariant morphisms $H\to(\Q/\Z)'(1)$. The computation of the cup-product map
\[\Pic(U_K/H)^{\otimes 2}\longrightarrow\rmH^1_\et(U_K/H,(\Q/\Z)'(1))^{\otimes 2}\longrightarrow \rmH^3_\et(U_K/H,(\Q/\Z)'(2))\]
can be conducted as follows: for any pair $(\phi,\phi')\in M_{H,q}$, we may view these classes as elements of $\rmH^1(H\times \Gamma_K,(\Q/\Z)'(1))$ and compute the cup-product $\phi\smile \phi'\in \rmH^3(H\times \Gamma_K,(\Q/\Z)'(2))$, as detailed e.g.\ in \cite[Chapter V, §3]{browncohomology}. This returns a list $C_{H,q}:=\{c_1,\hdots,c_m\}$ of elements of $\rmH^3(H\times \Gamma_K,(\Q/\Z)'(2))$. 

It then suffices to compute the image of these elements along the corestriction map $\rmH^3(H\times \Gamma_K,(\Q/\Z)'(2))\to \rmH^3(G\times \Gamma_k,(\Q/\Z)'(2))$ associated to the open subgroup $H\times \Gamma_K\subset G\times\Gamma_k$; a detailed treatment is provided in \cite[Chapter I, §5]{neukirchschmidtwingberg}. This returns a list $C'_{H,q}:=\{c'_1,\hdots, c'_m\}$ of elements of $\rmH^3(G\times \Gamma_k,(\Q/\Z)'(2))$. Varying the choice of $q$ and $H$ and concatenating the lists $C'_{H,q}$, we obtain a list $B:=\{\beta_1,\hdots, \beta_r\}$ of elements of $\rmH^3(G\times \Gamma_k,(\Q/\Z)'(2))$; for any given class $\beta_i\in B$, we may then discard all the $\beta_j$'s with $j>i$ that are equal to $\beta_i$. We hence obtain a list which describes the group of permutation negligible classes.\end{proof}

\begin{proof}[Proof of Theorem \ref{algorithmecorpsfinis}] By Lemma \ref{algorithmeH3}, there exists an algorithm that returns a finite list $A:=\{\alpha_1,\hdots,\alpha_j\}$ as well as operations between these elements that determine the structure of the group $\rmH^3_\et(U/G,(\Q/\Z)'(2))$. For a given $\alpha\in A$, we may apply the algorithm of Lemma \ref{algorithmeresidus} in order to compute the residues $\partial_{\mu,D}(\alpha)$. If one of these residues is non-zero, we may already stop this algorithm, discard $\alpha$ and replace it with the next element of $A$. Otherwise, $\alpha$ is globally unramified. By doing so for every element of $A$, we construct a list $A'\subset A$ of elements that correspond to all the globally unramified classes of $\rmH^3_\et(U/G,(\Q/\Z)'(2))$, together with its inherited group structure. 

Lemma \ref{algorithmenegligeables} provides a list $B:=\{\beta_1,\hdots,\beta_r\}$ which (up to rearranging) is contained in $A'$ and coincides as a set with that of permutation negligible classes in $\rmH^3_\et(U/G,(\Q/\Z)'(2))$. We may then specify the group structure on the elements of $B$ as inherited from that of $A'$. Quotienting the group arising from $A'$ by the one associated to $B$ returns a finite group isomorphic to $\rmH^3_\nr(k(U/G),(\Q/\Z)'(2))$.
\end{proof}

\section{Counterexamples to Noether's problem over $\Q$} \label{paragraphenoether}

\subsection{An example of Serre}\label{paragrapheexempleserre} Let $Q_{16}$ be the generalised quaternion group of order $16$. A well-known result of Serre \cite[Theorem 34.7]{garibaldimerkurjevserre} asserts that Noether's problem has a negative answer for the field of invariants $\Q(V)^{Q_{16}}$ associated to any faithful $\Q$\nobreakdash-representation $Q_{16}\to \GL(V)$. The proof of Serre invokes Rost's invariant for the simple simply-connected $\Q$-group $\mathrm{Spin}_{7,\Q}$ combined with an adequate embedding $Q_{16}\hookrightarrow \mathrm{Spin}_{7,\Q}$. Notably, he shows that there exists a degree $3$ normalised cohomological invariant of $\mathrm{Spin}_{7,\Q}$ whose restriction to $Q_{16}$ is unramified and non-zero. The obtained element corresponds to an unramified cohomology class in $\rmH^3_\nr(\Q(V)^{Q_{16}},\Q/\Z(2))$ which is not a constant; see e.g.\ \cite[Appendix C]{garibaldimerkurjevserre} for a proof of this claim. 

Let $c\ge 3$, fix a classifying pair $(U\subset\A^m_\Q,c)$ for $Q_{16}$, and denote $X:=U/Q_{16}$. Equivalently, Serre's proof provides a non-trivial element of $\rmH^{3,\norm}_\nr(\Q(X),\Q/\Z(2))$. In this section, we aim at describing this class explicitly and independently from Serre's approach thanks to Theorem \ref{theoremeprincipal}. We notably show that it is the only non-zero class, namely, the group $\rmH^{3,\norm}_\nr(\Q(X),\Q/\Z(2))$ has order $2$ and it is generated by this element. 

\subsubsection{Noether's problem for generalised quaternion groups}\label{paragrapheexempleserre2} More generally, let $Q_{2^n}$ be the generalised quaternion group of order $2^n$ for $n\ge 3$, given by the presentation
\[Q_{2^n}:=\Big\langle r,s~\Big|~r^{2^{n-1}}=1,~s^2=r^{2^{n-2}},~srs^{-1}=r^{-1}\Big\rangle.\]
As detailed in \cite[§33.28]{garibaldimerkurjevserre}, one can see that Serre's argument for $Q_{16}$ cannot be extended to $Q_{2^n}$ for any $n\ge 5$: indeed, one cannot associate a cohomological invariant to $Q_{2^n}$ through the same construction as in \emph{loc.\ cit}.\ Our approach allows us to generalise our computation for $Q_{16}$ to all generalised quaternion groups $Q_{2^n}$ for $n\ge 3$, thereby settling Noether's problem over $\Q$ in the negative for all such groups. Moreover, we observe that the third unramified cohomology for $Q_{2^n}$ is reduced to the constants as soon as we extend scalars to e.g.\ $\Q(\zeta_{2^{n-2}})$, hence shedding light on the known results in the literature; see Remark~\ref{remarquequaternionsobstructionnoether} for more detail. Let $c\ge 3$, fix a classifying pair $(U\subset\A^m_\Q,c)$ for $Q_{2^n}$, and denote $X:=U/Q_{2^n}$.

\begin{thm}\label{Q_16}
    Let $E:=K(\zeta_{2^{n-1}}+\zeta_{2^{n-1}}^{-1})$. There is an isomorphism
    \[\rmH^{3,\norm}_\nr(K(X),\Q/\Z(2)) \simeq\begin{cases} \Z/2 & \text{ if }K\neq E\text{ and }(-1,-1)_E\neq 0\in \Br(E),\\ 0 & \text{ otherwise}.\end{cases}\]
    More precisely, the map $\rmH^{3,\norm}_\nr(K(X),\Q/\Z(2))\to \Coker\big(\CH^2(X\otimes_\Q K)\to \CH^2(X\otimes_\Q \overline{K})^{\Gamma_K}\big)$ is always injective, and it is an isomorphism in the first case above. Consequently, in the first case, the group $\rmH^{3,\norm}_\nr(K(X),\Q/\Z(2))$ is generated by the image of $[E:K]\cdot c_2(\rho)$ where $\rho:Q_{2^n}\to \SL_2(\overline{K})$ denotes any irreducible and faithful $\overline{K}$-representation. 
\end{thm}

\subsection{Proof of Theorem \ref{Q_16}} The proof of Theorem~\ref{Q_16} is carried out through the following lemmas, whose combination yields the desired computation.

\subsubsection{An upper bound for the unramified cohomology} We first show that for any field $K\supset \Q$, the group $\rmH^{3,\norm}_\nr(K(X),\Q/\Z(2))$ is always at most a group of order $2$.

\begin{lem}\label{Sha2Q_16} For any field extension $K\supset \Q$, we have 

\[\Sha^2_\nr(K,\widehat{Q}^\ab_{2^n}(1))\simeq\begin{cases} (\Z/2)^2 & \text{ if } (-1,-1)_K\neq 0,\\ 0 & \text{ else.}\end{cases}\] More precisely, in the first case this group is generated by $(0,(-1,-1)_K)$ and $((-1,-1)_K,0)$ viewed as elements of $ \rmH^2(K,\widehat{Q}^\ab_{2^n}(1))\simeq\Br(K)[2]\times \Br(K)[2]$.
\end{lem}
\begin{proof} Since the exponent of $Q_{2^n}$ is equal to $2^{n-1}$, the action of the cyclotomic character of $\Gamma_K$ on $Q_{2^n}/(\mathrm{conjugacy})$ factors through that of $\Gal(K(\zeta_{2^{n-1}})/K)$. On the other hand, because $Q_{2^n}^\ab=\langle\overline{r}\rangle\times\langle\overline{s}\rangle\simeq \Z/2\oplus \Z/2$, we see that $Q_{2^n}^\ab(i)\simeq\mmu_2\oplus\mmu_2$ as a $\Gamma_K$-module for every $i\in\Z$, so that $ \rmH^2(K,\widehat{Q}^\ab_{2^n}(1))\simeq\Br(K)[2]\times \Br(K)[2]$. Let us therefore write $\alpha=(a,b)$ for a given element of $\rmH^2(K,\widehat{Q}^\ab_{2^n}(1))$. If $(\mu,D)\in{_\mathrm{any}}\mathscr{D}_\mathrm{{\mathrm{pro-}p'}}$, we have that $\overline{\mu}\in \rmH^0(D,Q_{2^n}^\ab(-1))\simeq\rmH^0(D,\mmu_2\oplus \mmu_2)$ is $\Gamma_K$-invariant. Therefore, if $\rho:D\to \Gamma_K$ is the associated morphism, we have that $\rho^*(\alpha)\smile \overline{\mu}=\rho^*(\alpha\smile \overline{\mu})$ where $\overline{\mu}$ is alternatively seen as an element of $\rmH^0(K,Q_{2^n}^\ab(-1))$. If $\overline{\mu}=0$, then there is nothing to show. We claim that $\rho^*((-1,-1)_K)=0$ for any $D$ such that $\overline{\mu}\neq 0$. The semi\nobreakdash-commutation relation yields a morphism $\psi:D\to N_{Q_{2^n}}(\Im(\mu))$; here, if $\gamma\in D$, we have that $\psi(\gamma)$ acts on a generator $\tau\in\Im(\mu)$ as $\chi(\rho(\gamma))\text{ mod }m$ where $m:=|\Im(\mu)|$. Otherwise, we have either $N_{Q_{2^n}}(\Im(\mu))=Q_8$ or $N_{Q_{2^n}}(\Im(\mu))=Q_{2^n}$. In both cases, $N_{Q_{2^n}}(\Im(\mu))$ is the unique non-split extension
\[1\longrightarrow \mmu_m^{-}\longrightarrow N_{Q_{2^n}}(\Im(\mu))\overset{\pi}{\longrightarrow} \Z/2\longrightarrow 1\]
where $\mmu_m^{-}$ denotes the $\Z/2$-module $\Z/m$ equipped with the action by inversion $\zeta\mapsto \zeta^{-1}$. Let $\epsilon\in \rmH^2(\Z/2,\mmu_m^{-})$ denote the class of the above extension. If $L\subset \overline{K}$ denotes the subfield fixed by the image of $D\to \Gamma_K$ and $c$ is the character of the extension $L\subset L(i)$, then $\pi\circ\psi=c\circ \rho$. By virtue of the periodicity of Tate cohomology \cite[Proposition~1.7.1]{neukirchschmidtwingberg} applied to the inclusion $\mmu_m\subset L(i)^*$, and since the class corresponding to $-1$ in $\Br(L(i)/L)$ is exactly $(-1,-1)_L$, we infer that the image of $c^*(\epsilon)$ in $\Br(L)$ is precisely $(-1,-1)_L$. Hence, $\rho^*((-1,-1)_K)=(c\circ\rho)^*(\epsilon)\in \rmH^2(D,\Q/\Z(1))$. Since $c\circ\rho=\pi\circ\psi$ and the action of $D$ on $\mmu_m$ coincides with its action on $\Im(\mu)$ via $\psi$, we obtain that $\rho^*((-1,-1)_K)=\psi^*(\pi^*(\epsilon))$, where $\pi^*(\epsilon)\in \rmH^2(N_{Q_{2^n}}(\Im(\mu)),\mmu_m^{-})$. The latter class is zero since the pullback of this extension along the projection $\pi$ is endowed with the diagonal section: hence, the desired claim follows. As a consequence, the subgroup generated by $(0,(-1,-1)_K)$ and $((-1,-1)_K,0)$ in $\rmH^2(K,Q_{2^n}^\ab(1))$ is contained in $\Sha^2_\nr(K,Q_{2^n}^\ab(1))$.

Let us now show that any class $\alpha=(a,b)\in \Sha^2_\nr(K,Q_{2^n}^\ab(1))$ must lie in the aforementioned subgroup. Let $N:=\langle r^{2^{n-3}},s\rangle\simeq Q_8$, let $\pi:N\to N/\langle s\rangle\simeq \Z/2$ be the natural projection and $c:\Gamma_K\to \Z/2$ be the character of the extension $K\subset K(i)$. Fix an epimorphism $\mu:\widehat{\Z}(1)\twoheadrightarrow \langle s\rangle $ and let $D:=\{ (\gamma,\sigma)\in N\times\Gamma_K~|~\pi(\gamma)=c(\sigma)\}$. One can check that $(\mu,D)\in {_\mathrm{any}}\mathscr{D}_\mathrm{{\mathrm{pro-}p'}}$. Moreover, $D$ is by definition the pullback of the extension 
\[1\longrightarrow \langle s\rangle\longrightarrow N\longrightarrow \Z/2\longrightarrow 1\]
along $c$. If we denote by $\epsilon\in \rmH^2(\Z/2,\langle s\rangle)$ the class of the latter extension then $c^*(\epsilon)$ yields an element of $\rmH^2(K,\langle s \rangle)\simeq \rmH^2(K,\mmu_4)=\Br(K)[4]$ (here $\Gamma_K$ acts on $\langle s\rangle$ by conjugation in $D$, i.e.\ by inversion through $c$) which coincides with $(-1,-1)_K$. An analysis of the Hochschild-Serre spectral sequence $E_2^{p,q}=\rmH^p(K,\rmH^q(\mmu_4,\Q/\Z(1)))\Rightarrow \rmH^{p+q}(D,\Q/\Z(1))$ shows that the kernel of the natural map $\rho^*:\Br(K) \to \rmH^2(D,\Q/\Z(1))$ is precisely the image of the transgression map $\rmH^0(K,\Hom(\mmu_4,\Q/\Z(1)))\to \Br(K)$, which is generated by $(-1,-1)_K$. Again, since $\overline{\mu}$ is $\Gamma_K$-invariant then $\rho^*(\alpha)\smile\overline{\mu}=\rho^*(\alpha\smile \overline{\mu})=\rho^*(b)\in \rmH^2(D,\Q/\Z(1))$. The condition that the latter class is zero hence forces $b$ to lie in the subgroup of $\Br(K)$ generated by $(-1,-1)_K$. If we now replace $\mu$ by an epimorphism $\widehat{\Z}(1)\twoheadrightarrow \langle rs\rangle$ (so that $N=\langle r^{2^{n-3}},rs\rangle\simeq Q_8$), we obtain that $a+b$ and hence $a$ must lie in the subgroup generated by $(-1,-1)_K$ as well. As a consequence, $\alpha\in \Sha^2_\nr(K,Q_{2^n}^\ab(1))$ must lie in the subgroup generated by $((-1,-1)_K,0)$ and $(0,(-1,-1)_K)$. 
\end{proof}

\begin{rmk}\label{remarqueSha2snrQ_16} In $Q_{2^n}/(\mathrm{conjugacy})$, the classes $[s]$ and $[sr]$ are always fixed under the cyclotomic action. Since $[s]=[sr^{2k}]$ and $[sr]=[sr^{2k+1}]$ for every $k\in \Z$ and given the equalities $s^2=r^{2^{n-2}}$ and $srs^{-1}=r^{-1}$, we have that $[s]=[s^3]$ and $[rs]=[(rs)^3]$: indeed, $s^3=sr^{2^{n-2}}$ and $(sr)^3=sr^{2^{n-2}+1}$. Suppose that $\alpha\in \Sha^2_\mathrm{snr}(K,\widehat{Q}^\ab_{2^n}(1))$ and let $L$ be a quasi-local field containing $K$. It follows that $\alpha$ must be orthogonal to the images of $r$ and $rs$ in $\rmH^0(L,Q_{2^n}^\ab(-1))=\langle\overline{r}\rangle\times\langle\overline{s}\rangle$, hence to the whole group $\rmH^0(L,Q^\ab_{2^n}(-1))$. In view of the equality $\rmH^2(L,\widehat{Q}_{2^n}^\ab(1))=\Br(L)[2]\oplus \Br(L)[2]$, we infer that the cup-product pairing $\rmH^2(L,\widehat{Q}_{2^n}^\ab(1))\times \rmH^2(L,Q_{2^n}^\ab(-1))\to \Br(L)$ is perfect. Thus, $\alpha$ dies in $\rmH^2(L,\widehat{Q}^\ab_{16}(1))$ for any such $L$, and the group $\Sha^2_\mathrm{snr}(K,\widehat{Q}^\ab_{2^n}(1))$ vanishes.
\end{rmk}

\begin{lem}\label{pasdecontributionduSha2nr} There is an equality 
\[\Ker\Big(\CH^2(X\otimes_\Q K)\to \CH^2(X\otimes_\Q \overline{K})\Big)=\Sha^2_\nr(K,\widehat{Q}_{2^n}^\ab(1))\]
of subgroups of $\rmH^3_\et(X\otimes_\Q K,\Q/\Z(2))$. 
\end{lem}
\begin{proof} Thanks to \cite[pp.\ 98-101]{browncohomology}, we know that the group $\rmH^2(Q_{2^n},\Q/\Z(2))$ is trivial, so that the equality $F^1\rmH^3_\mathrm{gnr}(X\otimes_\Q K,\Q/\Z(2))=F^2\rmH^3_\mathrm{gnr}(X\otimes_\Q K,\Q/\Z(2))$ holds. Therefore, we obtain an inclusion $\Ker\big(\CH^2(X\otimes_\Q K)\to \CH^2(X\otimes_\Q \overline{K})\big)\subset F^2\rmH^3_\mathrm{gnr}(X\otimes_\Q K ,\Q/\Z(2))$. As the image of $\CH^2(B_K Q_{2^n})$ in $F^2\rmH^3_\mathrm{gnr}(X\otimes_\Q K,\Q/\Z(2))$ does not meet $\rmH^3(K,\Q/\Z(2))$ non-trivially, we infer that $\Ker\big(\CH^2(X\otimes_\Q K)\to \CH^2(X\otimes_\Q \overline{K})\big)$ injects into $\Sha^2_\nr(K,\widehat{Q}_{2^n}^\ab(1))$. The aim from now on is to prove the reverse inclusion. 

By virtue of Lemma \ref{Sha2Q_16}, we may identify $\Sha^2_\nr(K,\widehat{Q}_{2^n}^\ab(1))$ with $\langle(-1,-1)_K\rangle\otimes \widehat{Q}_{2^n}^\ab$. Let us fix some $\alpha=(-1,-1)_K\otimes \xi\in \Sha^2_\nr(K,\widehat{Q}_{2^n}^\ab(1))$, where $\xi\in \widehat{Q}_{2^n}^\ab$ is viewed as a $1$-dimensional $K$-representation. Consider the cycle $z_\xi:=c_1(\xi)^2-\epsilon\cdot 2c_2(\psi)\in \CH^2(B_K Q_{2^n})$ where $\psi$ denotes the base change to $K$ of the irreducible faithful $\Q$\nobreakdash-representation of $Q_{2^n}$ of rank $2^{n-1}$ and such that $c_1(\psi)=0$. Here, $\epsilon$ is set to be equal to $0$ if $c_1(\xi)^2$ vanishes in $\CH^2(B_{\overline{K}}Q_{2^n})$ and equal to $1$ otherwise: this ensures that $z_\xi$ lies in $\Ker\big(\CH^2(B_K Q_{2^n})\to \CH^2(B_{\overline{K}} Q_{2^n})\big)$. We claim that $\alpha=\mathrm{cl}_{2,\mathrm{AJ}}(z_\xi)$. To prove this, we will first perform the computation after reducing modulo $2$, i.e.\ show that $\mathrm{cl}_{2,\mathrm{AJ}}(c_1(\xi)^2-\epsilon\cdot 2c_2(\psi))\equiv(-1,-1)_K\otimes \xi^2\text{ mod }2$: indeed, reduction mod $2$ being injective on $\rmH^2(K,\widehat{Q}_{2^n}^\ab(1))$ (as $\widehat{Q}_{2^n}^\ab$ is $2$-torsion), this will yield the desired equality. 

Note that $\mathrm{cl}_{2,\mathrm{AJ}}(c_1(\xi)^2-\epsilon\cdot 2c_2(\psi))\equiv \mathrm{cl}_{2,\mathrm{AJ}}(c_1(\xi)^2)\text{ mod }2$. The character $\xi:Q_{2^n}\to \mmu_2$ induces an evident morphism $f:X\otimes_\Q K\to (U\otimes_\Q K)/\mmu_2$ such that $\xi=f^*\mathcal{L}$ where $\mathcal{L}$ is the tautological line bundle, as well as a map $\CH^2(B_K\mmu_2)\to \CH^2(B_K Q_{2^n})$. Notably, we have $c_1(\xi)^2=f^*(c_1(\mathcal{L})^2)$, so it suffices to compute $\mathrm{cl}_{2,\mathrm{AJ}}(c_1(\mathcal{L})^2)$ on $(U\otimes_\Q K)/\mmu_2$ instead. If we view $[-1]\in \rmH^1(K,\mmu_2)$, we know that $\rmH^2_\et((U\otimes_\Q K)/\mmu_2,\mmu_2)\simeq \rmH^2(\mmu_2\times\Gamma_K,\mmu_2)$ naturally identifies with $\rmH^2(K,\mmu_2)\oplus \rmH^1(K,\mmu_2)x\oplus (\Z/2)x^2$ where $x\in \rmH^1_\et((U\otimes_\Q K)/\mmu_2,\mmu_2)$ is the class of the tautological $\mmu_2$-torsor $U\otimes_\Q K \to (U\otimes_\Q K)/\mmu_2$; see \cite[Theorem 6.10]{voevodskyreducedpowers} for details. Moreover, we have $\mathrm{cl}_{2,\mathrm{AJ}}(c_1(\mathcal{L}))=x^2+[-1]\smile x$ because this element is the Bockstein of the short exact sequence $0\to \mmu_2\to\mmu_4\to \mmu_2\to 0$ applied to $x$: indeed, by virtue of \emph{loc.\ cit., ibidem} the latter Bockstein is equal to $\mathrm{Sq}^1(x)+[-1]\smile x$ since $\Gamma_K$ acts on $\mmu_4$ by inversion through $[-1]$. Hence, we have an equality 
\[\mathrm{cl}_{2,\mathrm{AJ}}(c_1(\mathcal{L})^2)=(x^2+[-1]\smile x)^2=x^4+(-1,-1)_K \smile x^2.\]
Applying $f^*$, we obtain that $\mathrm{cl}_{2,\mathrm{AJ}}(c_1(\xi)^2)\equiv \xi^4+(-1,-1)_K\otimes \xi^2$. Since $\xi^4$ vanishes in $\rmH^4(Q_{2^n},\Z/2)$ \cite[pp.\ 98-101]{browncohomology}, we infer that $\mathrm{cl}_{2,\mathrm{AJ}}(c_1(\xi)^2)\equiv (-1,-1)_K\otimes \xi^2\text{ mod }2$, as desired.
\end{proof}

\begin{lem}\label{CH2Q_16} Let $E:= K(\zeta_{2^{n-1}}+\zeta_{2^{n-1}}^{-1})$. The cokernel of the base change map \[\CH^2(B_{K}Q_{2^n}) \to \CH^2(B_{\overline{K}}Q_{2^n})^{\Gamma_K}\] is non-zero if and only if $(-1,-1)_E\neq 0\in \Br(E)$. In the latter case, it is a group of order~$2$, generated by the image of $[E:K]\cdot c_2(\rho)$ where $\rho:Q_{2^n}\to \SL_2(\overline{K})$ denotes any irreducible and faithful $\overline{K}$-representation. 
\end{lem}
\begin{proof} Let $\rho:Q_{2^n}\to \SL_2(\overline{K})$ be the faithful $2$-dimensional representation given by 
\[\rho(r):=\begin{pmatrix}\zeta_{2^{n-1}} & 0 \\ 0 & \zeta_{2^{n-1}}^{-1}\end{pmatrix}~~\text{  and  }~~\rho(s):=\begin{pmatrix}0 & 1 \\ -1 & 0\end{pmatrix}.\]
The inclusion of the cyclic subgroup of order $2^{n-1}$ generated by $r$ induces a restriction map $\CH^2(B_{\overline{K}}Q_{2^n})\to \CH^2(B_{\overline{K}}\Z/2^{n-1})$. This map sends $c_2(\rho)$ to $c_2(\chi\oplus \chi^{-1})$ where $\chi$ denotes the faithful character $\Z/2^{n-1}\to \overline{K}^*$ that sends a generator to $\zeta_{2^{n-1}}$. The Whitney formula \cite[Theorem 3.2]{fulton} yields $c_2(\chi\oplus\chi^{-1})=-c_1(\chi)\smile c_1(\chi)\in \CH^2(B_{\overline{K}}\Z/2^{n-1})$, and the latter is a generator of $\CH^2(B_{\overline{K}}\Z/2^{n-1})$, see e.g.\ \cite[Theorem 2.10]{totaro}. 

Let $G_n:=\Gal(K(\zeta_{2^{n-1}})/K)$ and let $C_n$ be the image of $G_n$ in $(\Z/2^{n-1})^*/\langle \pm 1\rangle$ through the cyclotomic character (where $-1$ corresponds to the automorphism $\tau:\zeta_{2^{n-1}}\mapsto \zeta_{2^{n-1}}^{-1}$). We thus have an isomorphism $C_n\simeq\Gal(E/K)$. We see that for any $\sigma\in G_n$, the image of $c_2(^\sigma\rho)$ in $\CH^2(B_{\overline{K}}\Z/2^{n-1})$ is a generator. By virtue of the exact sequence from Proposition~\ref{suiteexactefondamentale}, the group $\CH^2(B_{\overline{K}}Q_{2^n})$ embeds into $\rmH^3(Q_{2^n},\Q/\Z(2))\simeq \mmu_{2^n}^{\otimes2}$; see \cite[pp.\ 98-101]{browncohomology} for the latter isomorphism. In view of the following commutative diagram
\[\begin{tikzcd}
	{\CH^2(B_{\overline{K}}Q_{2^n})} & {\rmH^3(Q_{2^n},\Q/\Z(2))} & {\mmu_{2^n}^{\otimes2}} \\
	{\CH^2(B_{\overline{K}}\Z/2^{n-1})} & {\rmH^3(\Z/2^{n-1},\Q/\Z(2))} & {\mmu_{2^{n-1}}^{\otimes2},}
	\arrow[hook, from=1-1, to=1-2]
	\arrow[from=1-1, to=2-1]
	\arrow["\sim", from=1-2, to=1-3]
	\arrow[from=1-2, to=2-2]
	\arrow[from=1-3, to=2-3]
	\arrow["\sim"', from=2-1, to=2-2]
	\arrow["\sim"', from=2-2, to=2-3]
\end{tikzcd}\]
we infer that since $c_2(^\sigma\rho)$ maps to a generator of $\CH^2(B_{\overline{K}}\Z/2^{n-1})$, then $\CH^2(B_{\overline{K}}Q_{2^n})$ must be generated by $c_2(^\sigma\rho)$ (notably, the map $\CH^2(B_{\overline{K}}Q_{2^n})\to\rmH^3(Q_{2^n},\Q/\Z(2))$ is an isomorphism). Let $v_K$ denote the $2$-adic valuation of $|\rmH^0(K,\Q/\Z(2))\{2\}|$ and $e_K:=\min(n,v_K)$. Taking $\Gamma_K$\nobreakdash-invariants, and since $[E:K]=|C_n|=2^{n-e_K}$, we deduce that $\CH^2(B_{\overline{K}}Q_{2^n})^{\Gamma_K}\simeq \mmu_{2^{e_K}}^{\otimes2}$ is generated by $[E:K]\cdot c_2(^\sigma\rho)$. Recall that the Schur index of a $\overline{K}$\nobreakdash-representation $V$ over $K$ is the multiplicity of this representation in the smallest $K$\nobreakdash-representation $W$ such that $W\otimes_K \overline{K}$ contains $V$; see \cite[§12.2]{serrereplin} for details. Since $^\sigma\rho$ has real character values and is not real valued, then it is quaternionic. Hence, its Schur index over $K$ is equal to $2$ if $(-1,-1)_E\neq 0\in \Br(E)$ and it is equal to $1$ otherwise, see \cite[Theorem 10.12 and Corollary 10.4]{isaacs}. Notably, in the first case we have $i\notin E$ so $\tau\in G_n$ and the corresponding irreducible $K$\nobreakdash-representation $\psi$ is $(2\cdot |G_n|)$\nobreakdash-dimensional and splits over $\overline{K}$ as 
\[\psi\otimes_{K}\overline{K}\simeq \bigoplus_{\sigma\in G_n} {^{\sigma}\rho}\simeq \bigoplus_{\sigma'\in C_n} ({^{\sigma'}\rho})^{\oplus 2}.\] 
In view of the following commutative diagram 
\[\begin{tikzcd}
	{F^2\mathrm{Rep}_K(Q_{2^n})} & {F^2\mathrm{Rep}_{\overline{K}}(Q_{2^n})} \\
	{\CH^2(B_{K}Q_{2^n})} & {\CH^2(B_{\overline{K}}Q_{2^n})}
	\arrow[from=1-1, to=1-2]
	\arrow["{c_2}"', twoheadrightarrow, from=1-1, to=2-1]
	\arrow["{c_2}", twoheadrightarrow, from=1-2, to=2-2]
	\arrow[from=2-1, to=2-2]
\end{tikzcd}\]
we infer that the image of the base change map $\CH^2(B_K Q_{2^n})\to\CH^2(B_{\overline{K}}Q_{2^n})$ must be generated by $c_2(\psi\otimes_K\overline{K})$. Taking Chern classes in the above isomorphism, and since ${^{\sigma}\rho}$ has trivial determinant for every $\sigma\in G_n$, the Whitney formula shows that $c_2(\psi\otimes_{K}\overline{K})=\sum_{\sigma\in G_n}c_2({^{\sigma}\rho})=2\cdot \sum_{\sigma'\in C_n}c_2({^{\sigma'}\rho})$. 

As each $c_2({^{\sigma}\rho})$ is a generator of $\CH^2(B_{\overline{K}}Q_{2^n})$ and $|C_n|=[E:K]$ we may apply the previous calculations to show that $\sum_{\sigma'\in C_n}c_2({^{\sigma'}\rho})$ is a generator of the $\Gamma_K$-invariant subgroup of $\CH^2(B_{\overline{K}}Q_{2^n})$. It follows that $c_2(\psi\otimes_{K}\overline{K})$ generates a subgroup of index~$2$ in $\CH^2(B_{\overline{K}}Q_{2^n})^{\Gamma_K}=\langle [E:K]\cdot c_2(\rho)\rangle\simeq\mmu_{2^{e_K}}^{\otimes 2}$. Therefore, the cokernel of the map $\CH^2(B_{K}Q_{2^n})\to\CH^2(B_{\overline{K}}Q_{2^n})^{\Gamma_K}$ has order $2$, as desired. By \cite[Lemma 2.8, (ii)]{bailleul}, any irreducible and faithful $\overline{K}$-representation of $Q_{2^n}$ is isomorphic to $^\sigma\rho$ for some $\sigma \in \Gal(\Q(\zeta_{2^{n-1}})/\Q)$. Since $\langle[E:K]\cdot c_2(\rho)\rangle= \langle [E:K]\cdot c_2({^\sigma}\rho)\rangle$ for every $\sigma \in \Gal(\Q(\zeta_{2^{n-1}})/\Q)$, the second part of the statement of the lemma follows. 

If otherwise the Schur index of $\rho$ is equal to $1$, i.e.\ if $(-1,-1)_E=0$ in $\Br(E)$, then $\psi\otimes_{K}\overline{K}\simeq\bigoplus_{\sigma'\in C_n} {^{\sigma'}\rho}$ and $c_2(\psi\otimes_K\overline{K})=\sum_{\sigma'\in C_n}c_2(^{\sigma'}\rho)$ generates the subgroup of $\Gamma_K$\nobreakdash-invariants in $\CH^2(B_{\overline{K}}Q_{2^n})$. Hence, the cokernel of $\CH^2(B_{K}Q_{2^n}) \to \CH^2(B_{\overline{K}}Q_{2^n})^{\Gamma_K}$ is zero.
\end{proof}

\begin{rmk}\label{remarquecorpsdescaractèresCH2} (i) Let $K\subset \overline{\Q}$ be the subfield generated by the character values of any irreducible and faithful representation $\rho:Q_{2^n}\hookrightarrow \SL_2(\overline{\Q})$. It is known that $K$ coincides with $\Q(\zeta_{2^{n-1}}+\zeta_{2^{n-1}}^{-1})$, see \cite[Lemma 2.8 and Lemma 2.9]{bailleul}; the latter is precisely the subfield of $\Q(\zeta_{2^{n-1}})$ that is fixed by the complex conjugation \cite[pp. 15--16]{washington}. Hence, it is Galois over $\Q$ and its Galois group corresponds to the quotient $C_n$ of $\Gal(\Q(\zeta_{2^{n-1}})/\Q)$. As in the proof of Lemma~\ref{CH2Q_16}, we infer that over this field, the $\Q$-representation $\psi$ splits as $\psi\otimes_\Q K\simeq \bigoplus_{\sigma'\in C_n}\psi_{\sigma'}$ where each $\psi_{\sigma'}$ is a $4$-dimensional $K$-representation such that $\psi\otimes_K\overline{\Q}\simeq (^{\sigma'}\rho)^{\oplus 2}$ (the Schur index is by definition equal to $2$ over $K$, cf.\ \cite[§12.2]{serrereplin}). One can check that $\mmu_{2^n}^{\otimes 2}$ is fixed by $\Gamma_K$, so that $\CH^2(B_{\overline{\Q}}Q_{2^n})^{\Gamma_K}=\CH^2(B_{\overline{\Q}}Q_{2^n})$ is generated by $c_2(\rho)$. The image of the restriction map $\CH^2(B_K Q_{2^n})\to\CH^2(B_{\overline{\Q}}Q_{2^n})$ is generated by $2c_2(^{\sigma'}\rho)$ for any $\sigma'\in C_n$ and thus the cokernel of this map has order $2$ as well.

(ii) If $K=\Q(i)$, then the Schur indices of all the faithful representations of $Q_{2^n}$ drop to $1$. Consequently, the image of the restriction map $\CH^2(B_K Q_{2^n})\to\CH^2(B_{\overline{\Q}}Q_{2^n})$ is generated by $\sum_{\sigma'\in C_n}c_2({^{\sigma'}\rho})$, the latter of which has order $8$ in $\CH^2(B_{\overline{\Q}}Q_{2^n})$. Since the subgroup of $\mmu_{2^n}^{\otimes 2}$ fixed by $\Gamma_K$ is $\mmu_8^{\otimes 2}$, this implies that the cokernel of $\CH^2(B_K Q_{2^n})\to\CH^2(B_{\overline{\Q}}Q_{2^n})^{\Gamma_K}$ is zero in this case.  
\end{rmk} 
\subsubsection{An unramified cohomology class via the Rost invariant}\label{invariantdeRost} Assume that $(-1,-1)_E\neq 0$ in $\Br(E)$, so that $i\notin E$ and $E(\zeta_{2^{n-1}})=E(i)$. Let $A$ be the quaternion algebra corresponding to $(-1,-1)_E$ and $I,J\in A$ be generators. Identifying $E(\zeta_{2^{n-1}})$ with $E(I)$, we infer that the mapping $r\mapsto \zeta_{2^{n-1}}$, $s\mapsto J$ defines an embedding of $K$-groups $\iota:Q_{2^n}\hookrightarrow \SL_1(A)(E)=H(K)$, where $H:=\mathrm{Res}_{E/K}(\SL_1(A))$. We let $\mathcal{R}:=\iota^*R_H\in \rmH^3_\et(X\otimes_\Q K,\Q/\Z(2))$ where $R_H$ denotes the Rost invariant of $H$ in the sense of \cite[Part II, §9]{garibaldimerkurjevserre}. More precisely, let $T$ be a $Q_{2^n}$-torsor over a field $L\supset K$, given by $\varphi:\Gamma_L\to Q_{2^n}$, and write $L\otimes_KE=\prod_jL_j$ where the $L_j$'s are fields containing $L$. Then $\rmH^1(L,H)=\prod_j\rmH^1(L_j,\SL_1(A_{L_j}))$, and $\rmH^1(L_j,\SL_1(A_{L_j}))\simeq L_j^*/\mathrm{Nrd}(A_{L_j}^*)$ by virtue of Hilbert's Theorem 90 for $\GL_1(A)$ (here $\mathrm{Nrd}:A^*\to E^*$ denotes the reduced norm, see e.g.\ \cite[§1.5]{CTS21} for details). Concretely, choose $b\in(A\otimes_E\overline L)^*$ such that $\iota(\varphi(\sigma))=b^{-1}\sigma(b)$ for $\sigma\in\Gamma_{L_j}$, and let $\lambda_{T,j}:=\mathrm{Nrd}(b)$, which lies in $L_j^*$ and is well defined modulo $\mathrm{Nrd}(A_{L_j}^*)$. By \cite[Part~II, §11]{garibaldimerkurjevserre}, the Rost invariant of $\SL_1(A_{L_j})$ is, up to a sign, equal to $\lambda\mapsto[\lambda]\smile[A_{L_j}]$ where $[\lambda]$ is viewed as an element of $\rmH^1(L_j,\widehat{\Z}(1))$. Hence,
\[\mathcal{R}(T)=\pm\sum_j\mathrm{Cores}_{L_j/L}\Big([\lambda_{T,j}]\smile[A_{L_j}]\Big)\in \rmH^3(L,\Q/\Z(2)).\]
Thanks to Serre's criterion \cite[Part~I, \S33]{garibaldimerkurjevserre}, it suffices to compute for any field $L\supset K$ and any  $Q_{2^n}$-torsor $T$ over $F:=L((t))$ given by $\varphi:\Gamma_F\cong\widehat{\Z}(1)\rtimes\Gamma_L\to Q_{2^n}$, the image of $\mathcal{R}(T)$ along the residue map $\partial_F:\rmH^3(F,\Q/\Z(2))\to \Br(L)$. Let $\tau\in Q_{2^n}$ be the image of a generator $\zeta$ of the inertia subgroup of $\Gamma_F$. The factors $F_j=L_j((t))$ of $F\otimes_KE$ are unramified over $F$ and $[A]$ is unramified, so
\[
\partial_F(\mathcal{R}(T))=\pm\sum_j\mathrm{Cores}_{L_j/L}\big(v(\lambda_{T,j})\,[A_{L_j}]\big),
\]
where $v$ denotes the $t$-adic valuation on $F_j$. Therefore, everything reduces to computing $v(\lambda_{T,j})=v(\mathrm{Nrd}(b))$, where $b$ trivialises $\iota\circ\varphi$ over $F_j$. In each case below, we construct $b$ explicitly. Fix $j$ and write $F$ for $F_j$ and $L$ for $L_j$. Recall that in $\Gamma_F\cong\widehat{\Z}(1)\rtimes\Gamma_L$, we have $(1,\sigma)(\zeta,1)(1,\sigma)^{-1}=(\zeta^{\chi(\sigma)},1)$, so that $\varphi(\sigma)\,\tau\,\varphi(\sigma)^{-1}=\tau^{\chi(\sigma)}$ in $Q_{2^n}$.
\begin{lem}\label{Lemmenormeréduite}
Let $F\subset F'$ be a finite unramified Galois extension and $\psi$ a $1$-cocycle of $\Gal(F'/F)$ with values in $(A\otimes_E\mathscr O_{F'})^*$. Then $\psi(\sigma)=a^{-1}\sigma(a)$ for some $a\in(A\otimes_E\mathscr O_{F'})^*$. In particular, $\mathrm{Nrd}(a)\in\mathscr O_{F'}^*$.
\end{lem}

\begin{proof}
Since $E\subset\mathscr{O}_F$, we see that $A\otimes_E\mathscr{O}_F$ is an Azumaya algebra over $\mathscr{O}_F$, so that $\GL_1(A)$ is a smooth group scheme over $\mathscr O_F$. By virtue of Hensel's lemma \cite[Exposé XXIV, Proposition 8.1]{SGA3}, the natural map $\rmH^1(\mathscr{O}_F,\GL_1(A))\to \rmH^1(L,\GL_1(A_L))$ is injective, and the latter group vanishes thanks to Hilbert's Theorem~90.\end{proof}

\begin{prop}\label{prop:res} If $\tau$ is the non-trivial central element, then $v(\lambda_{T,j})=1$ for all $j$, so $\partial_F(\mathcal{R}(T))=\pm[E:K]\cdot(-1,-1)_L$. If $\tau$ is not central, then $v(\lambda_{T,j})=0$ for all $j$, so $\partial_F(\mathcal{R}(T))=0$.
\end{prop}

\begin{proof} If $\tau=-1$, let $\varphi_0:=\varphi|_{\Gamma_L}$. Viewed as a homomorphism $\Gamma_F\to\Gamma_L\to\mathrm Q_{2^n}$, it is trivial on the inertia subgroup, so by Lemma \ref{Lemmenormeréduite} there exists a finite unramified extension $F\subset F'$ and some $a\in(A\otimes_{E}\mathscr{O}_{F'})^*$ such that $\iota(\varphi_0(\sigma))=a^{-1}\,{}^{\sigma}\!a$ for all $\sigma\in\Gamma_L$. In particular, $\mathrm{Nrd}(a)$ is a unit. If we let $b:=a\,t^{1/2}$, then $t^{1/2}$ is central in $A\otimes_{E}\overline F$ and fixed by $\Gamma_L$, so $b^{-1}\,{}^{\sigma}b=\iota(\varphi(\sigma))$ for all $\sigma\in\Gamma_L$. Moreover, since $\zeta$ fixes $a$ and ${}^{\zeta}(t^{1/2})=-t^{1/2}$, we get ${}^{\zeta}b=a\cdot(-t^{1/2})=b\cdot(-1)=b\,\iota(\varphi(\zeta))$. Since $\zeta$ and $\Gamma_L$ topologically generate $\Gamma_F$, we infer that $\iota(\varphi(\sigma))=b^{-1}\,{}^{\sigma}b$ for all $\sigma\in\Gamma_F$. On the other hand the reduced norm of a scalar $x$ is $x^2$, so $\lambda=\mathrm{Nrd}(a)\,t$ and $v(\lambda)=1$. Summing over $j$ and using the fact that $[A_{L_j}]=\mathrm{Res}_{L_j/L}(-1,-1)_L$, we infer that the residue is equal to $\mathrm{Cores}_{L\otimes E/L}(\mathrm{Res}_{L\otimes E/L}((-1,-1)_L))=[E:K]\cdot(-1,-1)_L\in\Br(L)$.

If now $\tau$ is not central, then it must have order $c\ge4$, and some power of $\tau$ squares to $-1$. Hence, identifying $Q_{2^n}$ with its image under $\iota$, we see that $E':=E(\tau)\subset A$ is a quadratic extension of $E$ isomorphic to $E(i)$. As it is a maximal subfield of $A$, the centraliser of $\tau$ in $A$ coincides with $E'$. By the semi-commutation relation $\varphi(\sigma)\tau\varphi(\sigma)^{-1}=\tau^{\chi(\sigma)}$, we see that $\varphi(\sigma)\in(E')^*$ if and only if $\chi(\sigma)\equiv1\bmod c$. Otherwise, $\chi(\sigma)\equiv-1\bmod c$ and $\varphi(\sigma)$ acts on $E'$ through the only non-trivial automorphism. Therefore, the subgroup $\Gamma_0:=\{\sigma\in\Gamma_F\mid\chi(\sigma)\equiv1\bmod c\}$ of $\Gamma_F$ is normal of index $\le2$ and contains $\widehat{\Z}(1)$, and its fixed field $F''=F(\zeta_c)$ is unramified over $F$. The restriction of $\phi$ to $\Gamma_0$ is a homomorphism to the norm-one roots of unity of $E'$. Since $i\in F''$, we have $E'\otimes_{E}F''\simeq F''\times F''$, and $\varphi|_{\Gamma_0}=(\chi_1,\chi_1^{-1})$ for some character $\chi_1:\Gamma_0\to\mu_\infty(F'')$. By Kummer theory, $\chi_1(\sigma)={}^{\sigma}\beta/\beta$ with $\beta=w^{1/m}$ and $w\in F''^*$. Put $b_0:=(\beta,\beta^{-1})\in(E'\otimes\overline{F})^*$. Then $b_0^{-1}\,{}^{\sigma}b_0=\varphi(\sigma)$ for $\sigma\in\Gamma_0$, and $\mathrm{Nrd}(b_0)=\beta\beta^{-1}=1$. Let $\varphi'(\sigma):=b_0\,\varphi(\sigma)\,({}^{\sigma}b_0)^{-1}$, which is a cocycle cohomologous to $\varphi$. It is trivial on the normal subgroup $\Gamma_0$, so it factors through $\Gal(F''/F)$, with values in $(A\otimes F'')^*$. If $F''=F$, then $\varphi'=1$. Otherwise, for $\sigma\notin\Gamma_0$, write $\varphi(\sigma)=\xi h$, with $h\in Q_{2^n}$ a fixed element normalising $E'$ and $\xi$ a root of unity in $E'$. In the coordinates $E'\otimes\overline{F}=\overline{F}\times\overline{F}$, both $h$ and $\sigma$ exchange the two factors, and a direct computation gives $\varphi'(\sigma)=(\xi_1u_\sigma^{-1},\ \xi_2u_\sigma)\,h$, where $u_\sigma:={}^{\sigma}\beta/\beta$. Since $\sigma$ preserves the valuation $v$, we see that $u_\sigma$ is a unit, and so $\varphi'$ takes values in $(A\otimes\mathscr{O}_{F''})^*$. By Lemma \ref{Lemmenormeréduite}, we get that $\varphi'(\sigma)=a^{-1}\,{}^{\sigma}a$ for some $a$ integral. Thus $b:=a\,b_0$ satisfies $\varphi(\sigma)=b^{-1}\,{}^{\sigma}b$ for all $\sigma\in\Gamma_F$, and $\lambda=\mathrm{Nrd}(a)\mathrm{Nrd}(b_0)=\mathrm{Nrd}(a)$ is a unit. This shows that $v(\lambda)=0$.\end{proof}
\begin{lem}\label{imagerostCH2}
The image of $\mathcal{R}$ in $\mathrm{Coker}(\CH^2(X\otimes_\Q K)\to \CH^2(X\otimes_\Q\overline{K}))$ is equal to the class of $\pm\big(\sum_{\sigma\in\Gal(E/K)}\chi(\sigma)^2\big)c_2(\rho)$, which is an odd multiple of $[E:K]\cdot c_2(\rho)$.
\end{lem}

\begin{proof} The Rost invariant of $\SL_{2,K}$ corresponds to $\pm c_2(\eta)$ where $\eta:\SL_{2,K}\to \GL_{2,K}$ denotes the standard representation, see \cite[Example 3.5 and Theorem 3.9]{merkurjevcohomologicalinvariants}. As the restriction to $\overline{K}$ of a corestriction is the sum of the Galois conjugates, this hence gives the formula. If we write the image of $\Gamma_K$ in $(\Z/2^{n-1})^*$ as $\{\pm1\}\times(1+2^{n-1-s}\Z)$ where $2^s=[E:K]$ and $s\le n-3$, then $\sum_{x<2^s}(1+2^{n-1-s}x)^2\equiv2^s(1+2^{n-1-s})\pmod{2^n}$ if $s\ge1$, and the sum is equal to $1$ if $s=0$.
\end{proof}
\begin{proof}[Proof of Theorem \ref{Q_16}] It is well-known that $X\otimes_\Q \overline{K}$ is rational, see e.g.\ \cite[Theorem~3.1]{hukangpgroups}, so that $\rmH^{3,\norm}_\nr(K(X),\Q/\Z(2))=\rmH^{3,\norm}_{\nr,\mathrm{alg}}(K(X),\Q/\Z(2))$. Thanks to  \cite[pp.\ 98-101]{browncohomology}, we know that the group $\rmH^2(Q_{2^n},\Q/\Z(2))$ is trivial, so that $F^1\rmH^3_\mathrm{gnr}(X\otimes_\Q K,\Q/\Z(2))=F^2\rmH^3_\mathrm{gnr}(X\otimes_\Q K,\Q/\Z(2))$. 
By Proposition \ref{suiteexactefondamentale} and Lemma \ref{pasdecontributionduSha2nr}, we see that the latter group has trivial image in $\rmH^{3,\norm}_\nr(K(X),\Q/\Z(2))$. On the other hand, Lemma~\ref{CH2Q_16} shows that $\Coker\big(\CH^2(X\otimes_\Q K)\to \CH^2(X\otimes_\Q \overline{K})^{\Gamma_K}\big)$ is trivial if $(-1,-1)_E= 0\in \Br(E)$, otherwise it has order $2$ and is generated by the image of $\beta:=[E:K]\cdot c_2(\rho)$. In the first case, the previous observations directly imply that $\rmH^{3,\norm}_\nr(K(X),\Q/\Z(2))$ vanishes.

Suppose now that $(-1,-1)_E\neq 0$ in $\Br(E)$. If $K\neq E$, then $[E:K]$ is even, and by virtue of Proposition \ref{prop:res}, all the residue values of $\mathcal{R}$ must vanish. Hence, Serre's criterion implies that the image of $\mathcal{R}$ in $\rmH^3(K(X),\Q/\Z(2))$ is unramified. Moreover Lemma \ref{imagerostCH2} shows that the image of $\mathcal{R}$ in $\Coker\big(\CH^2(X\otimes_\Q K)\to \CH^2(X\otimes_\Q \overline{K})^{\Gamma_K}\big)$ coincides with that of $\beta$, so its image in $\rmH^{3,\norm}_\nr(K(X),\Q/\Z(2))$ is non-zero.

If otherwise $K=E$, then $\beta=c_2(\rho)$. If we let $Z:=\langle r^{2^{n-2}}\rangle$ be the center of $Q_{2^n}$, $\mu:\widehat{\Z}(1)\to Z$ be a surjective morphism and $D:=Z\times\Gamma'$ with $\Gamma'\subset\Gamma_K$ any pro-$2$-Sylow subgroup, then $(\mu,D)\in{_\mathrm{any}}\mathscr{D}_\mathrm{{\mathrm{pro-}p'}}$ and the associated residue $\partial_{\mu,D}(\mathcal{R})\in \rmH^2(D,\Q/\Z(1))$ is the image of $(-1,-1)_K\neq0$ along the composite map $\Br(K)\to \Br(\overline{K}^{\Gamma'})\to \rmH^2(D,\Q/\Z(1))$. The latter element is non-zero because the natural projection $D\to\Gamma'$ admits a splitting and $\Br(K)[2]\to\Br(\overline{K}^{\Gamma'})[2]$ is injective. Hence, $\mathcal{R}$ is not globally unramified in this case. Moreover, any other lift of the generator of $\Coker\big(\CH^2(X\otimes_\Q K)\to \CH^2(X\otimes_\Q \overline{K})^{\Gamma_K}\big)$ differs from $\pm\mathcal R$ by $\mathrm{cl}_{2,\mathrm{AJ}}(z)+\alpha$ for some $z\in\CH^2(X\otimes_\Q K)$ and some $\alpha\in F^2\rmH^3_\et(X\otimes_\Q K,\Q/\Z(2))$. Such classes have trivial residue at this pair, since $\mathrm{cl}_{2,\mathrm{AJ}}(z)$ is always globally unramified and $\bar\mu=0\in \rmH^0(D,Q_{2^n}^\ab(-1))$ as $Z\subset[\mathrm Q_{2^n},\mathrm Q_{2^n}]$. Hence, no lift is globally unramified, and $H^{3,\mathrm{nor}}_{\mathrm{nr}}(K(X),\mathbf Q/\mathbf Z(2))=0$.
\end{proof}

\begin{rmk} For $K=\Q$ and $n=3$, i.e.\ $Q_{2^n}=Q_8$, the last part of the proof of Theorem~\ref{Q_16} notably explains the vanishing of $\rmH^{3,\norm}_\nr(\Q(X),\Q/\Z(2))$, consistent with the classical result of Gröbner \cite{gröbner}. Indeed, in this case we have $\Q=\Q(\zeta_4+\zeta_4^{-1})=E$ and $[E:\Q]=1$. The class $\beta=c_2(\rho)$ descends to $\Q$ but the central residue is non-zero since $(-1,-1)_\Q\neq 0$ in $\Br(\Q)$: taking for instance $Z:=\langle r^2\rangle$, $\mu:\widehat{\Z}(1)\to Z$ surjective and $D:=Z\times\Gal(\C/\R)$, the associated residue value $\partial_{\mu,D}(\mathcal{R})$ is the non-zero class $(-1,-1)_{\R}\in\Br(\R)$.
\end{rmk}

\begin{rmk} In \cite[§33.28]{garibaldimerkurjevserre}, Serre explains how to use Rost's invariant in order to construct a degree $3$ unramified cohomology class for $Q_{16}$. Let us compare our construction with his. Let $c\ge 3$, fix a classifying pair $(U\subset \A^m_\Q,c)$ for $\mathrm{Spin_{7,\Q}}$ and let $Y:=U/\mathrm{Spin_{7,\Q}}$. Consider an embedding $Q_{16}\hookrightarrow  \widetilde{\mathfrak{A}}_7\hookrightarrow \mathrm{Spin_{7,\Q}}$ where $\widetilde{\mathfrak{A}}_7$ is the unique non-trivial $2$-cover of the alternating group $\mathfrak{A}_7$. Since $(U\subset \A^m_\Q,c)$ is also a classifying pair for $\widetilde{\mathfrak{A}}_7$ and $Q_{16}$, we may consider $Y':=U/\widetilde{\mathfrak{A}}_7$ and $X:=U/Q_{16}$. The evident morphism $X\to Y$ induces a commutative diagram:
\[\small\begin{tikzcd}
	0 & {\rmH^3(\Q,\Q/\Z(2))} & {\rmH^0(Y,\mathscr{H}_\et^3(\Q/\Z(2)))} & {\Coker\Big(\CH^2(Y)\to \CH^2(\overline{Y})^{\Gamma_\Q}\Big)} & 0 \\
	0 & {\rmH^3(\Q,\Q/\Z(2))} & {\rmH^0(X,\mathscr{H}_\et^3(\Q/\Z(2)))} & {\Coker\Big(\CH^2(X)\to \CH^2(\overline{X})^{\Gamma_\Q}\Big),}
	\arrow[from=1-1, to=1-2]
	\arrow[from=1-2, to=1-3]
	\arrow[equals, from=1-2, to=2-2]
	\arrow[from=1-3, to=1-4]
	\arrow[from=1-3, to=2-3]
	\arrow[from=1-4, to=1-5]
	\arrow[from=1-4, to=2-4]
	\arrow[from=2-1, to=2-2]
	\arrow[from=2-2, to=2-3]
	\arrow[from=2-3, to=2-4]
\end{tikzcd}\]
where the top row is exact and the lower one is a complex. The surjectivity of the rightmost map in the top row holds by virtue of \cite[Proposition 9.5 and Theorem~9.11]{garibaldimerkurjevserre}. Moreover, the generator of this group is Rost's invariant, see \emph{loc.\ cit.}, \emph{ibidem}. Serre explains in \cite[§33.28]{garibaldimerkurjevserre} that the restriction of this element to $\rmH^0(Y',\mathscr{H}^3_\et(\Q/\Z(2))$ is in this case unramified and proves in \emph{loc.\ cit.}, Theorem 33.26 that its further restriction to $\rmH^0(X,\mathscr{H}^3_\et(\Q/\Z(2)))$ is non-zero. It turns out that the restriction of the generator of $\Coker\big(\CH^2(Y)\to \CH^2(\overline{Y})^{\Gamma_{\Q}}\big)$ to $\Coker\big(\CH^2(X)\to \CH^2(\overline{X})^{\Gamma_{\Q}}\big)$ coincides with the image of $2c_2(\rho)\in\CH^2(\overline{X})^{\Gamma_{\Q}}$, because the unique central element of order $2$ of $Q_{16}$ maps identically to that of $\mathrm{Spin}_{7,\Q}$. I thank Pierre Guillot for pointing out this argument to me.

\end{rmk}

\begin{rmk}\label{remarquequaternionsobstructionnoether} We discuss two cases in which the unramified cohomology class appearing in Theorem \ref{Q_16} vanishes or ramifies, and relate them to known results in the literature.

(i) Let $K$ be either $\Q(\sqrt{2})$ or $\Q(i)$. In \cite[Theorem 1.4]{kangnoetherdihedral}, the author proves that $X$ becomes rational over $K$, because $\Q(\zeta_8)$ is a cyclic extension of $K$. The vanishing of the obstruction lying in $\rmH^3_\nr(\Q(X),\Q/\Z(2))$ along the natural restriction map $\rmH^3_\nr(\Q(X),\Q/\Z(2))\to \rmH^3_\nr(K(X),\Q/\Z(2))$ appears explicitly in this case. Indeed, if $K=\Q(i)$, then Remark \ref{remarquecorpsdescaractèresCH2}, (ii) implies that the cokernel of $\CH^2(B_KQ_{16}) \to \CH^2(B_{\overline{\Q}}Q_{16})^{\Gamma_\Q}$ is trivial. If $K=\Q(\sqrt{2})$, then Remark \ref{remarquecorpsdescaractèresCH2}, (i) shows that this cokernel has order $2$ and $\CH^2(B_{\overline{\Q}}Q_{16})^{\Gamma_\Q}=\CH^2(B_{\overline{\Q}}Q_{16})$ is generated by $c_2(\rho)$ for any irreducible and faithful $\overline{\Q}$\nobreakdash-representation $\rho:Q_{16}\hookrightarrow \SL_2(\overline{\Q})$. The same line of arguments as those given in the beginning of the proof of Theorem \ref{Q_16} shows that this Chern class descends to $K$, however the rest of the proof fails since the restriction of $c_2(\rho)$ to a proper subgroup is no longer zero. This highlights the fact that this class is ramified and hence cannot contribute to $\rmH^3_\nr(K(X),\Q/\Z(2))$.

(ii) More generally, Remark \ref{remarquecorpsdescaractèresCH2}, (ii) shows that the map $\CH^2(B_K Q_{2^n})\to \CH^2(B_{\overline{\Q}} Q_{2^n})^{\Gamma_K}$ is onto as soon as $(-1,-1)_E=0$ in $\Br(E)$, e.g.\ when $i\in K$. This is consistent with the known positive results on Noether’s problem for generalised quaternion groups over extensions of $\Q$, see e.g.\ \cite[Theorem~3.2]{hukangpgroups}, which require $K$ to contain $\Q(\zeta_{2^{n-2}})$; it suggests that these results might hold under weaker assumptions, for instance over $K=\Q(i)$, where degree three unramified cohomology yields no obstruction. On the other hand, Remark \ref{remarquecorpsdescaractèresCH2}, (i) shows that if $K$ is the subfield of $\overline{\Q}$ generated by the character values of an irreducible faithful representation $\rho:Q_{2^n}\hookrightarrow\mathrm{SL}_2(\bar{\mathbf Q})$, then the above cokernel has order $2$ and is generated by $c_2(\rho)$. As explained above in the case of $Q_{16}$, this class descends to $K$ but it is ramified by virtue of the proof of Theorem \ref{Q_16} (the central residue being non-zero for this class), so that unramified cohomology of degree three yields no obstruction. It would thus be interesting to investigate whether Noether’s problem for $Q_{2^n}$ holds over this field.
\end{rmk}

\section{$p$-adic function fields}\label{paragraphehomogène}

In this section, we discuss the case of homogeneous spaces of $\SL_{n,k}$ when $k$ is the function field of a $p$-adic curve. Notably, we compare our formulas for the unramified cohomology in degree three with those appearing in work of Harari-Scheiderer-Szamuely \cite{hararischeidererszamuely}.

\subsection{Homogeneous spaces over $p$-adic function fields}  Let $C$ be a smooth, proper and geometrically integral curve over a $p$-adic field, i.e.\ a finite extension of $\Q_p$. We denote by $C^{(1)}$ the set of closed points of $C$ and by $k$ its function field. If $v\in C^{(1)}$, we denote by $k_v$ the completion of $k$ with respect to the induced discrete valuation, resp.\ $\kappa_v$ its residue field (which coincides with the residue field of $\mathscr{O}_{C,v}$). It follows that $k_v$ is a $2$-dimensional local field, that is, a complete discretely valued field whose residue field is a classical non-archimedean local field. By \cite[§II.4.3]{serrecg}, we know that there exists an isomorphism $k_v\simeq \kappa_v((t))$, where the latter denotes the field of Laurent series in one variable over $\kappa_v$. In particular, if we denote by $k_v^\nr$ the isomorphic image of $\kappa_v((t^{1/\infty}))$ in $k_v$, then the inclusion of $\Gamma_{k_v^\nr}$ into $\Gamma_k$ induces a splitting of the following short exact sequence:
\[1\longrightarrow I_v\longrightarrow \Gamma_{k_v}\longrightarrow \Gamma_{\kappa_v}\longrightarrow 1,\]
where $I_v$ denotes the inertia subgroup of $\Gamma_{k_v}$. In particular, for any finite $\Gamma_{k_v}$-module $M$ on which $I_v$ acts trivially, the Hochschild-Serre spectral sequence gives rise to a residue map $\rmH^3(k_v,M)\to \rmH^2(\kappa_v,M(-1))$ which is an isomorphism. In particular, there is a local invariant map:
\[\mathrm{inv}_v:\rmH^3(k_v,\Q/\Z(2))\overset{\sim}{\longrightarrow}\Q/\Z.\] 
Subsequently, there exist two split short exact sequences
\begin{equation}\label{scindagepadique1}0\longrightarrow \rmH^1(\kappa_v,M)\longrightarrow\rmH^1(k_v,M)\longrightarrow \rmH^0(\kappa_v,M(-1))\longrightarrow 0\end{equation}

and
\begin{equation}\label{scindagepadique2}0\longrightarrow \rmH^2(\kappa_v,M)\longrightarrow \rmH^2(k_v,M)\longrightarrow \rmH^1(\kappa_v,M(-1))\longrightarrow 0.\end{equation}
Moreover, if we denote $M':=\Hom(M,\Q/\Z(2))$, then the cup-product pairings 
\[\rmH^i(k_v,M)\times \rmH^{3-i}(k_v,M')\longrightarrow \rmH^3(k_v,\Q/\Z(2))\simeq\Q/\Z\]
are perfect for $0\le i \le 3$, see \cite[I.2.17]{milneadt}. We denote by $\Sha^2_\omega(k,M)$ the subgroup of $\rmH^2(k,M)$ consisting of those classes which vanish along the natural restriction map $\rmH^2(k,M)\to \rmH^2(k_v,M)$ for all but finitely many points $v\in C^{(1)}$. If $G$ is a finite $k$-group, and $v\in C^{(1)}$, we denote by $\mathcal{N}_v(G)$ the subgroup of $\rmH^0(\kappa_v,G^\ab(-1))$ consisting of the images of those elements $\tau\in G(\overline{k})$ which satisfy the equality ${^\sigma}[\tau]=[\tau^{\chi(\sigma)}]\in G(\overline{k})/(\mathrm{conjugacy})$ for any $\sigma\in \Gamma_{\kappa_v}$.
\begin{thm}\label{pointslocauxcorpsdefonctions}Let $G$ be a finite $k$-group. If a class lies in $\Sha^2_\nr(k,\widehat{G}^\ab(1))$, then for almost every point $v\in C^{(1)}$, its image in $\rmH^2(k_v,\widehat{G}^\ab(1))$ is orthogonal, with respect to the cup-product pairing \[\rmH^2(k_v,\widehat{G}^\ab(1))\times\rmH^1(k_v,G^\ab)\to \Q/\Z,\] 
to the inverse image of $\mathcal{N}_v(G)$ in $\rmH^1(k_v,G^\ab)$.
\end{thm}
\begin{proof} Let $\alpha\in \Sha^2_\nr(k,\widehat{G}^\ab(1))$. We can find a suitable Zariski open subset $U\subset C$, a finite étale group scheme $\mathscr{G}$ over $U$ and a class $\beta\in \rmH^2_\et(U,\widehat{\mathscr{G}}^\ab(1))$ whose image in $\rmH^2(k,\widehat{G}^\ab(1))$ is $\alpha$; see \cite[Chapter III, Lemma 1.16]{milneetale}. It follows that for any point $v\in U$, the image of $\alpha$ along $\rmH^2(k,\widehat{G}^\ab(1))\to \rmH^2(k_v,\widehat{G}^\ab(1))$ factors through the natural map $\rmH^2_\et(\mathscr{O}_{U,v},\widehat{\mathcal{G}}^\ab(1))\to \rmH^2(k_v,\widehat{G}^\ab(1))$. By smooth base change, see \emph{op.\ cit.}, Chapter VI, Corollary 4.2, this map coincides with the pullback $\rmH^2(\kappa_v,\widehat{G}^\ab(1))\to \rmH^2(k_v,\widehat{G}^\ab(1))$ along the projection $\Gamma_{k_v}\to \Gamma_{\kappa_v}$. Since the inertia subgroup $I_v$ acts trivially on $G$, then the short exact sequences \eqref{scindagepadique1} and \eqref{scindagepadique2} naturally assemble into the following adjunction diagram:
\[\begin{tikzcd}
	{\rmH^2(k_v,\widehat{G}^\ab(1))} & [-3em] \times & [-3em]{\rmH^1(k_v,G^\ab)} & {\Q/\Z} \\
	{\rmH^2(\kappa_v,\widehat{G}^\ab(1))} & \times & {\rmH^0(\kappa_v,G^\ab(-1))} & {\Q/\Z,}
	\arrow[from=1-3, to=1-4]
	\arrow[two heads, from=1-3, to=2-3]
	\arrow[equal, from=1-4, to=2-4]
	\arrow[hook, from=2-1, to=1-1]
	\arrow[from=2-3, to=2-4]
\end{tikzcd}\]
where the horizontal arrows are induced by the cup-product pairing. Since $k_v^\nr$ is a field of cohomological dimension $2$ containing $k$ and $\alpha\in \Sha^2_\nr(k,\widehat{G}^\ab(1))$, we see that the image of this class in $\rmH^2(k_v^\nr,\widehat{G}^\ab(1))\simeq\rmH^2(\kappa_v,\widehat{G}^\ab(1))$ must be orthogonal to $\rmH^0(k_v^\nr,G^\ab(-1))=\rmH^0(\kappa_v,G^\ab(-1))$ with respect to the cup-product pairing. By virtue of the above diagram, this is equivalent to the image of $\alpha$ in $\rmH^2(k_v,\widehat{G}^\ab(1))$ being orthogonal to the inverse image of $\mathcal{N}_v(G)$ in $\rmH^1(k_v,\widehat{G}^\ab)$. As $C$ is an irreducible curve, the complement of $U$ consists of finitely many points, hence the result.
\end{proof}

\begin{cor}\label{comparaisonsha2} Suppose in addition that $G$ is abelian. Then there is an equality
\[\Sha^2_\mathrm{nr}(k,\widehat{G}(1))=\Sha^2_\omega(k,\widehat{G}(1))\]
of subgroups of $\rmH^2(k,\widehat{G}(1))$.
\end{cor}
\begin{proof} Let $\alpha\in \Sha^2_\nr(k,\widehat{G}(1))$ and let $v\in C^{(1)}$ be a point for which the image of $\alpha$ in $\rmH^2(k_v,\widehat{G}(1))$ is orthogonal to the inverse image of $\mathcal{N}_v(G)$ in $\rmH^1(k_v,G)$. As explained in Remark \ref{remarquenonramificationalgébrique}, we see that $\mathcal{N}_v(G)$ coincides with the whole group $\rmH^0(\kappa_v,G(-1))$; therefore, the image of $\alpha$ is orthogonal to $\rmH^1(k_v,G)$. As the cup-product pairing $\rmH^2(k_v,\widehat{G}(1))\times \rmH^1(k_v,G)\to \Q/\Z$ is perfect, this implies that $\alpha$ dies in $\rmH^2(k_v,\widehat{G}(1))$, hence the first inclusion. Conversely, fix an embedding $G\hookrightarrow \SL_{n,k}$ for some $n\ge 1$ and denote $X:=\SL_{n,k}/G$. By virtue of \cite[Theorem 4.2]{hararischeidererszamuely} the image of $\Sha^2_\omega(k,\widehat{G}(1))$ along the composite map
\[\rmH^2(k,\widehat{G}(1))\longrightarrow  \rmH^3_\et(X,\Q/\Z(2))\longrightarrow \rmH^0(X,\mathscr{H}^3_\et(\Q/\Z(2)))\]
is contained in $\rmH^3_\nr(k(X),\Q/\Z(2))$. Applying Corollary~\ref{nonramificationalgébrique}, we infer that $\Sha^2_\omega(k,\widehat{G}(1))$ must be contained in $\Sha^2_\nr(k,\widehat{G}(1))$. Therefore, both groups are equal.
\end{proof}
\begin{rmk} More generally, if we assume that $G$ is a finite $k$-group such that the abelianisation map $G(\overline{k})\to G^\ab(\overline{k})$ admits a $\Gamma_k$-equivariant section, then the first part of the proof of Corollary \ref{comparaisonsha2} yields an inclusion $\Sha^2_\mathrm{nr}(k,\widehat{G}^\ab(1))\subset\Sha^2_\omega(k,\widehat{G}^\ab(1))$ (which we expect to be an equality). Indeed, if $\alpha\in \Sha^2_\nr(k,\widehat{G}^\ab(1))$ and $v\in C^{(1)}$ is a point for which the image of $\alpha$ in $\rmH^2(k_v,\widehat{G}^\ab(1))$ is orthogonal to the inverse image of $\mathcal{N}_v(G)$ in $\rmH^1(k_v,G^\ab)$, then pulling back the prescribed section in cohomology yields a section of the pushforward map $\rmH^0(\kappa_v,G(-1))\to \rmH^0(\kappa_v,G^\ab(-1))$. We readily see that $\mathcal{N}_v(G)$ coincides with the whole group $\rmH^0(\kappa_v,G^\ab(-1))$, and therefore the image of $\alpha$ is orthogonal to $\rmH^1(k_v,G^\ab)$. As the cup-product pairing $\rmH^2(k_v,\widehat{G}^\ab(1))\times \rmH^1(k_v,G^\ab)\to \Q/\Z$ is perfect, this implies that $\alpha$ dies in $\rmH^2(k_v,\widehat{G}^\ab(1))$ as desired.
\end{rmk}
\begin{thm}\label{contrexempleZ/8} Let $k=\Q_2(t)$ and $G=\Z/8$. Let $c\ge 3$, fix a classifying pair $(U\subset\A^m_k,c)$ for $G$, and denote $X:=U/G$. There exist infinitely many distinct classes in $\Sha^2_\mathrm{nr}(k,\widehat{G}(1))$ whose image along the natural map 
\[\Sha^2_\mathrm{nr}(k,\widehat{G}(1))\longrightarrow \rmH^{3,\norm}_\nr(k(X),\Q/\Z(2))\]
provide a reciprocity obstruction to weak approximation on $X$.
\end{thm}
\begin{proof} Let $c\ge 3$, fix a classifying pair $(U_0\subset\A^m_{\Q_2},c)$ for $G$ and a smooth compactification $X_0$ of $U_0/G$. By Wang's counter-example to Grunwald's theorem \cite{wang}, there exists a class $\alpha\in \Br_1(X_0)$ such that the evaluation map $\mathrm{ev}_\alpha:X_0(\Q_2)\to \Q/\Z$ is non-constant: the image of $\alpha$ in $\Br_1(X_0)/\Br(\Q_2)\hookrightarrow \rmH^1(\Q_2,\mmu_8)\simeq \Q_2^*/\Q_2^{*8}$ is given by the class of $16\in \Q_2^*$. Let $k=\Q_2(t)$ and denote $X:=X_0\otimes_{\Q_2} k$. We may view $\alpha$ as an element of $\Br(\Q_2(X_0))$ and let $\beta:=\alpha\smile [t]\in \rmH^3(k(X),\Q/\Z(2))$, where $[t]\in \rmH^1(k,\widehat{\Z}(1))$. This class $\beta$, interpreted as a cohomology class of the function field of the $\Q_2$-variety $X_0\times\P^1_{\Q_2}$, restricts to an unramified class on the generic fibre of the projection $X_0\times\P^1_{\Q_2}\to \P^1_{\Q_2}$. Indeed, $\alpha$ is unramified and $[t]$ comes from the base field $k$ (hence is trivially unramified). By construction, $\beta$ is also algebraic and non-zero: its residue with respect to $t=0$ is $\alpha$, which does not come from an element of $\Br(\Q_2)$. Moreover, it lifts to a non-zero cohomology class $\widetilde{\beta}\in F^2\rmH^3_\et((U_0\otimes_{\Q_2}k)/G,\Q/\Z(2))$ which comes from the element $[16]\smile [t]\in \rmH^2(k,\mmu_8^{\otimes2})$. Thanks to Proposition \ref{nonramificationalgébrique}, it follows that the latter actually lies in $\Sha^2_\mathrm{nr}(k,\mmu_8^{\otimes2})$. 

Let now $P_0,Q_0$ be two $\Q_2$-points of $X_0$ on which $\alpha$ takes distinct values. Denote $v_1:=0$ and $v_2:=\infty$ viewed as $\Q_2$-points of $\P^1_{\Q_2}$ and let $P$ be a $k_{v_1}$-point of $X$ which specialises into $P_0$, resp.\ let $Q$ be a $k_{v_2}$-point which specialises into $Q_0$. We then have 
\[\mathrm{inv}_{v_1}(\mathrm{ev}_\beta(P))+\mathrm{inv}_{v_2}(\mathrm{ev}_\beta(Q))\neq 0\in \Q/\Z,\] 
and $\mathrm{ev}_\beta(R)=0$ for any $k_v$-point $R$ of $X$ with $v$ a place distinct from $v_1$ and $v_2$. We thus see that $\beta$ provides a reciprocity obstruction to weak approximation in the sense of \cite{hararischeidererszamuely} on the $k$-variety $X$. By replacing $[t]$ with $[t-a]\in \rmH^1(k,\widehat{\Z}(1))$ for $a\in \Q_2$, we construct in a similar way classes $[16]\smile [t-a]\in \Sha^2_\nr(k,\mmu_8^{\otimes 2})$ whose image $\beta_a\in \rmH^3_\nr(k(X),\Q/\Z(2))$ is non-zero. As $\Q_2$ is uncountably infinite, it follows by construction that there cannot exist a single element $a_0\in \Q_2$ on which every $\beta_a$ has trivial residue. Hence, infinitely many $\beta_a$'s are distinct in $\rmH^3_\nr(k(X),\Q/\Z(2))$ and produce a reciprocity obstruction to weak approximation on $X$.
\end{proof}

\begin{rmk} In \cite[§3.4, Question 3]{LinhArithmetics}, the author asks whether the group $\Sha^2_\omega(k,\mmu_{2^m}^{\otimes2})$ is trivial for $m\ge 3$. By combining Corollary \ref{comparaisonsha2} with the example of Theorem \ref{contrexempleZ/8}, we see that the answer is negative for $m=3$. One can easily extend the above arguments to get a similar result when $m\ge 4$.
\end{rmk}

\begin{rmk} Following a suggestion of Federico Scavia, we can combine the example of Theorem \ref{contrexempleZ/8} with ideas of Ekedahl to construct a smooth and projective variety over $k=\Q_2(t)$ for which there exist infinitely many distinct degree three cohomology classes which are not of coniveau $\ge 1$. Indeed, if we let $c\ge3$ and $(U\subset \A^m_k,c)$ be a classifying pair for the group $G=\Z/8$, resp.\ if we let $G\hookrightarrow \SL_{n,k}$ be an embedding of $k$-groups for some $n\ge 1$, then we may find a smooth and projective $k$-variety $X$ which is stably birationally equivalent to both $U/G$ and $\SL_{n,k}/G$ and such that there is a commutative diagram:
\[\begin{tikzcd}
	{\rmH^3_\et(U/G,\Q/\Z(2))} & {\rmH^3(k(U/G),\Q/\Z(2))} \\
	{\rmH^3_\et(X,\Q/\Z(2))} & {\rmH^3(k(X),\Q/\Z(2)),}
	\arrow[from=1-1, to=1-2]
	\arrow["\wr"', from=1-1, to=2-1]
	\arrow[from=1-2, to=2-2]
	\arrow[from=2-1, to=2-2]
\end{tikzcd}\]
where the leftmost vertical map is an isomorphism; see \cite[§1.1]{ekedahlapproximation} for details. Since $U/G$ and $X$ are stably birationally equivalent then the righmost vertical arrow induces an isomorphism $\rmH^3_\nr(k(U/G),\Q/\Z(2))\simeq \rmH^3_\nr(k(X),\Q/\Z(2))$. On the other hand, if we denote by $N^\bullet\rmH^3_\et(X,\Q/\Z(2))$ the coniveau filtration on $\rmH^3_\et(X,\Q/\Z(2))$ as defined in \cite[§9.7]{brauerIII}, then there is a short exact sequence
\[0\longrightarrow N^1\rmH^3_\et(X,\Q/\Z(2))\longrightarrow \rmH^3_\et(X,\Q/\Z(2))\longrightarrow \rmH^3_\nr(k(X),\Q/\Z(2)).\]
Consequently, since the cohomology classes constructed in Theorem \ref{contrexempleZ/8} map to distinct non-zero unramified elements in $\rmH^3(k(U/G),\Q/\Z(2))$, then they yield infinitely many non-zero elements in the quotient $\rmH^3_\et(X,\Q/\Z(2))/N^1\rmH^3_\et(X,\Q/\Z(2))$. 

\end{rmk}

\bibliographystyle{amsalpha}
\bibliography{abcd}
\end{document}